\pdfoutput=1
\newif\ifpersonal
\newif\ifarxiv
\personaltrue 
\arxivtrue 
\RequirePackage[l2tabu, orthodox]{nag} 
\documentclass[10pt,a4paper,reqno,oneside]{amsart}
\usepackage{calligra}
\usepackage{amsmath,amsthm,amssymb,mathrsfs,mathtools,bm,eucal,tensor} 
\usepackage[scaled]{beramono,berasans}
\usepackage{enumerate,comment,braket,xspace,tikz-cd} 
\usepackage[all,cmtip]{xy} 
\usepackage[utf8]{inputenc} 
\usepackage[T1]{fontenc} 
\usepackage{lmodern}
\usepackage{scalerel}
\usepackage{varwidth}

\definecolor{linkcolor}{HTML}{005050}

\usepackage[centering,vscale=0.7,hscale=0.7]{geometry}
\usepackage{hyperref}
\usepackage[capitalize]{cleveref}
\usepackage{graphicx}
\usepackage{xparse}
\usepackage{url}
\usepackage[toc,page]{appendix}
\DeclareMathAlphabet{\mathcal}{OMS}{cmsy}{m}{n} 

\usepackage[backend=biber, style=alphabetic, sorting = nyt , giveninits=true,doi=false,url=false,maxbibnames=99]{biblatex}
\usepackage{vmargin}
\setpapersize{A4}
\setmarginsrb{25mm}{10mm}{25mm}{10mm}%
{12mm}{10mm}{5mm}{10mm}

\theoremstyle{plain}
\newtheorem{thm-intro}{Theorem}
\newtheorem{thm}{Theorem}[section]

\newtheorem*{thm*}{Theorem}
\newtheorem{lem}[thm]{Lemma}

\newtheorem*{lem*}{Lemma}
\newtheorem{prop}[thm]{Proposition}
\newtheorem{proposition}[thm]{Proposition}

\newtheorem{conj}[thm]{Conjecture}
\newtheorem{Question}{Question}

\newtheorem{cor}[thm]{Corollary}

\newtheorem{cor-intro}[thm-intro]{Corollary}

\theoremstyle{definition}
\newtheorem{definition}[thm]{Definition}
\newtheorem{defin}[thm]{Definition}

\newtheorem{defin-intro}[thm-intro]{Definition}

\theoremstyle{remark}
\newtheorem{rem}{Remark}

\newtheorem{example}[thm]{Example}

\newtheorem{eg-intro}[thm-intro]{Example}
\newtheorem{rem-intro}[thm-intro]{Remark}
\numberwithin{equation}{section}

\ifpersonal
\newcommand*{\personal}[1]{\textcolor[rgb]{0.6,0.6,1}{(Personal: #1)}}
\newcommand*{\todo}[1]{\textcolor{red}{(Todo: #1)}}
\else
\newcommand*{\personal}[1]{\ignorespaces}
\newcommand*{\todo}[1]{\ignorespaces}
\fi

\newcommand{\Tot}{\textbf{Tot}}

\newcommand{\Biggg}{\bBigg@{3.5}}
\newcommand{\C}{\mathbb C}
\newcommand{\CP}{\mathbb{CP}}

\newcommand{\Q}{\mathbb Q}
\newcommand{\R}{\mathbb R}
\newcommand{\Z}{\mathbb Z}
\newcommand{\N}{\mathbb N}

\newcommand{\bbk}{\mathbb K}

\newcommand{\sw}{\mathsf w}
\newcommand{\sE}{\mathsf E}

\newcommand{\cF}{\mathcal F}

\newcommand{\cO}{\mathcal O}

\newcommand{\cY}{\mathcal Y}

\DeclareFontFamily{U}{BOONDOX-calo}{\skewchar\font=45 }
\DeclareFontShape{U}{BOONDOX-calo}{m}{n}{<-> s*[1.05] BOONDOX-r-calo}{}
\DeclareFontShape{U}{BOONDOX-calo}{b}{n}{<-> s*[1.05] BOONDOX-b-calo}{}
\DeclareMathAlphabet{\mathcalboondox}{U}{BOONDOX-calo}{m}{n}

\newcommand{\bbA}{\mathbb A}

\newcommand{\bbP}{\mathbb P}

\newcommand{\bbH}{\mathbb H}

\newcommand{\Gr}{\text{Gr}}

\makeatletter
\let\save@mathaccent\mathaccent
\newcommand*\if@single[3]{%
	\setbox0\hbox{${\mathaccent"0362{#1}}^H$}%
	\setbox2\hbox{${\mathaccent"0362{\kern0pt#1}}^H$}%
	\ifdim\ht0=\ht2 #3\else #2\fi
}
\newcommand*\rel@kern[1]{\kern#1\dimexpr\macc@kerna}
\newcommand*\widebar[1]{\@ifnextchar^{{\wide@bar{#1}{0}}}{\wide@bar{#1}{1}}}
\newcommand*\wide@bar[2]{\if@single{#1}{\wide@bar@{#1}{#2}{1}}{\wide@bar@{#1}{#2}{2}}}
\newcommand*\wide@bar@[3]{%
	\begingroup
	\def\mathaccent##1##2{%
		\let\mathaccent\save@mathaccent
		\if#32 \let\macc@nucleus\first@char \fi
		\setbox\z@\hbox{$\macc@style{\macc@nucleus}_{}$}%
		\setbox\tw@\hbox{$\macc@style{\macc@nucleus}{}_{}$}%
		\dimen@\wd\tw@
		\advance\dimen@-\wd\z@
		\divide\dimen@ 3
		\@tempdima\wd\tw@
		\advance\@tempdima-\scriptspace
		\divide\@tempdima 10
		\advance\dimen@-\@tempdima
		\ifdim\dimen@>\z@ \dimen@0pt\fi
		\rel@kern{0.6}\kern-\dimen@
		\if#31
		\overline{\rel@kern{-0.6}\kern\dimen@\macc@nucleus\rel@kern{0.4}\kern\dimen@}%
		\advance\dimen@0.4\dimexpr\macc@kerna
		\let\final@kern#2%
		\ifdim\dimen@<\z@ \let\final@kern1\fi
		\if\final@kern1 \kern-\dimen@\fi
		\else
		\overline{\rel@kern{-0.6}\kern\dimen@#1}%
		\fi
	}%
	\macc@depth\@ne
	\let\math@bgroup\@empty \let\math@egroup\macc@set@skewchar
	\mathsurround\z@ \frozen@everymath{\mathgroup\macc@group\relax}%
	\macc@set@skewchar\relax
	\let\mathaccentV\macc@nested@a
	\if#31
	\macc@nested@a\relax111{#1}%
	\else
	\def\gobble@till@marker##1\endmarker{}%
	\futurelet\first@char\gobble@till@marker#1\endmarker
	\ifcat\noexpand\first@char A\else
	\def\first@char{}%
	\fi
	\macc@nested@a\relax111{\first@char}%
	\fi
	\endgroup
}
\makeatother

\newcommand{\Coh}{\mathrm{Coh}}

\newcommand{\Aff}{\mathrm{Aff}}

\newcommand{\Top}{\mathcal T\mathrm{op}}

\newcommand{\Fuk}{\mathrm{Fuk}}
\newcommand{\FS}{\mathrm{FS}}

\newcommand{\vol}{\mathsf{vol}}

\usetikzlibrary{decorations.markings} 
\tikzset{
  closed/.style = {decoration = {markings, mark = at position 0.5 with { \node[transform shape, xscale = .8, yscale=.4] {/}; } }, postaction = {decorate} },
  open/.style = {decoration = {markings, mark = at position 0.5 with { \node[transform shape, scale = .7] {$\circ$}; } }, postaction = {decorate} }
}

\DeclareMathOperator{\Bl}{Bl}

\DeclareMathOperator{\Coker}{Coker}

\DeclareMathOperator{\cone}{cone}

\DeclareMathOperator{\crit}{crit}

\DeclareMathOperator{\Hom}{Hom}
\DeclareMathOperator{\Image}{Im}

\DeclareMathOperator{\Ker}{Ker}

\DeclareMathOperator{\Spec}{Spec}
\DeclareMathOperator{\Span}{Span}

\DeclareMathOperator{\supp}{supp}

\title{Mirror P=W Conjecture and Extended Fano/Landau-Ginzburg Correspondence}

\author{Sukjoo Lee}

\address{Department of Mathematics, University of Edinburgh, EH9 3FD, UK}
\email{Sukjoo.Lee@ed.ac.uk}

\begin{document}

\maketitle
\begin{abstract}
The mirror P=W conjecture, formulated by Harder-Katzarkov-Przyjalkowski \cite{PWmirror}, predicts a correspondence between weight and perverse filtrations in the context of mirror symmetry. In this paper, we revisit this conjecture through the lens of mirror symmetry for a Fano pair $(X,D)$, where $X$ is a smooth Fano variety and $D$ is a simple normal crossing divisor. We introduce its mirror object as a multi-potential analogue of a Landau-Ginzburg (LG) model, which we call the hybrid LG model. This model is expected to capture the mirrors of all irreducible components of $D$. We study the topological aspects, particularly the perverse filtration, and the Hodge theory of hybrid LG models, building upon the work of Katzarkov-Kontsevich-Pantev \cite{KKPbogomolov}. As an application, we discover an interesting upper bound on the multiplicativity of the perverse filtration for a hybrid LG model. Additionally, we propose a relative version of the homological mirror symmetry conjecture and explain how the mirror P=W conjecture naturally emerges from it.

\end{abstract}
\tableofcontents

\section{Introduction}

\subsection{Motivation: the mirror P=W conjecture} 
Traditionally mirror symmetry relates two compact Kähler $n$-dimensional Calabi-Yau manifolds $X$ and $X^\vee$: the complex (algebraic) geometry of $X$ (B-side) is equivalent to the symplectic geometry of $X^\vee$ (A-side) and vice versa \cite{MirrorBook}\cite{MirrorBook2Clay}. The categorial formulation of mirror symmetry known as the homological mirror symmetry (HMS for short) conjecture was proposed by M. Kontsevich \cite{Kontsevichmirror}. According to this conjecture, the derived category of coherent sheaves on $X$, ${D^b\Coh}(X)$, is equivalent to the (derived) Fukaya category of $X^\vee$, $\mathrm{Fuk}(X^\vee)$ \cite{SeidelFukayaLeftbook}\cite{FOOO}. On the other hand, the most basic form of mirror symmetry is the symmetry of Hodge numbers: $h^{p,q}(X)=h^{n-p,q}(X^\vee)$ for all $p,q \geq 0$. This simple relation becomes more involved when we attempt to extend mirror symmetry to non-compact Kähler Calabi-Yau's $U$ and $U^\vee$. In this case, Hodge-theoretic data are refined by two filtrations: Deligne's canonical weight filtration $W_\bullet$ and the perverse Leray filtration $P_\bullet$ associated to the canonical affinization map. Incorporating these refinements, we define a perverse-mixed Hodge polynomial of $U$ as 
\[
PW_U(u,t,w,p):=\sum_{a,b,r,s}(\dim_\C \Gr_F^a\Gr_{s+b}^W \Gr^P_{s+r}(H^s(U,\C)))u^at^sw^bp^r
\]
where $F^\bullet$ is the Hodge filtration (Definition \ref{def:PMH poly}).
\begin{conj}\label{conj:mirrorP=W}(Mirror P=W Conjecture) 
     Assume that two $n$-dimensional log Calabi-Yau varieties $U$ and $U^\vee$ are mirror to each other. Then we have the following polynomial identity
     \begin{equation}\label{eq:mirrorP=W}
          PW_U(u^{-1}t^{-2}, t,p,w)u^nt^n=PW_{U^\vee}(u,t,w,p)
     \end{equation}
 \end{conj}   
 
  In this paper, we will focus on the case where the perverse filtration on $U$ becomes trivial. In this case, the relation (\ref{eq:mirrorP=W}) becomes simpler so that we can see the relation between the weight filtration on $H^\bullet(U)$ and the perverse Leray filtration on $H^\bullet(U^\vee)$ more directly. By taking $w=1$ and changing the variables, the equation (\ref{eq:mirrorP=W}) becomes
    \begin{equation}\label{eq:mirror P=W oneside} \dim_\C\Gr_F^q\Gr^W_{p+q+r}H^{p+q}(U)=\dim_\C\Gr_F^{n-q}\Gr^P_{n+p-q+r}H^{n+p-q}(U^\vee)
    \end{equation}
    for $p,q$ and $r$. To prove the identities (\ref{eq:mirror P=W oneside}), it is natural to try and realize each filtration as a limiting filtration on the abutment of a suitable spectral sequence. This requires a choice of extra data: a compactification $(X,D)$ of $U$ for the weight filtration and a morphism to the affine space $h:U^\vee \to \C^N$ for the perverse filtration. Therefore, the natural question we should answer is

    \begin{Question}\label{Question A}
        How do we describe the relation between the choices of $(X,D)$ and $h:U^\vee \to \C^N$ in the context of mirror symmetry?
    \end{Question}
        
       Of course, even if we answer Question \ref{Question A}, it is not immediately clear why the two spectral sequences should be isomorphic, particularly because they have very different natures. Indeed, the spectral sequence for the weight filtration is described in terms of the combinatorics of the compactification $(X,D)$ (a Čech-type sequence) while the perverse filtration is described in terms of a general flag (see Section \ref{sec:mirrorPW}). This leads to the following question:
        
        \begin{Question}\label{Question B}
             Is there a combinatorial (Čech-type) description of the perverse Leray filtration? How does it appear in the context of mirror symmetry?
        \end{Question}
         We will focus on the case when $U$ admits a compactification $(X,D)$ where $X$ is a smooth Fano manifold and $D$ is an anti-canonical divisor. Mirror symmetry of such a Fano pair has been studied extensively under the \textit{Fano/Landau-Ginzburg (Fano/LG) correspondence}\cite{KKPHodgemirror}\cite{KKPbogomolov}.\footnote{In this paper, we consider a Fano pair on the B-side and its dual LG model on the A-side.} The mirror dual of the pair $(X,D)$ is a Landau-Ginzburg (LG) model, a pair $(Y,\omega, \mathsf{w}:Y \to \C)$ where $Y$ is an $n$-dimensional Kähler Calabi-Yau manifold with the Kähler form $\omega$, and $\mathsf{w}$ is a holomorphic map that is a locally trivial fibration with smooth fibers near infinity.\footnote{Here, $Y$ is the same as $U^\vee$. We use $Y$ to emphasize that we are studying the mirror symmetry of Fano pairs.} The map $\sw:Y \to \C$ is called a \textit{LG potential} and its generic fiber, denoted by $Y_{sm}$, is expected to be mirror to the anti-canonical divisor $D$. One of the key features of the Fano/LG correspondence is the functoriality between mirror pairs $(X,\sw:Y \to \C)$ and $(D, Y_{sm})$, which would capture the mirror relation between $U$ and $Y$. These mirror relations are summarized in Table 1
         \begin{center}
\begin{tabular}{ |c|c|c|c| } 
 \hline
 B-side & $X$ & $D$ & $U$ \\ 
 \hline
 A-side & $\mathsf{w}:Y \to \C$ & $ Y_{sm}$ & $Y$ \\ 
 \hline
\end{tabular}
\end{center}
\begin{center}
    (Table 1. Fano/LG correspondence)
\end{center}
\noindent 
         
        On the cohomology level, the functoriality can be described as morphisms between relevant cohomology groups. For example, on the B-side, there is Gysin morphism $\iota_!:H^*(D) \to H^*(X)$ where $\iota:D \hookrightarrow X$ is the canonical inclusion. The corresponding morphism on the A-side is  expected to be the connecting homomorphism  $\rho:H^*(Y_{sm}) \to H^{*+1}(Y, Y_{sm})$. When $D$ is smooth, then Gysin morphism $\iota_!$ becomes the differential of the $E_1$-page of the spectral sequence for the weight filtration on $H^*(U)$. Similarly, when $Y_{sm}$ is proper, the connecting homomorphism $\rho$ becomes the differential of the $E_1$-page of the spectral sequence for the perverse Leray filtration on $H^*(Y)$. Note that the two cases occur at the same time since the mirror symmetry predicts that the smoothness of $D$ is reflected in the properness of the LG potential $\sw:Y \to \C$. Therefore, in this case, one can see that the mirror P=W conjecture follows from the functoriality in the Fano/LG correspondence. The precise relationship between the two morphisms $\iota_!$ and $\rho$ as well as how it implies the mirror P=W conjecture are discussed in Section \ref{sec:Fano/LG}. 
        
        However, when $D$ has more than one component, such a formulation of the functoriality loses information about the mirror symmetry of each irreducible component of $D$, and hence will not be enough to produce the mirror P=W correspondence. Motivated by the SYZ construction \cite{SYZ}, one can capture the mirror object of each irreducible component of $D$ by considering a multi-potential, each coordinate of which corresponds to counting holomorphic disks touching each irreducible component of $D$ \cite{AurouxFanoSYZ}\cite[Section 5.3]{AurouxFanoSYZLag}. 
        
    \subsection{Hybrid LG model}
    In the introduction, we mainly consider the case where $D$ has two irreducible components because this case is sophisticated enough to deliver the ideas of a multi-potential and turns out to be a building block for the general case. Let $X$ be a smooth Fano manifold and $D$ be a simple normal crossing anti-canonical divisor $D=D_1 \cup D_2$ with $D_{12}=D_1 \cap D_2$ with $(D_1, D_{12})$ and $(D_2,D_{12})$ being again smooth Fano pairs. In this case, the SYZ mirror construction yields a pair of potentials $h=(h_1, h_2):(Y,\omega) \to \C^2$ where $h_i$ counts disks touching each irreducible component $D_i$. It is furthermore expected that the mirror of the Fano pair $(D_1, D_{12})$ (resp. $(D_2, D_{12})$) is an ordinary LG model $(Y_1, \sw_1:=h_2|_{Y_1})$ (resp. $Y_2, \sw_2:=h_1|_{Y_2}$) where $Y_i$ is a generic fiber of $h_i$ for $i=1,2$. We denote a generic fiber of $h$ by $Y_{12}$ which is mirror to $D_{12}$. One can also obtain an ordinary LG potential $\mathsf{w}:=h_1+h_2:Y \to \C$, which is now non-proper, by composing with the sum map $\Sigma:\C^2 \to \C$. This is expected to be mirror to the Fano pair $(X,D)$ in the classical sense of Fano/LG correspondence introduced before. These mirror relations are summarized in Table 2.\\
    
    \begin{center}
\begin{tabular}{ |c|c|c|c|c| } 
 \hline
 B-side & $(X,D)$ & $(D_1, D_{12})$ & $(D_2, D_{12})$ & $D_{12}$ \\ 
 \hline
 A-side & $(\mathsf{w}:Y \to \C, Y_{sm})$ & $(\sw_1:Y_1 \to \C, Y_{12})$ & $(\sw_2:Y_2 \to \C, Y_{12})$ & $Y_{12}$ \\ 
 \hline
\end{tabular}
\end{center}
\begin{center}\label{Table2}
    (Table 2. Extended Fano/LG correspondence)
\end{center}

\noindent To guarantee that the A-side objects in Table 2 are ordinary LG models, one should put some geometric restrictions on the multi-potential $h=(h_1,h_2):Y \to \C^2$.\footnote{The reference to the $A$-side is somewhat misleading, because we mainly discuss topological and algebro-geometric aspects, leaving the symplectic geometric aspects in the speculative discussion in Section 5.}

\begin{defin}(Definition \ref{def: weak hybrid LG model general})
    A \textbf{hybrid Landau-Ginzburg (LG) model of rank 2} is a triple $(Y,\omega, h=(h_1, h_2):Y \to \C^2)$ where
    \begin{enumerate}
        \item $(Y,\omega)$ is a $n$-dimensional complex Kähler Calabi-Yau manifold with a Kähler form $\omega \in \Omega^2(Y)$;
        \item $h:Y \to \C^2$ is a proper (surjective) holomorphic map that satisfies the following compatible local trivializations: 
        \begin{enumerate}
            \item $h:Y \to \C^2$ is a locally trivial fibration near infinity corner $(\infty, \infty)$ with smooth fibers;
            \item for $\{i,j\}=\{1,2\}$, the map $h_i:Y \to \C$ is a locally trivial fibration near infinity with smooth fibers which also preserves the induced map $h_j$.
        \end{enumerate}
    \end{enumerate}
\end{defin}

Note that the condition $(2)-(b)$ implies that the map $h_i:Y \to \C$ is a locally trivial fibration of ordinary LG models $(Y_i, \sw_i)$ near infinity. Furthermore, the induced map $\sw:=h_1+h_2:Y \to \C$ indeed becomes a locally trivial fibration near infinity with smooth fibers. This follows from what we call the \textit{gluing property}:

\begin{prop}\label{prop:Gluing property intro}(Proposition \ref{prop:Gluing Property in general})
    Let $(Y, \omega, h:Y \to \C^2)$ be a hybrid LG model and $H$ be a generic line in the base $\C^2$, which is not parallel to any coordinate lines. There exists an open cover $H=U_1 \cup U_2$ such that for $i=1,2$, the induced map $h^{-1}(U_i) \to U_i$ is  isotopic to the induced LG potential $\sw_i:Y_i \to \C$.
\end{prop}

\noindent Heuristically, Proposition \ref{prop:Gluing property intro} implies that the induced fibration $h_{Y_{sm}}:=h|_{Y_{sm}}:Y_{sm} \to \C$ is glued from two ordinary LG models $(Y_1, \sw_1)$ and $(Y_2,\sw_2)$. This is the mirror counterpart of the structure that the anti-canonical divisor $D$ is the union of the pairs $(D_1, D_{12})$ and $(D_2, D_{12})$. In Section \ref{sec:general construction}, when $X$ is a smooth Fano complete intersection, we construct a mirror hybrid LG model by applying the Hori-Vafa construction and taking a suitable fiberwise compactification introduced in \cite{Przyjalalgorithm}.

On the cohomology level, the gluing property yields an isomorphism of cohomology groups
\begin{equation*}
    H^*(Y_{sm}, Y_{12}) \cong H^*(Y_1, Y_{12}) \oplus H^*(Y_2, Y_{12})
\end{equation*} 
For $i=1,2$, this provides a well-defined morphism $\rho^i_{Y}:H^{*-1}(Y_i, Y_{12}) \to H^*(Y, Y_{sm})$ by composing with the connecting homomorphism of the long exact sequence of cohomology groups associated to the triple $(Y_{12}, Y_{sm}, Y)$. Therefore, we have a commutative diagram of cohomology groups
\begin{equation}\label{eq:A-side rk2 intro cohomology}
        \begin{tikzcd}
    H^{*-2}(Y_{12})\arrow[r, "\rho^{12}_1"] \arrow[d, "\rho^{12}_2"] & H^{*-1}(Y_1, Y_{12}) \arrow[d, "\rho^1_{Y}"] \\
    H^{*-1}(Y_2, Y_{12})\arrow[r, "\rho^2_{Y}"] & H^{*}(Y, Y_{sm})
    \end{tikzcd}
    \end{equation}
where $\rho^{12}_i$ is the connecting homomorphism in the long exact sequence of cohomology groups of the pair $(Y_{12}, Y_i)$. This diagram also respects the monodromy operations associated to the hybrid LG model. Consider the monodromy $T_1$ and $T_2$ of a generic fiber $h^{-1}(t_1, t_2)$ of $h:Y \to \C^2$ that correspond to a loop around infinity $(e^{2\pi it}\cdot t_1, t_2)$ and $(t_1, e^{2\pi it}\cdot  t_2)$ for $(0 \leq t \leq 1)$, respectively. They induce the monodromy actions on the diagram (\ref{eq:A-side rk2 intro cohomology}), still denoted the same. We also write the composition as $T=T_1 \circ T_2$. If the hybrid LG model $(Y, \omega, h:Y \to \C^2)$ is mirror to the Fano pair $(X,D)$, then it is predicted that all such monodromy actions are indeed unipotent. This allows the definition of the monodromy weight filtration $W_\bullet^T$. We refer to such a hybrid LG model of this kind as being of \textit{Fano type} (see Definition \ref{def: Fanotypeintro} for the precise definition).

\begin{thm}\label{thm:intro a-side}(Theorem \ref{thm: A-side cubical diagram coh}\ref{thm:mainhodge})
	    Let $(Y,\omega, h:Y \to \C^2)$ be a hybrid LG model introduced above. Then for $-n \leq a \leq n$, there is a well-defined diagram of cohomology groups 
	    \begin{equation*}
        \begin{tikzcd}
    H^{a+n-2}(Y_{12})\arrow[r, "\rho^{12}_1"] \arrow[d, "\rho^{12}_2"] & H^{a+n-1}(Y_1, Y_{12}) \arrow[d, "\rho^1_{Y}"] \\
    H^{a+n-1}(Y_2, Y_{12})\arrow[r, "\rho^2_{Y}"] & H^{a+n}(Y, Y_{sm})
    \end{tikzcd}
    \end{equation*}
	    which is compatible with the monodromy actions $T_1$, $T_2$, and $T$. By applying Mayer-Vietoris sign rule (\ref{eq:MV sign rule}), it gives rise to the $E_1$-page of the spectral sequence for the perverse Leray filtration associated to $h:Y \to \C^2$ on $H^*(Y)$. In particular, if the hybrid LG model is of Fano type, then the filtration $W_\bullet^T$ on the $E_1$-page  is compatible with Deligne's canonical weight filtration $W_\bullet$ on $H^*(Y)$.
	\end{thm}
The proof is divided into two steps: First, we describe the combinatorial nature of the perverse Leray filtration. Second, we construct the relevant limiting mixed Hodge structures to elevate the statement to the Hodge-theoretic level, a topic we will discuss in the next section. For the first part, recall that the gluing property tells us that there is no serious topological distinction between the fibrations on $Y_{sm}$ and $Y_1 \sqcup Y_2$ where $Y_1 \sqcup Y_2$ is the normal crossing union of $Y_1$ and $Y_2$. It suggests to consider a flag of subvarieties $(Y_{12} \subset Y_1 \sqcup Y_2)$ where $Y_1 \sqcup Y_2$ is the simple normal crossing union of $Y_1$ and $Y_2$ over $Y_{12}$. Then the diagram above produces the $E_1$-page of the spectral sequence for the flag filtration given by $(Y_{12} \subset Y_1 \sqcup Y_2)$. This provides an answer for Question \ref{Question B}.
    \begin{prop}\label{prop:comb perv intro}(Proposition \ref{prop:two-perverse-same})
            The flag filtration associated to $(Y_{12} \subset Y_1 \sqcup Y_2)$ is the same as the perverse Leray filtration associated to a hybrid LG potential $h:Y \to \C^2$.
    \end{prop}
While the gluing property is the key to obtain combinatorial description of the perverse Leray filtration for a hybrid LG model, there is an algebraic way to describe it once the genericity assumption on $Y_i$ is further imposed. We treat this case in Section \ref{sec:comb perv}. 
   
\subsection{Hodge theory of hybrid LG models}  
We start with extending the work of Katzarkov-Kontsevich-Pantev \cite{KKPbogomolov} in the hybrid setting. Again, we mainly discuss a hybrid LG model $(Y, h:Y \to \C^2)$ with two potentials for the purpose of the introduction. The key ingredient is a choice of an appropriate compactification $((Z, D_Z), f=(f_1,f_2):Z \to \bbP^1 \times \bbP^1)$ of a given hybrid LG model $(Y, h:Y \to \C^2)$, where $Z$ is smooth and $D_Z$ is a simple normal crossing boundary divisor. A non-trivial condition we impose on the compactification is that the morphism $f$ is semi-stable at $(\infty, \infty) \in \bbP^1 \times \bbP^1$, which is allowable because the local triviality condition of $h$ near infinity corner $(\infty, \infty)$. We also call the compactified LG model \textit{tame} if the pole divisor  $f_i^{-1}(\infty)$ is reduced for $i=1,2$.

\begin{defin}\label{def: Fanotypeintro}
    A hybrid LG model $(Y, h:Y \to \C^2)$ is of \textbf{Fano type} if it admits a tame compactification.
\end{defin}

Such a compactification provides two operators on the logarithmic de Rham complex $(\Omega_Z^\bullet(\log D_Z), d)$: the wedge product with $df_1$ and the wedge product with $df_2$. This allows us to consider the subcomplex, denoted by $(\Omega^\bullet_Z(\log D_Z, f_1, f_2),d)$, that is preserved by both $\wedge df_1$ and $\wedge df_2$. Not only that, by taking the pushout along the canonical inclusions $j_i:(\Omega^\bullet_Z(\log D_Z, f_1, f_2),d) \hookrightarrow (\Omega^\bullet_Z(\log D_Z, f_i),d)$ for $i=1,2$ as illustrated in the diagram (\ref{diag:pushout intro 2 comp}), we obtain another subcomplex, denoted by $(\Omega^\bullet_Z(\log D_Z, f_{12}),d)$. This subcomplex is generated by forms that are preserved by either $\wedge df_1$ or $\wedge df_2$. In other words, we have the pullback and pushout diagram of complexes as follows:

    \begin{equation}\label{diag:pushout intro 2 comp}
	\begin{tikzcd}
	\Omega^\bullet_Z(\log D_Z, f_1, f_2) \arrow[r, hook, "j_1"] \arrow[d, hook, "j_2"] & \Omega^\bullet_Z(\log D_Z, f_1) \arrow[d, "i_1"] \\
	\Omega^\bullet_Z(\log D_Z, f_2) \arrow[r, "i_2"] & \Omega^\bullet_Z(\log D_Z,f_{12}) 
	\end{tikzcd}
	\end{equation}
	
\noindent We will show that the Hodge-to-de Rham spectral sequence for both subcomplexes degenerates at the $E_1$-page (Proposition \ref{prop: degeneracy of Hodge filtration general case}), allowing us to define two different sets of LG Hodge numbers.
    \begin{defin}\label{def:f-adapted Hodge intro}
	      Let $((Z,D_Z),f:Z \to \bbP^1 \times \bbP^1)$ be a compactified hybrid LG model of $(Y, h:Y \to \C^2)$. We define two LG Hodge numbers as follows: For $p,q \geq 0$,
	      \begin{equation*}
	      \begin{aligned}
	        & f^{p,q}(Y,h):=\dim_\C \bbH^q(Z, \Omega_Z^p(\log D_Z, f_{12})) \\
	        & f^{p,q}(Y,h_1,h_2):=\dim_\C \bbH^q(Z, \Omega_Z^p(\log D_Z, f_1, f_2))
	      \end{aligned}
	      \end{equation*}
	 \end{defin}

From now on, we assume that the hybrid LG model is of Fano type. In this case, we will show that the cohomology $\bbH^k(Z,\Omega^\bullet_Z(\log D_Z, f_1, f_2))$ (resp. $\bbH^k(Z,\Omega^\bullet_Z(\log D_Z, f_{12}))$) is the Gauss-Manin limit of the relative cohomology group $H^k(Y, Y_1 \sqcup Y_2)$ (resp. $H^k(Y, Y_{12})$) (see Proposition \ref{prop: limit of fadpated complex general case}). By following \cite{ShamotoHodgeTate}, we also construct a limiting mixed Hodge structure on the relative cohomolgoy groups, denoted by $(H^k(Y, Y_{1, \infty}\sqcup Y_{2, \infty}), W_\bullet^{h_1, h_2}, F^\bullet_{h_1, h_2})$ and $(H^k(Y, Y_{12, \infty}), W_\bullet^h, F^\bullet_h)$, so that the above limits are indeed compatible with Hodge filtrations. Now the question is the relation between the set of LG Hodge numbers in Definition \ref{def:f-adapted Hodge intro} and another set of LG Hodge numbers, the complex dimension of the associated graded pieces of the monodromy weight filtrations.

    \begin{conj}(Extended KKP Conjecture)(Conjecture \ref{conj:Extended KKP})\label{conj:Extended KKP-intro}
    \begin{enumerate}
        \item Let $(Y, h:Y \to \C^2)$ be a hybrid LG model of Fano type. There are identifications of the LG Hodge numbers:
        \begin{equation*}
        \begin{aligned}
            & f^{p,q}(Y,h)= \dim_\C \Gr^{W^h}_{2p}H^{p+q}(Y, Y_{12, \infty};\C)\\
            & f^{p,q}(Y,h_1, h_2)=\dim_\C \Gr^{W^{h_1, h_2}}_{2p}H^{p+q}(Y, Y_{1, \infty}\cup Y_{2,\infty};\C)
        \end{aligned}
        \end{equation*}
        for non-negative integers $p,q \geq 0$.
        \item Let $(Y, \sw:Y \to \C)$ be the associated ordinary LG model of the hybrid LG model $(Y, h:Y \to \C^2)$. Then there are identifications of the LG Hodge numbers:
        \begin{equation*}
            f^{p,q}(Y, \sw)=f^{p,q}(Y, h_1, h_2)
        \end{equation*}
        for non-negative integers $p,q \geq 0$. 
    \end{enumerate}
    \end{conj}
    The second part of Conjecture \ref{conj:Extended KKP-intro} explains in what sense it extends the original KKP conjecture (Conjecture \ref{conj:KKP Hodge numbers smooth}). Meanwhile, the first part of Conjecture \ref{conj:Extended KKP-intro} is motivated by a relative version of homological mirror symmetry for a Fano pair (see Section \ref{sec:Extended Fano/LG}). We provide a proof for Conjecture \ref{conj:Extended KKP-intro} under certain assumptions inspired by mirror symmetry.
    
    \begin{thm}
        Let $(Y, h:Y \to \C^2)$ be a hybrid LG model of Fano type. 
        
        \begin{enumerate}
            \item (Theorem \ref{prop: Hodge-Tate implies KKP hybrid setting:General} If Deligne's canonical mixed Hodge structures on $H^k(Y)$ and $H^k(Y_i)$ are Hodge-Tate for $i=1,2$, $k \geq 0$, then the first part of Conjecture \ref{conj:Extended KKP-intro} holds. 
            \item (Theorem \ref{thm:extended KKP second}) If we further assume that Deligne's canonical mixed Hodge structures on $H^k(Y, Y_{sm})$ is Hodge Tate for $k \geq 0$, then the second part of Conjecture \ref{conj:Extended KKP-intro} holds as well.
        \end{enumerate}
    \end{thm}

    We conclude this section by introducing a de Rham theoretic description of the perverse filtration and its application. Consider the filtration on $\Omega_Z^\bullet(\log D_Z)$ induced by the inclusions
\begin{equation*}
    \Omega_Z^\bullet(\log D_Z, f_1, f_2) \hookrightarrow \Omega_Z^\bullet(\log D_Z, f_{12}) \hookrightarrow \Omega_Z^\bullet(\log D_Z)
\end{equation*}
and define an increasing filtration $G^f_\bullet\bbH^k(Z, \Omega_Z^\bullet(\log D_Z))$ by 
\begin{equation*}
        \begin{aligned}
            &G^f_kH^k(Y;\C):=\Image(\bbH^k(Z, \Omega_Z^\bullet(\log D_Z, f_1, f_2) \to \bbH^k(Z,\Omega_Z^\bullet(\log D_Z)) \\
            &G^f_{k+1}H^k(Y;\C):=\Image(\bbH^k(Z, \Omega_Z^\bullet(\log D_Z, f_{12}) \to \bbH^k(Z,\Omega_Z^\bullet(\log D_Z))
        \end{aligned}
        \end{equation*}
    Combining this with Proposition \ref{prop:comb perv intro}, we have the de Rham description of the perverse Leray filtration.
\begin{thm}(Theorem \ref{thm: two fil same general})
    For $k \geq 0$, the two filtrations $P^h_\bullet$ and $G^f_\bullet$ on $H^k(Y;\C)$ are the same. 
\end{thm}

All the arguments above are extended to the general case. A compactified hybrid LG model is given by $f=(f_1, \dots, f_N):Z \to (\bbP^1)^N$ where $N$ is the number of irreducible components of the anti-canonical divisor $D$. Especially regarding the perverse filtration, the logarithmic de Rham complex $(\Omega_Z^\bullet(\log D_Z), d)$ is filtered by the number of differentials $df_i \wedge$ the forms are preserved by. As an application, we have the following statement about the weak multiplicativity of the perverse filtration. 
    \begin{thm}(Theorem \ref{thm:multiplicativity})
    For $k_1, k_2 \geq 0$, the perverse Leray filtration $P_\bullet$ associated to a hybrid LG model $(Y, h:Y \to \C^N)$ satisfies the following relation under the cup product
    \begin{equation*}
        \cup: P_{k_1+a_1}H^{k_1}(Y;\C) \otimes P_{k_2+a_2}H^{k_2}(Y;\C) \to P_{k_1+k_2+\min(a_1, a_2)}H^{k_1+k_2}(Y;\C)
    \end{equation*}
    where $0 \leq a_1, a_2 \leq N$.
    \end{thm}

\subsection{Extended Fano/LG correspondence}
In Section \ref{sec:Extended Fano/LG}, we will discuss the extended Fano/LG correspondence more, providing specific formulations of the relative version of mirror symmetry. Particularly, we will focus on $B$-to-$A$ mirror symmetry, where the Fano pair resides on the $B$-side and the LG model on the $A$-side. We should note that this section is speculative, as significant work is required to develop a precise treatment for the symplectic geometry of hybrid LG models. The primary aim here is to propose a relative version of HMS conjecture for the extended Fano/LG correspondence, with the expectation that it will recover the mirror P=W relation, as outlined below:
\begin{equation*}
    \fbox{\begin{varwidth}{\textwidth}
\centering
Relative Homological \\
Mirror Symmetry
\end{varwidth}} \xrightarrow{HH_a}
\fbox{\begin{varwidth}{\textwidth}
\centering
Relative Hodge \\
Mirror Symmetry
\end{varwidth}} \xrightarrow{signs}
\fbox{\begin{varwidth}{\textwidth}
\centering
Mirror P=W conjecture
\end{varwidth}}
\end{equation*}
\noindent Here $HH_a(-)$ means taking Hochschild homology for $-n \leq a \leq n$ and $signs$ means putting the signs following the Mayer-Vietoris sign rule \eqref{eq:MV sign rule} to recover the $E_1$-page of the spectral sequence for the weight and perverse filtration. For more details on the rank 2 case, we refer to Section \ref{sec:extended fano lg rank2}, building upon the previous discussion.

 \begin{rem}
The term '(extended) Fano/LG correspondence' may sound misleading because one could relax the ampleness condition on a pair $(X,D)$ when discussing mirror symmetry. The minimal condition required for the SYZ proposal \cite{AurouxFanoSYZ}\cite{AurouxFanoSYZLag} is the existence of a Calabi-Yau member $D$ in the linear system $|-K_X|$. Depending on the specific aspect of mirror symmetry under consideration, additional conditions on $-K_X$ may be necessary. For instance, when $X$ is quasi-Fano, a mirror hybrid LG potential is expected to be a Calabi-Yau fibration over a polydisk, not merely a collection of holomorphic functions. With these generalizations in mind, we use the term '(extended) Fano/LG correspondence' liberally.

\end{rem}

\subsection{Outline}
In Section 2, we provide basic preliminaries on Deligne's Hodge theory and the perverse Leray filtration. Section 3 introduces the notion of a hybrid LG model, the primary focus of this paper, and examines its topological properties. Section 4 explores the Hodge theory of hybrid LG models, introducing two different LG Hodge numbers and extending the conjecture of Katzarkov-Kontsevich-Pantev\cite{KKPbogomolov}. Section 5 contains speculative discussions and proposes a relative version of homological mirror symmetry for the extended Fano/LG correspondence. We also explain how the mirror P=W conjecture can be deduced from this relative version. In the Appendix, we review the theory of cohomological mixed Hodge complexes, which will be employed in Section 4. Additionally, we provide a combinatorial description of the perverse Leray filtration, which holds independent interest.
\subsection{Acknowledgement}
I extend my sincere gratitude to Tony Pantev for suggesting this project and providing his consistent support. My appreciation goes to Andrew Harder for sharing his expertise on the topic of hybrid LG models and for his work \cite{Harderunpublished}. Special thanks to Jeff Hicks, Andrew Hanlon, and Sheel Ganatra for explaining their work on the construction of partially wrapped Fukaya categories associated to LG multi-potentials. Additionally, I would like to thank Cheuk-Yu Mak, Nick Sheridan, Ron Donagi, and Charles Doran for their willingness to engage in many meaningful discussions. I am indebted to the anonymous referees for their valuable assistance in reorganizing and providing insightful comments to enhance this paper. Lastly, I gratefully acknowledge the financial support received from the Leverhulme Trust.

\section{The mirror P=W conjecture}\label{sec:mirrorPW}
\subsection{The weight filtration}\label{sec:weight}
We review Deligne's construction of a $\Q$-mixed Hodge structure on the cohomology of a quasi-projective variety by following the exposition in \cite{peters2008mixed}. Let $U$ be a smooth quasi-projective variety over $\C$ and assume that we have a good compactification $(X,D)$. Recall that a pair $(X,D)$ is called a good compactification of $U$ if $X$ is a smooth and compact variety and $D$ is a simple normal crossing divisor. Let $j:U \to X$ be a natural inclusion. Consider the logarithmic de Rham complex
\begin{equation*}
    \Omega^{\bullet}_X(\log D) \subset j_*\Omega^{\bullet}_U
\end{equation*}
Locally at $p \in D$ with an open neighborhood $V \subset X$ with coordinates $(z_1, \cdots, z_n)$ in which $D$ is given by $z_1\cdots z_k=0$, one can see
\begin{equation*}
    \begin{aligned}
    &\Omega^1_X(\log D)_p=\cO_{X,p}\frac{dz_1}{z_1} \oplus \cdots \oplus \cO_{X,p}\frac{dz_k}{z_k} \oplus \cO_{X,p}dz_{k+1} \oplus \cdots \oplus \cO_{X,p}dz_n\\
&\Omega^r_X(\log D)_p = \bigwedge^r\Omega^1_X(\log D)_p
    \end{aligned}
\end{equation*}

\noindent There are two natural filtrations on the logarithmic de Rham complex $(\Omega^{\bullet}_X(\log D),d)$. 
    \begin{enumerate}
        
        \item The decreasing Hodge filtration $F^\bullet$ on $\Omega^{\bullet}_X(\log D)$ is defined by \begin{equation*}
            F^p\Omega^\bullet_X(\log D):= \Omega^{\geq p}_X(\log D) 
        \end{equation*}
        \item The increasing weight filtration $W_\bullet$ on $\Omega^{\bullet}_X(\log D)$ is defined by 
        \[W_m\Omega^r_X(\log D):= \begin{cases}
        0 &  m<0\\
        \Omega^r_X(\log D) &  m \geq r\\
        \Omega^{r-m}_X\wedge \Omega^m_X(\log D) &  0 \leq m \leq r
        \end{cases}
        \]
    \end{enumerate}
\begin{thm}\cite[Theorem 4.2]{peters2008mixed}
    \begin{enumerate}
        \item The logarithmic de Rham complex $\Omega^{\bullet}_X(\log D)$ is quasi-isomorphic to $j_*\Omega^{\bullet}_U$:
        \begin{equation*}
         H^k(U;\C)=\mathbb{H}^k(X, \Omega^{\bullet}_X(\log D))
        \end{equation*}.
        \item The decreasing filtration $F^\bullet$ on $\Omega^{\bullet}_X(\log D)$ induces the filtration in cohomology
        \begin{equation*}
          F^pH^k(U;\C)=\Image(\mathbb{H}^k(X, F^p\Omega^{\bullet}_X(\log D)) \to H^k(U;\C))
        \end{equation*}
        which is called the \textbf{Hodge filtration} on $H^\bullet(U)$. Similarly, the increasing filtration $W_\bullet$ on $\Omega_X^\bullet( \log D)$ induces the filtration in cohomology
        \begin{equation*}
          W_mH^k(U;\C)=\Image(\mathbb{H}^k(X, W_{m-k}\Omega^{\bullet}_X(\log D)) \to H^k(U;\C))
        \end{equation*}
        which is called the \textbf{weight filtration} on $H^\bullet(U)$. In particular, the weight filtration can be defined over the field of rational numbers $\Q$ and we denote it by $W^\Q_\bullet$.
        \item The package $(\Omega^{\bullet}_X(\log D), W^\Q_{\bullet}, F^{\bullet})$ gives a $\Q$-mixed Hodge structure on $H^k(U;\C)$.
    \end{enumerate}
\end{thm}
The key properties of these two filtrations are the degeneration of the associated spectral sequences. More precisely, we have 

\begin{proposition}\cite[Theorem 4.2, Proposition 4.3]{peters2008mixed}
    \begin{enumerate}
    \item The spectral sequence for $(\mathbb{H}(X, \Omega^{\bullet}_X(\log D)), F^{\bullet})$ whose $E_1$-page is given by 
    \begin{equation*}
         E_1^{p,q}=\mathbb{H}^{p+q}(X, \Gr_F^p\Omega^{\bullet}_X(\log D))
    \end{equation*}
    degenerates at the $E_1$-page. Thus, we have
    \begin{equation*}
        \Gr_{F}^p\mathbb{H}^{p+q}(X, \Omega^{\bullet}_X(\log D))=\mathbb{H}^{p+q}(X, \Gr_F^p\Omega^{\bullet}_X(\log D)).
    \end{equation*}
    \item The spectral sequence for $(\mathbb{H}(X, \Omega^{\bullet}_X(\log D)), W_{\bullet})$ whose $E_1$-page is given by 
    \begin{equation*}
        E_1^{-m, k+m}=\mathbb{H}^k(X, \Gr^W_m\Omega^{\bullet}_X(\log D)) 
    \end{equation*}
    degenerates at the $E_2$-page and the differential $d_1:E_1^{-m, k+m} \to E_1^{-m+1, k+m}$ is strictly compatible with the filtration $F_\bullet$. In other words, 
    \begin{equation*}
        E_2^{-m, k+m}=E_{\infty}^{-m, k+m}=\Gr^{W}_{m+k}\mathbb{H}^k(X, \Omega^{\bullet}_X(\log D)).
    \end{equation*}
    \end{enumerate}
\end{proposition}
Next, we recall the notion of the Tate twist. 
\begin{defin}
Let $V=(V_\Q, W_\bullet, F^\bullet)$ be a $\Q$-mixed Hodge structure. For $m \in \Z$, we define the $m$-th \textbf{Tate twist} of $V$ by setting  $V(m):=(V_\Q(m), W(m)_\bullet, F(m)^\bullet)$ where $V_\Q(m):=(2\pi i)^mV_\Q$ and 
\begin{equation*}
    W(m)_k:=W_{k+2m} \quad F(m)^p:=F^{m+p}
\end{equation*}
for all $k$ and $p$.
\end{defin}

To compute the mixed Hodge structure, we provide a geometric description of the spectral sequence. For a given normal crossing divisor $D$, let's denote the irreducible components of $D$ by $D_i$. We set $D(k)$ to be the disjoint union of $k$-tuple intersections of the components of $D$ and $D(0)=X$. Also, for $I=(i_1, \cdots, i_m)$ and $J=(i_1, \cdots, \Hat{i_j}, \cdots, i_m)$, there are inclusion maps 
\begin{align*}
    & \iota^I_J:D_I \hookrightarrow D_J\\
    \iota^m_j=\bigoplus_{|I|=m}& \iota_J^I:D(m) \hookrightarrow D(m-1)
\end{align*}
which induces canonical Gysin maps on the level of cohomology. Therefore, we have
\begin{equation}\label{eq:MV sign rule}
    \gamma_m=\bigoplus^m_{j=1}(-1)^{j-1}(\iota^m_j)_!:H^{k-m}(D(m))(-m) \rightarrow H^{k-m+2}(D(m-1))(-m+1) 
\end{equation}
where $(-)_!$ is the Gysin morphism. Here we call this sign convention by \textit{the Mayer-Vietoris sign rule}, which is unique up to $\pm1$. Under the residue map, this gives a geometric description of the differential $d_1:E_1^{-m, k+m} \to E_1^{-m+1, k+m}$ as follows:
\begin{proposition}\cite[Proposition 4.7]{peters2008mixed}\label{prop:CD d_1 spectral}
The following diagram is commutative
\begin{equation}\label{eq:E1 spectral weight}
    \begin{tikzcd}
        E_1^{-m, k+m} \arrow{r}{res_m} \arrow{d}{d_1} & H^{k-m}(D(m);\C)(-m) \arrow{d}{-r_m} \\
        E_1^{-m+1, k+m} \arrow{r}{res_{m-1}} & H^{k-m+2}(D(m-1);\C)(-m-1)
    \end{tikzcd}
\end{equation}
    where $res_m$ is the residue map for all $m \geq 0$.
\end{proposition}

Note that all morphisms in the diagram (\ref{eq:E1 spectral weight}) are compatible with the Hodge filtration $F^\bullet$. This description provides several computational tools as well as functorial properties of mixed Hodge structures under geometric morphisms. For more details, we refer the reader to \cite{peters2008mixed}. We introduce one more piece of  terminology, which will be used later.

\begin{defin}
Let $U$ be a quasi-projective variety. A $\Q$-mixed Hodge structure on $H^k(U)$ is called \textbf{Hodge-Tate} if the weight $2l$-Hodge structure on the associated graded pieces $\Gr^W_{2l}H^k(U)$ is concentrated at $(l,l)$ for $l \geq 0$ and $\Gr^W_{2l+1}H^k(U)=0$.
\end{defin}

\begin{example}
   Let $U$ be a $n$-dimensional torus $(\C^*)^n$. It admits a good compactification $(\CP^n, D)$ where $D$ is the toric anti-canonical divisor. Then the $\Q$-mixed Hodge structure on $H^k(U)$ is given by both the long exact sequence of the pair $(\CP^n, D)$ and Poincar\'e duality:
   \begin{equation*}
   \begin{aligned}
       \cdots \subset 0=W^{\Q}_{2k-1} \subset W^{\Q}_{2k}=H^k(U;\Q)\\
       \cdots \subset 0=F^{k+1} \subset F^k=H^k(U;\C).
   \end{aligned}
   \end{equation*}
   In other words, the only non-trivial associated graded piece is $\Gr^F_k\Gr_{W}^{2k}H^k(U;\C)=\C^{\binom{n}{k}}$ for all $k \geq 0$. Using Tate twists, one can write down the $\Q$-mixed Hodge structure as 
   \begin{equation*}
       H^k(U;\Q) \cong \Q(-k)^{\binom{n}{k}}
   \end{equation*}
   and it is clearly of Hodge-Tate type. 
\end{example}

\begin{example}
   Let $X_n$ be a del Pezzo surface of degree $9-n$ which is a blow up of $\CP^2$ at general $n (\leq 8)$ points. Take a smooth anti-canonical divisor $E \subset X_n$. By adjunction, the complement $U_n=X_n\setminus E$ is Calabi-Yau. The $\Q$-mixed Hodge structure on $H^k(U_n)$ is given by 
   \begin{enumerate}
       \item $H^0(U_n;\Q) \cong \Q(0)$.
       \item $H^2(U_n;\Q)$ sits in a short exact sequence
       \begin{equation*}
           0 \to \Q(-1)^{n-1} \to H^2(U_n;\Q) \to H^1(E;\Q)(-1) \to 0
       \end{equation*}
       \item $H^i(U_n;\Q)=0$ for $i \neq 0, 2$
   \end{enumerate}
   It implies that the $\Q$-mixed Hodge structure on $H^2(U_n;\Q)$ is of Hodge-Tate type with non-trivial weight graded pieces $\Gr_2^W$ and $\Gr_4^W$.
\end{example}

\subsection{The perverse filtration}
The notion of perversity was invented by Goresky-MacPherson \cite{GoreskyMacphersonI} \cite{GoreskyMacphersonII} to capture the singular behavior of an algebraic variety or sheaves via cohomology theories. In particular, it plays a crucial role in understanding the topology of algebraic maps, and for more details, we refer to \cite{Cataldotopologyofalgebraic}. Let $X$ be an algebraic variety or scheme over $\C$. In the case where $X$ is singular, a sheaf $F$ on $X$ behaves unexpectedly over the singular locus, making it challenging to  understand the ring structure on cohomology. To resolve this issue, instead of studying sheaves on $X$, one can introduce the notion of constructible sheaves which becomes locally constant over each singularity stratum. These sheaves form a well-defined triangulated category, denoted by $D^b_c(X)$, referred to as the (bounded) derived category of constructible sheaves over $\C$.\footnote{One can work with constructible sheaves over $\Q$ as well in this section.}
\begin{definition}
    Let $X$ be an algebraic variety (or scheme) and $D_c^b(X)$ be a derived category of constructible sheaves on $X$. An object $K^{\bullet} \in D^b_c(X)$ is called a \textbf{perverse sheaf} if it satisfies following two dual conditions:
    \begin{enumerate}
        \item (Support condition) $\dim \supp(\mathcal{H}^i(K^{\bullet})) \leq -i$
        \item (Cosupport condition)
        $\dim \supp(\mathcal{H}^i(\mathbb{D}K^{\bullet})) \leq i$ where $\mathbb{D}:D^b_c(X) \to D^b_c(X)$ is the Verdier dualizing functor. 
    \end{enumerate}
\end{definition}

Verdier's dualizing functor on $D^b_c(X)$ is defined as  $\mathbb{D}=\Hom_{\cO_X}(-, p^{!}(\C_{pt}))$ where $p:X \to pt$ is a trivial map. We call $p^{!}(\C_{pt})$ the dualizing complex of $X$. In particular, if $X$ is non-singular of complex dimension $n$, the the dualizing complex is given by $\C_X[2n]$. Note that
the subcategory $\mathcal{P}(X)$ of perverse sheaves on $X$ is an abelian category. Also, the support and cosupport conditions induce the so-called perverse $t$-structure $({}^\mathfrak{p}D^{b, \geq 0}_c(X), {}^\mathfrak{p}D^{b, \leq 0}_c(Y))$ on $D^b_c(X)$ whose heart is  $\mathcal{P}(X)$. Explicitly, this is given by 
\begin{enumerate}
    \item $K^{\bullet} \in {}^\mathfrak{p}D^{b, \leq 0}_c(X)$ if and only if $K$ satisfies the support condition. Also, ${}^\mathfrak{p}D^{b, \leq n}_c(X):={}^\mathfrak{p} D^{b, \leq 0}_c(X)[-n]$
    \item $K^{\bullet} \in {}^\mathfrak{p} D^{b, \geq 0}_c(X)$ if and only if $K$ satisfies the cosupport condition. Also, ${}^\mathfrak{p}D^{b, \geq n}_c(X):={}^\mathfrak{p} D^{b, \geq 0}_c(X)[-n]$
\end{enumerate}
We denote ${}^\mathfrak{p}\tau_{\leq n}:D_c^b(X) \to {}^\mathfrak{p} D^{b, \leq n}_c(X)$ (resp. ${}^\mathfrak{p}\tau_{\geq n}:D_c^b(X) \to {}^\mathfrak{p} D^{b, \geq n}_c(X)$) the natural truncation functor. This induces \textit{perverse cohomology functors} ${}^\mathfrak{p}\mathcal{H}:D^b_c(X) \to \mathcal{P}(X)$ defined by ${}^\mathfrak{p}\mathcal{H}^k:={}^\mathfrak{p} \tau_{\leq 0} \circ {}^\mathfrak{p}\tau_{\geq 0} \circ [k]$. Applying the perverse truncation, one can define the perverse filtration on the hypercohomology of a constructible sheaf $K^\bullet$ on $X$ as follows.
\begin{defin}
For any $K^\bullet \in D^b_c(X)$, the perverse filtration $P_\bullet$ on $\bbH^k(X, K^\bullet)$ is defined to be
\begin{equation*}
    P_b\bbH^k(X, K^\bullet):=\Image(\bbH^k(X, {}^\mathfrak{p}\tau_{\leq b}K^\bullet) \to \bbH^k(X, K^\bullet))
\end{equation*}
\end{defin}

Let $f:X \to Y$ be a morphism of smooth varieties. Then we define the perverse (Leray) filtration $P^f_\bullet$ associated to $f$ on the cohomology $H^\bullet(X, \C)$ by setting
\begin{equation*}
    P^f_lH^k(X,\C):=P_l\bbH^k(Y, Rf_*\C)
\end{equation*}
Due to the decomposition Theorem \cite{BBDperverse}\cite{Cataldotopologyofalgebraic}, the spectral sequence for the perverse filtration degenerates at the $E_2$ page when $f$ is proper. In other words, 
\begin{thm}\label{thm: perverse degenerate E_2}\cite[Corollary 14.41]{peters2008mixed}
    If $f$ is proper, then the perverse Leray spectral sequence degenerates at the $E_2$ page. In other words, we have 
    \begin{equation}
        \Gr^{P^f}_l\bbH^{k}(X,\C)=E_2^{k-l,l}=\bbH^{k-l}(Y, {}^\mathfrak{p}\mathcal{H}^l(Rf_*\C))
    \end{equation}
\end{thm}

\begin{thm}\cite{Cataldolef}\label{thm:cataldo}
    Let $f:X \to Y$ be a morphism of varieties with $Y$ being affine and $K^\bullet \in D_c^b(X)$ be a constructible sheaf. Then for a flag of subvarieties $Y_{\bullet} \subset Y$ obtained by the linear hyperplane sections in sufficiently general position, the perverse filtration is given by
    \begin{equation*}
        P^f_lH^*(X,K^\bullet)=\Ker\{(H^*(X,K^\bullet) \to H^*(X_{l-*+1}, K^\bullet|_{X_{l-*+1}})\}
    \end{equation*}
    where $Y_r \subset Y$ has the codimension $r$ and $X_\bullet \subset X$ is the pre-image flag. In particular, if $f$ is proper, then the spectral sequence for the filtration $P^f_\bullet$ degenerates at the $E_2$-page.
\end{thm}
In particular, the length of the perverse filtration depends on the dimension of $Y$,
\begin{equation*}
     0=P^f_{k-1} \subset P^f_{k} \subset \cdots \subset P^f_{\dim(Y)+k-1} \subset P^f_{\dim(Y)+k}=H^k(X,K^\bullet).
\end{equation*}
For any quasi-projective variety $U$, there exists a universal affinization map 
\begin{equation*}
    \mathsf{Aff}: U \to \Spec H^0(U, \cO_U)
\end{equation*}
which factors through any morphisms from $U$ to an affine variety. We denote $P_\bullet$ (drop the notation of the morphism) the perverse filtration associated to $\mathsf{Aff}$ and call it the \textbf{canonical perverse Leray filtration}. 
Moreover, the universal property of the affinization map implies that it factors through any morphism to an affine space $f:U \to \Spec(A)$. Therefore, if the affinization map $\mathsf{Aff}$ is proper, then $P_\bullet$ is the same as $P^h_\bullet$ for any proper morphism $h:U \to \Spec(A)$ because the induced morphism $\Spec H^0(U,\cO_U) \to \Spec(A)$ is finite.

\begin{example}
   Let $U$ be a $n$-dimensional affine torus $(\C^*)^n$. The affinization map is an isomorphism $\Aff: U \to (\C^*)^n$. By iteratively applying the Lefschetz hyperplane theorem, one can see that the only non-trivial associated graded pieces of the canonical perverse Leray filtration is $\Gr^P_{n}H^{k}(U;\C)=H^{k}(U;\C)$. 
\end{example}

\begin{example}
   Let $X$ be a smooth projective variety with a smooth ample divisor $D$. Since the complement $U:=X \setminus D$ is always affine, we can iteratively apply the Lefschetz hyperplane theorem to argue that $\Gr_i^PH^k(U;\C)$ is trivial except when $i=\dim(U)$. 
\end{example}

\subsection{The mirror P=W conjecture}
Let $U$ be a quasi-projective variety. Note that the canonical perverse Leray filtration on the cohomology of $U$ is compatible with Deligne's mixed Hodge structure. Therefore, we can define the \textit{perverse-mixed Hodge polynomial}, which encodes information about all the Hodge numbers that are refined with respect to both the canonical perverse and weight filtration. 
\begin{defin} \label{def:PMH poly}
For any (non-singular) quasi-projective variety $U$ over $\C$ of dimension $n$, we define a \textbf{perverse-mixed Hodge polynomial} of $U$ as follows
    \begin{equation*}
        PW_U(u,t,w,p)=\sum_{\substack{0 \leq s \leq 2n, r\geq 0,\\ 0 \leq b \leq s, 0 \leq a \leq s+b} }\dim_\C (\Gr_F^a\Gr_{s+b}^W \Gr^P_{s+r}H^s(U,\C))u^at^sw^bp^r
    \end{equation*}
\end{defin}

\begin{example}
From the previous examples, we have the perverse-mixed Hodge polynomials of $(\C^\ast)^n$ and $U_n$ as follows:
   \begin{enumerate}
        \item $PW_{(\C^*)^n}(u,t,w,p)=(utw+p)^n$.
        \item $PW_{U_n}(u,t,w,p)=p+uw^2p^2+ut^2w^2+(9-n)uw^2$.
    \end{enumerate}
\end{example}

\begin{conj}\label{conj:mirrorP=W}(Mirror P=W Conjecture) 
     Assume that two $n$-dimensional log Calabi-Yau varieties $U$ and $U^\vee$ are mirror to each other. Then we have the following polynomial identity:
     \begin{equation}\label{eq:PW conjecture full form}
          PW_U(u^{-1}t^{-2}, t,p,w)u^nt^n=PW_{U^\vee}(u,t,w,p)
     \end{equation}
 \end{conj}   
In case that $U$ is affine and the mixed Hodge structure on $H^*(U^\vee)$ is of Hodge-Tate type then the identity (\ref{eq:PW conjecture full form}) simplifies to the following:
\begin{equation*} \dim\Gr_F^q\Gr^W_{p+q+r}H^{p+q}(U)=\dim\Gr_F^{n-q}\Gr^P_{n+p-q+r}H^{n+p-q}(U^\vee)
\end{equation*}
for all $p,q$ and $r$. 
In this article, we specialize to this case, presuming that $U$ is the complement of a Fano pair $(X, D)$. The Hodge-Tateness assumption is also motivated by mirror symmetry (see Section \ref{sec:Extended Fano/LG}).
\section{Hybrid Landau-Ginzburg model}\label{sec:hybridLG}
We introduce the notion of hybrid Landau-Ginzburg (LG) models, which serve as a mirror object for a Fano pair $(X,D)$ where $X$ is a Fano manifold, and $D$ is a simple normal crossing anti-canonical divisor. In this section, we mainly discuss the topological aspects of hybrid LG models. 

We introduce some notations first. Let $h=(h_1, \dots, h_N):Y \to \C^N$ be a $N (\geq 1)$-tuple of holomorphic functions and $(z_1,\dots, z_N)$ be coordinates of the base $\C^N$. For each non-empty subset $I=\{i_1, \dots, i_l\} \subset \{1,\dots, N\}$, we write $h_I=(h_{i_1}, \dots, h_{i_l}):Y \to \C^{|I|}$ and $(z_{i_1}, \dots, z_{i_l})$ for the coordinate of the base $\C^{|I|}$, which implicitly determines the natural inclusion $\C^{|I|} \subset \C^N$. 

\begin{defin}\label{def: weak hybrid LG model general}
	A \textbf{hybrid Landau-Ginzburg (LG) model of rank $N$} is a triple \\
	$(Y, \omega, h=(h_1, h_2, \dots, h_N):Y \to \C^N)$ where 
	\begin{enumerate}
		\item $(Y, \omega)$ is a $n$-dimensional complex K\"ahler Calabi-Yau manifold with a K\"ahler form $\omega \in \Omega^2(Y)$;
		\item $h:Y \to \C^N$ is a proper surjective holomorphic map such that 
		\begin{enumerate}
		    \item (Local trivialization) there exists a constant $R>0$ such that for any non-empty subset $I \subset \{1, \dots, N\}$, the induced map $h_I:Y \to \C^{|I|}$ is a locally trivial fibration over the region $B_I:=\{|z_i| > R | i \in I\}$ with smooth fibers. Furthermore, over $B_I$ we have $v(h_j)=0$ for any horizontal vector field $v \in T^{h_I}Y$ associated to $h_I$ and $j \notin I$;
		    \item (Compatibility) for $I \subset J$,  such local trivializations are compatible under the natural inclusions $B_J \times \C^{N-|J|} \subset B_I \times \C^{N-|I|} \subset \C^N$. 
		\end{enumerate}
	\end{enumerate}
   	We call $h:Y \to \C^N$ a \textbf{hybrid Landau-Ginzburg (LG) potential}. 
\end{defin}
Here, the rank is defined to be the number of components of $h$. We often omit it from the name if it is clear from the context. Also, when $N=1$, this definition corresponds to the classical notion of a \textit{proper} Landau-Ginzburg model.

\begin{rem}
In the definition of a hybrid LG model, the role of the Kähler form $\omega$ is vacuous and one may state the condition $(1)$ as $Y$ being a  quasi-projective Calabi-Yau manifold. However, from the perspective of mirror symmetry, it is more natural to include a Kähler form in the definition, so we retain it in the notation. Discussion of the symplectic geometric aspects of hybrid LG models is postponed until Section \ref{sec:Extended Fano/LG}. One can also drop the condition of the surjectivity of $h:Y \to \C^N$ and consider all the morphisms in the definition as fibrations over the image. Then the arguments we will make in this section are valid for this case. 

\end{rem}

    Let us explain the conditions in Definition \ref{def: weak hybrid LG model general} more. For each $I \subset \{1, \dots, N\}$, the condition (a) is equivalent to the following: for $p \in U \subset B_I$ with $U$ open and simply-connected and $j \notin I$, we have the commutative diagram
    \begin{equation}\label{eq:local triv}
        \begin{tikzcd}
            h_I^{-1}(U) \arrow[r, "\iota_{I,U}", "\cong"'] \arrow[d, "h_j"] & h_I^{-1}(p) \times U \arrow[d, "h_j \times \mathsf{id}"] \\
            \C & \C \times U \arrow[l, "pr_1"]
        \end{tikzcd}
    \end{equation}
    which satisfies the usual compatibility conditions for trivializations. Also, let $B_I \times \C^{N-|I|} \hookrightarrow \C^N$ be the natural inclusion whose image is given by $\{|z_i|>R|i \in I\}\subset \C^N$. For $I \subset J$, consider the natural inclusions $B_J \times \C^{N-|J|} \hookrightarrow B_I \times \C^{N-|I|}$ and the natural projection $\pi^J_I:B_J \to B_I$. For $p \in V \subset B_J$ with $V$ open and simply-connected,  and $j \notin J$, we have the commutative diagram of local trivializations (\ref{eq:local triv})
    \begin{equation}
\begin{tikzcd}\label{diag: compatibility IJ}
     & h_J^{-1}(V) \arrow[dd, "h_j" near start] \arrow[rr, hook] \arrow[dl, "\iota_{J,V}", "\cong"'] && h_I^{-1}(\pi^J_I(V)) \arrow[dl, "\iota_{I,\pi_I^J(V)}", "\cong"'] \arrow[dd, "h_j" near start] \\
    h_J^{-1}(p) \times V \arrow[rr, hook] \arrow[dd, "h_j \times \mathsf{id}" near start] && h_I^{-1}(\pi_I^J(p)) \times \pi_I^J(V) \arrow[dd, "h_j \times \mathsf{id}" near start] & \\
    & \C \arrow[rr, swap, "\mathsf{id}" near start] && \C \\
    \C \times V \arrow[ur, "pr_1"] \arrow[rr, "\mathsf{id} \times \pi_I^J"] && \C \times \pi_I^J(V) \arrow[ur, "pr_I^J"] &  
\end{tikzcd}
\end{equation}
    which satisfies usual compatibility conditions for trivializations as well. If we write $D(I,J)$ for the diagram (\ref{diag: compatibility IJ}), we should have $D(J,K) \circ D(I,J)=D(I,K)$ for $I \subset J \subset K$. Let's write $Y_I$ for a generic fiber of $h_I:Y \to \C^{|I|}$ and $h_{Y_I}:Y_I \to \C^{N-|I|}$ for the restriction of $h$ into $Y_I$. Then the induced triple $(Y_I, \omega|_{Y_I}, h_{Y_I})$ is a hybrid LG model of rank $N-|I|$. From this point of view, the condition $(a)$ in Definition \ref{def: weak hybrid LG model general} is rephrased as the condition that $h_I:Y \to \C^{|I|}$ is a local trivialization of the induced hybrid LG models of rank $N-|I|$ near infinity.

    Associated to the hybrid LG model $(Y, \omega, h:Y \to \C^N)$, we define the ordinary LG model to be a triple $(Y, \omega, \sw:=\Sigma \circ h:Y \to \C)$ where $\Sigma:\C^N \to \C$ is the summation map. To justify this terminology, we need to show that $\sw:Y \to \C$ is locally trivial near infinity.

\begin{prop}\label{prop:Gluing Property in general}(Gluing property)
    Let $(Y, \omega, h:Y \to \C^N)$ be a hybrid LG model and $H$ be a generic hyperplane in the base $\C^N$, which is not parallel to any coordinate lines. There exists an open cover $\{U_i\}_{i=1}^N$ of $H$ such that for any non-empty subset $I \subset \{1, \dots, N\}$, the induced map $h^{-1}(U_I) \to U_I$ is isotopic to the induced hybrid LG potential $h_{Y_I}:Y_I \to \C^{N-|I|}$ which is linear along the base. 
\end{prop}

\begin{proof}
    Take a hyperplane $H=\{a_1z_1 + \cdots + a_Nz_N=M\}$ where $a_i \neq 0$ for all $i$. By changing the coordinate $z_i \mapsto z_i/a_i$, we reduce to the case where $a_i=1$ for all $i$. We also further reduce to the case when $M$ is real due to the rotational symmetry. Suppose that $M >NR$. First, note that $H \cap (\cap_{i=1}^N\{|z_i|\leq R\})=\emptyset$. Let $R_i=\{Re(z_i)>R\}$ and the simply-connected region
    \begin{equation*}
        U_i:=\{Re(z_i) > R\} \cap H=\{Re(z_1 + \cdots + \hat{z_i}+ \cdots + z_n)<M-R\} \cap H
    \end{equation*}
    for each $i$. Since $U_i \subset \{|z_i|>R\}$, one can project $U_i$ to the locus $\{z_i=2R\}$ inside the region $\{|z_i|<R\}$. The image of the projection is $V_i:=\{z_i=2R, Re(z_1+\cdots+\hat{z_i}+\cdots+z_N)<M-R\}$ which contains $\bigcap_{j \neq i}\{|z_j| \leq R\}$. Therefore, this projection identifies $h:h^{-1}(U_i) \to U_i$ with $h:h^{-1}(V_i) \to V_i$ due to the local triviality of the hybrid LG model. Moreover, the latter map is completed to $h_{Y_i}:Y_i \to \C^{N-1}$ by the inductive argument. In general, for each $I$, $U_I= \cap_{i \in I}U_i$ is non-empty and simply-connected. Since $U_I \subset \{|z_i|>R, i \in I\}$, one can apply the same argument to get the conclusion. 
\end{proof}

\begin{defin}
    Let $(Y, \omega, h:Y \to \C^N)$ be a hybrid LG model. We define the induced triple $(Y, \omega, \sw:Y \to \C)$ to be the ordinary LG model associated to the hybrid LG model $(Y, \omega, h)$ and denote a generic fiber of $\sw$ by $Y_{sm}$.
\end{defin}

On the cohomology level, the gluing property implies that the cohomology groups of $\pi^{-1}(H)$ are all (non-canonically) isomorphic to those of the normal crossing union of $Y_i$'s. We will use this fact to study the perverse Leray filtration associated to $h:Y \to \C^N$ on $H^*(Y)$. 

Consider a flag of affine subspaces in $\C^N$:
\begin{equation*}
    \mathfrak{H}:0=H_{-N-1} \subset H_{-N} \subset \cdots \subset H_{-1} \subset H_0=\C^N
\end{equation*}
which is generic enough to induce the perverse filtration (Theorem \ref{thm:cataldo}). Without loss of generality, we may further assume that each $H_{-l}$ is not parallel to any coordinate lines. We write $Y_{sm^{(l)}}$ for $h^{-1}(H_{-l})$ so that the flag filtration on $H^*(Y)$ induced by the flag of subvarieties 
\begin{equation*}
    0 \subset Y_{sm^{(N)}} \subset \cdots \subset Y_{sm^{(1)}}\subset Y
\end{equation*}
is the same as the perverse filtration $P^h_\bullet$. Then the $E_1$-page of the spectral sequence is given by
\begin{equation}\label{eq:E_1-spectral perverse}
        H^{a-N}(Y_{sm^{(N)}}) \xrightarrow{d_1} H^{a-N+1}(Y_{sm^{(N-1)}}, Y_{sm^{(N)}}) \xrightarrow{d_1} \cdots \xrightarrow{d_1} H^{a-1}(Y_{sm^{(1)}}, Y_{sm^{(2)}}) \xrightarrow{d_1} H^a(Y, Y_{sm^{(1)}})
\end{equation}

\noindent We use the same notation in the proof of Proposition \ref{prop:Gluing Property in general}. Take open (simply-connected) regions $\{R_i \subset \C^N|i=1, \dots, N\}$ which induces the open covering of $H_{-1}$, $\{U_i:=R_i \cap H_{-1}|i=1, \dots, N\}$ that yields the gluing property. Let $V_i=\{z_i=const\}$ be the region that $U_i$ projects to. Due to the generality of the flag, we may assume that $H_{-2} \cap U_i \subset U_i$ projects to a hyperplane that is contained in $\cup_{j \neq i}R_j \cap V_i$ for all $i$. As both $H_{-1}$ and $H_{-2}$ are not parallel to any coordinate lines, this can be done by scaling $M$ sufficiently large to place $H_{-2}$ far enough from the each coordinate line. It ensures that $Y_{sm}\cap h^{-1}(U_i)$ is isotopic to $Y_{i,sm}$ for each $i$. Inductively, for each $k$, we could assume that the collection of regions $\{R_i \subset \C^N|i=1, \dots, N\}$ yields the gluing property for $H_{-k}$ in $V_I:=\cap_{i \in I}V_i$ for any $|I|=k-1$. Then the gluing property implies the following:
\begin{equation*}
    Y_{sm^{(k)}} \cap h^{-1}(U_I) \cong \begin{cases}
    Y_{I,sm^{(k-|I|)}} &|I| < k \\
    Y_I & |I| \geq k
    \end{cases}
\end{equation*}

\begin{lem}\label{lem:gluing description of flags}
For any $a \geq 0, k \geq 1$, the relative cohomology $H^a(Y_{sm^{(k)}}, Y_{sm^{(k+1)}})$ is isomorphic to $\bigoplus_{|I|=k}H^a(Y_I, Y_{I,sm})$.
\end{lem}

\begin{proof}
Take the (simply-connected) open region $\{R_i \subset \C^N|i=1, \dots, N\}$ and the induced cover $\{U_i=R_i \cap H_{-1}\}$ as above. When $k=1$, the Mayer-Vietoris sequence with respect to the open cover $\{h^{-1}(U_i)\}$ and the gluing property imply that $H^a(Y_{sm^{(1)}}, Y_{sm^{(2)}}) \cong \bigoplus_{i=1}^N H^a(Y_i, Y_{i, sm^{(1)}})$ where $H^a(Y_i, Y_{i,sm^{(1)}}) \cong  H^a(Y_i,Y_{i,sm})$. In general, we apply the Mayer-Vietoris sequence to the cohomology group\\
$H^a(Y_{sm^{(k)}}, Y_{sm^{(k+1)}})$ with the induced open cover. The $E_1$-page of the spectral sequence is given by 
\begin{equation*}
    \bigoplus_{|I|=1}H^a(Y_{I,sm^{(k-1)}},Y_{I,sm^{(k)}}) \xrightarrow{d_1} \dots \xrightarrow{d_1} \bigoplus_{|I|=k-1}H^a(Y_{I,sm^{(1)}},Y_{I,sm^{(2)}}) \xrightarrow{d_1}  \bigoplus_{|I|=k}H^a(Y_{I}, Y_{I,sm^{(1)}}) \rightarrow 0
\end{equation*}
By induction, each direct summand is the direct sum of $H^a(Y_J, Y_{J,sm})$ for some $J$ with $|J|=k$. Then $d_1$ becomes the alternating sum of the identity morphisms where the signs are determined by the Mayer-Vietoris sign rule (\ref{eq:MV sign rule}). It follows from a simple combinatorial fact that this sequence is exact except at the first term, and $H^a(Y_{sm^{(k)}}, Y_{sm^{(k+1)}})=\ker(d_1) \cong \bigoplus_{|I|=k}H^a(Y_I, Y_{I,sm})$.
\end{proof}

For any $I \subset J$ with $|J|=|I|+1$, we write $\rho^{J}_I$ for the composition of morphisms
\begin{equation*}
    \rho^{J}_I:H^\bullet(Y_J, Y_{J, sm}) \hookrightarrow H^\bullet(Y_{I,sm}, Y_{I,sm^{(2)}}) \to  H^{\bullet+1}(Y_I, Y_{I, sm})
\end{equation*}
where the first one is induced by Lemma \ref{lem:gluing description of flags} and the second one is the connecting homomorphism of the long exact sequence of cohomology groups of the triple 
$(Y_I, Y_{I,sm}, Y_{I, sm^{(2)}})$. Since we choose an open cover globally, Lemma \ref{lem:gluing description of flags} allows one to rewrite the $E_1$-page of the spectral sequence (\ref{eq:E_1-spectral perverse}) as follows:
\begin{equation}\label{eq:E_1-flag}
        H^{a-N}(Y_{sm^{(N)}}) \xrightarrow{d_1} \bigoplus_{|I|=N-1}H^{a-N+1}(Y_I, Y_{I,sm}) \xrightarrow{d_1} \cdots \xrightarrow{d_1} \bigoplus_{|I|=1}H^{a-1}(Y_I, Y_{I,sm}) \xrightarrow{d_1} H^a(Y, Y_{sm})
\end{equation}
where the differential $d_1$ is the signed sum of the induced morphisms $\rho^{J}_I$'s.

On the other hand, consider a collection of coordinate hyperplanes $\{H_i|i=1, \cdots, N\}$ near infinity. For each $l \in \{1, 2, \cdots, N\}$, we take a simple normal crossing subvariety $H\{l\}:=\sqcup_{|I|=l}H_I$ where $H_I=\cap_{i \in I} H_i$ and define $Y\{l\}:=h^{-1}(H\{l\})$. It allows us to consider the flag of subvarieties $\mathcal{F}_h$ in $Y$:
\begin{equation*}
       \mathcal{F}_h: Y\{N\} \subset Y\{N-1\} \subset \cdots \subset Y\{1\} \subset Y\{0\}=Y
\end{equation*}
The $E_1$-page of the spectral sequence of the flag filtration associated to $\mathcal{F}_h$ on $H^*(Y)$ is given by
\begin{equation}\label{eq:E_1-spectral flag}
H^{a-N}(Y\{N\}) \xrightarrow{d_1} H^{a-N+1}(Y\{N-1\}, Y\{N\}) \xrightarrow{d_1} \cdots \xrightarrow{d_1} H^{a-1}(Y\{1\}, Y\{2\}) \xrightarrow{d_1} H^{a}(Y, Y\{1\})
\end{equation}
By applying the Mayer-Vietoris argument, we have $H^{k}(Y\{l\}, Y\{l+1\})\cong \bigoplus_{|I|=l}H^k(Y_I, Y(I))$. Since the triple $(Y_I, \omega|_{Y_I}, h_{Y_I}:Y_I \to \C^{N-|I|})$ is itself a hybrid LG model, we apply the gluing property to show that the sequence \eqref{eq:E_1-spectral perverse} is identified with the sequence \eqref{eq:E_1-spectral flag}. 

\begin{prop}\label{prop:two-perverse-same}
Let $(Y, \omega, h:Y \to \C^N)$ be a hybrid LG model. Then the flag filtration induced by $\mathcal{F}_h$ on $H^*(Y)$ is the same as the perverse filtration $P^h_\bullet$. In particular, the spectral sequence of the flag filtration degenerates at the $E_2$-page. 
\end{prop}

Consider the monodromy $T$ of a generic fiber $h^{-1}(t_1, \dots, t_n)$ of $h:Y \to \C^{N}$ induced by the loop $e^{2\pi it}\cdot(t_1, \dots, t_N)$ $(0 \leq t \leq 1)$. It induces the monodromy action on the relative cohomology groups $H^*(Y_I, Y(I))\cong H^*(Y_I, Y_{I,sm})$ for all $I$. We simply denote the action by the same $T$ and the induced monodromy weight filtration by $W_\bullet^T$. For later use, we also introduce other monodromy actions $T_j$ induced by the loop $(t_1, \dots, e^{2\pi it}\cdot t_j, \dots, t_N)$ $(0 \leq t \leq 1)$ along the $j$-th coordinate hyperplane for $j=1, \dots, n$. For each $I=(i_1, \dots, i_{l})$, we write $T_I=T_{i_1} \circ \cdots \circ T_{i_{l}}$. Clearly, $T=T_1 \circ \cdots \circ T_N$. 

\begin{thm}\label{thm: A-side cubical diagram coh}
    Let $(Y, \omega, h:Y \to \C^N)$ be a hybrid LG model. 
    \begin{enumerate}
        \item For $-n \leq a \leq n$, there is a well-defined cubical diagram 
    \begin{equation}\label{eq:diag cubical coh}
    \begin{aligned}
        H^{a+n-N}(Y\{N\}) &\xrightarrow{\psi_N} \bigoplus_{|I|=N-1}H^{a+n-|I|}(Y_I,Y_{I,sm}) \xrightarrow{\psi_{N-1}} \cdots \\
         &\xrightarrow{\psi_2} \bigoplus_{|I|=1} H^{a+n-|I|}(Y_I, Y_{I,sm}) \xrightarrow{\psi_1} H^{a+n}(Y, Y_{sm})
    \end{aligned}
    \end{equation}
    where $\psi_k$ is the collection of $\rho^{J}_{I}$ for $J \subset I$ with $|J|=|I|+1=k$. Moreover, for each index subset $K$, the monodromy action $T_K$ induces the auto-equivalence of the diagram.
    \item By applying the Mayer-Vietoris sign rule (\ref{eq:MV sign rule}), the diagram \eqref{eq:diag cubical coh} gives rise to the $E_1$-page of the spectral sequence for the perverse filtration associated to $h:Y \to \C^N$ on $H^*(Y)$.
    \end{enumerate}
\end{thm}
We will discuss mirror symmetry motivation of Theorem \ref{thm: A-side cubical diagram coh} in Section \ref{sec:extendedfano lg 5.2}. Additionally, Hodge-theoretic version of Theorem \ref{thm: A-side cubical diagram coh} will be presented as Theorem \ref{thm:mainhodge}.

We finish this section by constructing some examples of hybrid LG models. For this, we introduce mirror symmetry motivation first. Let $X$ be a smooth Fano manifold and $D$ be an effective simple normal crossing anti-canonical divisor with $N$ components $D_1, D_2, \dots, D_N$. For any index set $I=\{i_1, i_2, \cdots, i_m\} \subset \{1, 2, \cdots, N\}$, we define
\begin{equation*}
     D_I:=D_{i_1} \cap \cdots \cap D_{i_m}, \qquad
     D(I):=\left(\bigcup_{j \notin I}D_I\right) \cap D_j
\end{equation*}
For example, if $I=\{1\}$, then $D_{\{1\}}=D_1$ and $D(\{1\})=(D_2 \cup \cdots \cup D_k) \cap D_1$. We also assume that all pairs $(D_I, D(I))$ are Fano. For later use, we denote the normal crossing union of $l$-th intersections as $D{l}=\bigcup_{|I|=l}D_I$ for all $l \geq 0$.

\begin{defin}\label{def:mirrorhybrid}
	A hybrid LG model $(Y, \omega, h:Y \to \C^N)$ is mirror to $(X, D)$ if it satisfies the following mirror relations:
	\begin{enumerate}
	    \item the associated ordinary LG model $(Y, \omega, \sw:Y \to \C)$ is mirror to $(X, D)$;
		\item for $i=1,2, \dots, N$, a hybrid LG model $(Y_i, \omega|_{Y_i}, h_{Y_i}:Y_i \to \C^{N-1})$ is mirror to $(D_i, D(\{i\}))$.
	\end{enumerate}
	A mirror pair of this kind is called a \textbf{Fano mirror pair}, denoted as $(X,D)|(Y, \omega, h)$. 
\end{defin}
\noindent The precise meaning of being a mirror is ambiguous, and we will elaborate on this in Section \ref{sec:Extended Fano/LG}.

\subsection{Construction/Examples}\label{sec:general construction}
We provide some examples of hybrid LG models.

\noindent\textbf{(A) Product}.

For $i=1,2, \dots, N$, let $(X_i,D_{i,sm})$ be a smooth Fano pair, and $(Y_i, \omega_i, h_i:Y_i \to \C)$ be its mirror LG model. The Cartesian product $X:=X_1 \times \cdots \times X_N$ comes with the anti-canonical divisor $\sum_{i=1}^N X_1 \times \cdots \times D_{i,sm} \times \cdots \times X_N$ which has $N$ irreducible components. The corresponding mirror hybrid LG model is given by the triple $(Y, \omega, h:Y \to \C^N)$ where 
     \begin{equation*}
         \begin{aligned}
            & Y=Y_1 \times \cdots \times Y_N, \qquad  \omega=\sum_{i=1}^N pr_i^*\omega_i \\
            &h=(h_1, \cdots, h_N):Y \to \C^N
         \end{aligned}
     \end{equation*}
     and $pr_i:Y \to Y_i$ is the $i$-th projection. By construction, it is straightforward to see that $(Y, \omega, h)$ is a hybrid LG model expected to be mirror to $(X,D)$.

\noindent\textbf{(B) Smooth complete intersections in $\bbP^n$}.

First, consider a Fano pair $(\bbP^n, D)$ where $D=\cup_{i=1}^N D_i$ is a simple normal crossing anti-canonical divisor and $D_i$ is a smooth hypersurface of degree $k_i>0$. Note that $k_1+ \cdots +k_N=n+1$. The choice of such an anti-canonical divisor induces a decomposition of the Hori-Vafa potential for $\bbP^n$ into $N$ potentials. Namely, we have
    $h=(h_1, \cdots, h_N):(\C^*)^n \to \C^N$ where we identify $(\C^*)^n=\left\{\prod_{i=1}^N\prod_{j=1}^{k_i}x_{ij}=1\right\} \subset \C^n_{x_{ij}}$ and $h_i=x_{i,1}+\cdots+x_{i,k_i}$. The discriminant locus of $h:(\C^*)^n \to \C^N$ is given by
    \begin{equation*}
    \Delta_c=\{a_1^{k_1}\cdots a_N^{K_N}=\prod_{i=1}^N(k_i)^{k_i}\} \subset \C^N_{a_1, \dots, a_N}
    \end{equation*}
    Over $\Delta_c$, the singluar locus is given by $\{x_{i,j}=\frac{a_i}{k_i}| i=1,\dots, N\}$. Moreover, the morphism $h:(\C^*)^n \to \C^N$ compactifies to a family of complete intersections over $\C_{a_1, \dots, a_N}^N$ as follows:
    \begin{equation*}
        \mathcal{Y}:=\left\{\begin{array}{l}
        \sum_{j=1}^{k_i}x_{ij}=a_it_i \text{ for } i=1, \dots, N-1\\ 
        x_{N,1}\cdots, x_{N,k_N-1}(a_Nt_N-\sum_{j=1}^{k_N-1}x_{N,j})\prod_{i=1}^{N-1}\prod_{j=1}^{k_i}x_{ij}=t_1^{k_1} \dots t_N^{k_N}
        \end{array}\right\} \subset \bbP^{k_1}\times \cdots \times \bbP^{k_N-1} \times \C^{N}
    \end{equation*}
    We now resolve the singularities of $\cY$ by following the resolution procedure introduced in \cite{Przyjalalgorithm}. The authors construct the algorithm of finding the crepant resolution of the affine hyperplane $L_{\bar{d}, s}$ in \\
$\bbA(a_1, \dots, a_k, \lambda, x_1, \dots, x_s)$ which is of the form
		\begin{equation*}
		    a_1^{d_1} \cdots a_k^{d_k}=\lambda x_1\cdots x_s
		\end{equation*}
		where $\bar{d}=(d_1, \dots, d_k)$ and $d_i >0$ for all $i$. In our case, we have the family of the intersection of $N$ hypersurfaces parametrized by $\C^N$. One cannot directly apply the resolution procedure because $\cY$ is given as a complete intersection. However, since the $N-1$ hypersurfaces are linear, we can view $\cY$ having hypersurface singularities and apply the resolution procedure. For example, when $N=2$, we rewrite $\mathcal{Y}\subset \bbP^{k_1-1}\times \bbP^{k_2-1} \times \C^2_{a_1,a_2}$ by 
		\begin{equation*}
		 \mathcal{Y}=\left\{(a_1t_1-x_{1,2}-\cdots-x_{1,k_1})(a_2t_2-\sum_{j=1}^{k_2}x_{2,j})x_{1,2}\cdots x_{1,k_1}x_{2,1}\cdots x_{2,k_2} =t_1^{k_1}t_2^{k_1}\right\} 
		\end{equation*}
		The resulting variety, denoted by $Y$, is Calabi-Yau since each step in the resolution is crepant. Let $h:Y \to \C^N$ be the induced morphism. The fibration $h$ has (possibly) reducible fibers over $\{a_1\cdots a_N=0\} \subset \C^N$. In particular, since it is a normal crossing union of coordinate hyperplanes, it comes with the canonical stratification 
    \begin{equation*}
    \Delta_N(h) \subset \cdots \subset  \Delta_2(h) \subset \Delta_1(h)=\Delta(h)
\end{equation*}
where $\Delta_r(h)$ is the union of codimension $r-1$ intersections. By construction, the restriction of $h$ over $\Delta_r(h) \setminus \Delta_{r+1}(h)$ is trivial over each connected component so that the differential $dh|_{h^{-1}(\Delta_r(h)\setminus \Delta_{r+1}(h))}$ is surjective. Since the blow-up locus does not intersect with the singular locus over $\Delta_c$, the locus $\Delta_c$ is also a part of discriminant locus of $h:Y \to \C^N$ whose singular locus is isomorphic to $\Delta_c$ itself under $h$. Therefore, for each $I \subset \{1, \dots, N\}$, a generic hyperplane $H_I=\{a_i=const|
i \in I\} \subset \C^N$ intersects transversally with the discrimiant loci so that $h^{-1}(H_I)$ is smooth. We take $R$ to be any constant greater than $\max\{k_i|i=1, \dots N\}$ to achieve the local triviality of each map $h_I:Y \to \C^{|I|}$ over $B_I$. Modulo a choice of a Kähler form $\omega$, this becomes a hybrid LG model $(Y, \omega, h:Y \to \C^N)$.

Furthermore, this construction implies that for each subset $I$, the induced hybrid LG model $(Y_I, \omega|_{Y_I}, h_{Y_I}:Y_I \to \C^{N-|I|})$ is mirror to $(D_I, D(I))$ in the sense that $h_{Y_I}:Y_I \to \C^{N-|I|}$ is a fiberwise compactification of the Hori-Vafa mirror associated to $(D_I, D(I))$. In other words, a mirror hybrid LG model for any smooth complete intersections in the projective space $\bbP^n$ with a simple normal crossing anti-canonical divisor can be given by the above construction. 
    
    \begin{example}
        A typical example is a mirror for a smooth complete intersection in the projective space with the toric anti-canonical divisor. The mirror hybrid LG model is given by the Hori-Vafa mirror whose multi-potential consists of each monomial. In this case, we do not need to take a fiberwise compactification.  
    \end{example}

\begin{example}\label{ex:conic+line}
	Consider a Fano pair with $(\bbP^2,D)$ where $D=C\cup L$ is the union of a conic and a line. By the Hori-Vafa construction, modulo a choice of a Kähler form, one can get a hybrid LG model $(Y,\omega ,h:Y \to \C^2)$ such that 
	\begin{enumerate}
		
		\item The hybrid LG model $h:Y \to \C^2$ is a branched double cover of $\C^2$ whose discriminant locus is 
		\begin{equation*}
		\{a^2b=4\} \cup  \{b=0\}
		\end{equation*}
		
		\item More explicitly, 
		\begin{equation*}
		\begin{tikzcd}
		    Y:=\{x+y=az, z^2=bxy\} \subset \bbP^2_{x,y,z} \times \C_a \times \C_b \arrow[d] \\
		    \C_a \times \C_b
		\end{tikzcd}
		\end{equation*}
	
		\item The original LG potential $\sw=\Sigma \circ h:Y \to \C$ can be obtained by removing two sections of the elliptic fibration over $\bbP^1$ which is mirror to $(\bbP^2, D_{sm})$ (See Example \ref{ex:Mirror Symmetry surface}). This ordinary LG model is the same with one introduced in \cite{AurouxFanoSYZ}. 
	\end{enumerate}
\end{example}

\begin{example}\label{ex:cubic surface}
	Consider a Fano pair  $(\bbP^3, D)$ where $D=C \cup H$ is the union of cubic surface $C$ and linear hypersurface $H$. From the above construction, a hybrid LG model $(Y, \omega, h:Y \to \C^2)$ is an elliptic fibration over $\C^2$ whose discriminant locus is 
	    \begin{equation*}
		\{a^3b=27\} \cup \{b=0\}
		\end{equation*}
    and singular fibers are of type $A_2$ over each irreducible component. Note that over each generic coordinate line, the restriction of the hybrid LG potential $h:Y \to \C^2$ induces a mirror LG model to the corresponding divisor component. 
\end{example}

\section{Hodge Theory of hybrid Landau-Ginzburg models}\label{sec:deform}
In this section, we explore algebro-geometric aspects of hybrid LG models, particularly focusing on the Hodge theory of hybrid LG models. As the choice of a Kähler form does not play a crucial role here, we omit it from the notation. We will see that the combinatorial description of the perverse filtration naturally appears in this context, which turns out to be useful to study the behavior of the perverse filtration under the cup product (see Theorem \ref{thm:multiplicativity}).

\subsection{The classical case}
We begin by reviewing the Hodge theory of classical LG models discussed in \cite{KKPbogomolov, ShamotoHodgeTate, ESYIrregularHodge}. 

\begin{defin}\label{def:compactified LG classic}
    Let $(Y, \sw:Y \to \C)$ be a LG model. A compactification of $(Y, h:Y \to \C)$ is a datum $((Z, D_Z), f, \vol_Z)$ where
    \begin{enumerate}
        \item $Z$ is a smooth projective variety and $f:Z \to \bbP^1$ is a flat morphism such that $f|_Y=h$;
        \item the complement of $Y$ in $Z$ is a simple normal crossing anti-canonical divisor $D_Z=D^\mathsf{h} \cup D^\mathsf{v}$ where $D^\mathsf{h}$ and $D^\mathsf{v}=(f^{-1}(\infty))_{red}$ are horizontal and vertical divisors, respectively;
        \item $\vol_Z \in H^0(Z, K(*D_Z))$ is a nowhere vanishing meromorphic volume form with poles at most at $D_Z$.
    \end{enumerate}
    In particular, we call the compactified LG model $((Z, D_Z), f, \vol_Z)$ \textbf{tame} if it further satisfies the following conditions:
    \begin{enumerate}
        \item the critical locus $\crit(f)$ does not intersect with the horizontal divisors $D^h$;
        \item $f^{-1}(\infty)$ is reduced and the order of $\vol_Z$ at each irreducible component of $D$ is $-1$.
    \end{enumerate}
\end{defin}

\begin{defin}\label{def:fano type general}
    A LG model $(Y, \sw:Y \to \C)$ is of \textbf{Fano type} if it admits a tame compactification.
\end{defin}
We often simply write $((Z, D_Z), f:Z \to \bbP^1)$ for the compactified LG model since the choice of the volume form has no crucial role for the following discussion.  
Consider the logarithmic de Rham complex $(\Omega^\bullet_Z(\log D_Z),d)$. We define the subcomplex of $(\Omega^\bullet_Z(\log D_Z),d)$ which is preserved by $(-)\wedge df$, namely
	\begin{equation*}
	    \Omega^a_Z(\log D_Z, f):=\{u \in \Omega^a_Z(\log D_Z)| u \wedge df \in \Omega^{a+1}_Z(\log D_Z)\} 
	\end{equation*}
	for all $a \geq 0$. This subcomplex $(\Omega^\bullet_Z(\log D_Z, f), d)$ is called \textbf{$f$-adapted de Rham complex}. 
	
	\begin{example}
	Assume that we have a local chart $\C^2_{z_1, z_2}$ over $p \in D_Z$ and the morphism is given by $f(z_1, z_2)=\frac{1}{z_1z_2}$. Then we have the following local descriptions.
	\begin{equation*}
	    \begin{aligned}
	        \Omega^\bullet_Z(\log D_Z)& =\cO_Z \oplus \cO_Z\frac{dz_1}{z_1} \oplus \cO_Z\frac{dz_2}{z_2} \oplus \cO_Z\frac{dz_1dz_2}{z_1z_2} \\
	        \Omega^\bullet_Z(\log D_Z, f) &= (z_1z_2)\cO_Z \oplus (z_1z_2)\cO_Z\frac{dz_1}{z_1} \oplus (z_1z_2)\cO_Z\frac{dz_2}{z_2} \oplus \cO_Z\frac{dz_1dz_2}{z_1z_2}
	    \end{aligned}
	\end{equation*}
	\end{example}

    \begin{lem}\label{lem:local f-adapted smooth}
	The space of $f$-adapted de Rham differential forms of degree $a$,  $\Omega^a_Z(\log D_Z, f)$, is locally free of rank $\binom{n}{a}$ for all $a \geq 0$. Explicitly,
	\begin{equation*}
	\Omega_Z^a(\log D_Z, f)=\bigoplus_{p=0}^a\left[\frac{1}{f} \wedge^pW \oplus d\log f \wedge \left(\wedge^{p-1}W\right)\right]\bigotimes \wedge^{a-p}R
	\end{equation*}
	where $W$ is spanned by logarithmic 1-forms of the vertical part of $f:Z \to \bbP^1$ and $R$ is spanned by holomorphic 1-forms on $Y$ and logarithmic 1-forms associated to the horizontal part of $f:Z \to \bbP^1$.
	\end{lem}
	\begin{proof}
	See \cite[Lemma 2.12]{KKPbogomolov}
	\end{proof}
This local description allows to prove that Hodge-to-de Rham spectral sequence degenerates at the $E_1$-page. 
\begin{prop}
    The Hodge-to-de Rham spectral sequence for $(\Omega_Z^\bullet(\log D_Z, f), d)$ degenerates at the $E_1$-page. 
\end{prop}
\begin{proof}
    See \cite[Section 2]{KKPbogomolov} for the tame case, and \cite[Appendix C]{ESYIrregularHodge} for the general case. 
\end{proof}

\begin{defin}
    Let $(Y,\sw:Y \to \C)$ be a LG model with a compactification $((Z,D_Z),f)$. For non-negative integers $p,q \geq 0$, we define LG Hodge numbers 
    \begin{equation*}
        f^{p,q}(Y,\sw):=\dim_\C H^q(Z,\Omega_Z^p(\log D_Z,f))
    \end{equation*}
\end{defin}

Note that such LG Hodge numbers are independent of the choice of compactification \cite{YuIrregular, ESYIrregularHodge}

If we further assume that $((Z,D_Z),f)$ is tame, then the $f$-adapted de Rham complex $(\Omega^\bullet_Z(\log D_Z, f), d)$ is a limit of the relative de Rham complex $(\Omega^\bullet_Z(\log D_Z ,\text{ rel} f^{-1}(\rho)), d)$ as $\rho \to -\infty$. In particular, the Gauss-Manin parallel transport along $\rho \in \R_{<0}$ identifies $H^k(Y, Y_{sm};\C)$ with $\bbH^k(Z, (\Omega^\bullet_Z(\log D_Z, f), d))$ \cite[Claim 2.22]{KKPbogomolov}. Not only that, there is a
monodromy weight filtration $W_\bullet$ on $\bbH^k(Z, \Omega_Z^\bullet(\log D_Z, f))$ induced by the residue of the Gauss-Manin connection at infinity. On the other hand, one can consider the limiting mixed Hodge structure on $H^k(Y, Y_{sm};\Q)$ whose weight filtration is induced by the monodromy operator of a generic fiber $Y_{sm}$ around infinity \cite{ShamotoHodgeTate}. We denote it by $(H^k(Y, Y_{sm};\C),W^{\lim}_\bullet,  F_{\lim}^\bullet)$ (or $H^k(Y, Y_{sm, \infty})$ for short).
    
    \begin{prop}\label{prop: mixed hodge structure f-adapted smooth case}\cite[Section 3]{ShamotoHodgeTate}
        For $k \geq 0$, there exists a bi-filtered isomorphism  
        \begin{equation*}
            (\bbH^k(Z, \Omega_Z^\bullet(\log D_Z, f)), W_\bullet, F^\bullet) \cong (H^k(Y, Y_{sm};\C),W^{\lim}_\bullet,  F_{\lim}^\bullet)
        \end{equation*}
    \end{prop}
    
    \begin{proof}[Sketch of the proof]
   Consider the short exact sequence of complexes 
    \begin{equation}\label{eq:SES shamoto}
        0 \to \Omega_Z^\bullet(\log D_Z,f) \to \Omega_Z^\bullet(\log D_Z) \to \Omega_Z^\bullet(\log D_Z) \otimes \cO_{D_Z} \to 0
    \end{equation}
    Let $g:=1/f$ and consider a bi-filtered double complex $(A^{\bullet, \bullet}, d', d'', W^{\lim}_\bullet, F_{\lim}^\bullet)$ where 
    \begin{equation*}
        \begin{aligned}
        &A^{p,q}=\frac{\Omega^{p+q+1}_Z(\log D_Z)}{W_q(D_f)} \quad
        d'=d_{dR}:A^{p,q} \to A^{p+1, q}\quad 
        d''=g^{-1}dg:A^{p,q} \to A^{p,q+1} \\
        & W^{\lim}_rA^{p,q}=\frac{W_{r+q+1}(D_Z)\Omega^{p+q+1}_Z(\log D_Z)}{W_q(D_f)} \quad F^{\lim}_{p'}A^{p,q}=\frac{F^{p'}\Omega^{p+q+1}_Z(\log D_Z)}{W_q(D_f)}
        \end{aligned}
    \end{equation*} 
    Here $W_\bullet(D_Z)$\footnote{For future reference, we distinguish $D_f$ from $D_Z$, even though they represent the same divisor in the present case.} corresponds to Deligne's canonical weight filtration associated to the boundary divisor $D_Z$. The induced filtrations of $W^{\lim}_\bullet$ and $F_{\lim}^\bullet$ on the total complex $(\Tot(A^{\bullet, \bullet}), d=d'+d'')$ are also denoted the same.  Note that the complex $(\Tot(A^{\bullet, \bullet}), d)$ is quasi-isomorphic to $\Omega_Z^\bullet(\log D_Z) \otimes \cO_{D_Z}$ under the morphism
    \begin{equation*}
        \begin{aligned}
            \nu:\Omega^p_Z(\log D_Z) \otimes \cO_{D_Z}& \to A^{p,0} \\
            \tau & \mapsto [(-1)^pg^{-1}dg \wedge \tau \mod W_0]
        \end{aligned}
    \end{equation*}
     which induces a bi-filtered morphism 
    \begin{equation*}
        \nu:(\Omega_Z^\bullet(\log D_Z), W_\bullet, F^\bullet) \to (\Tot(A^{\bullet, \bullet}), W_\bullet^{\lim}, F^\bullet_{\lim})
    \end{equation*}
    whose kernel is $\Omega_Z^\bullet(\log D_Z, f)$ by construction. Therefore, there is a quasi-isomorphism $\Omega_Z^\bullet(\log D_Z, f) \to \cone(\nu)[-1]$ and it was shown that it induces the bi-filtered isomorphism on the cohomology level. 
    \end{proof}
    
    \begin{rem}
    In \cite{ShamotoHodgeTate}, the rational structure is constructed to promote the morphism $\nu$ to a morphism of $\Q$-mixed Hodge complexes\footnote{See Appendix \ref{sec:Appendix} for the definition of $\Q$-mixed Hodge complexes.}. Since the shifted cone of the morphism of $\Q$-mixed Hodge complexes is again a $\Q$-mixed Hodge complex, we have the limiting $\Q$-mixed Hodge structure on $H^a(Y, Y_{sm})$. Furthermore, we have a (strictly) filtered quasi-isomorphism $(\Omega_Z^\bullet(\log D_Z, f), F^\bullet) \to (\cone(\nu)[-1], F_{\cone}^\bullet)$ and the latter term is strict so that the spectral sequence of the Hodge filtration $F_\bullet$ on $\bbH^a(Z, \Omega_Z^\bullet(\log D_Z, f))$ degenerates at the $E_1$ page. 
    \end{rem}

Now we introduce another set of LG Hodge numbers. 
    \begin{defin}
        Let $(Y,\sw:Y \to \C)$ be a LG model of Fano type. For non-negative integers $p,q \geq 0$, we define LG Hodge numbers
        \begin{equation*}
        h^{p,q}(Y,\sw) :=\dim_\C \Gr_{2p}^{W^{\lim{}}}H^{p+q}(Y, Y_{sm,\infty})
        \end{equation*}
    \end{defin}
    
    We finish this section by introducing one of the main conjectures in \cite{KKPbogomolov} about the LG Hodge numbers. 
    
    \begin{conj}(\textbf{KKP conjecture})\cite[Conjecture 3.6]{KKPbogomolov}\label{conj:KKP Hodge numbers smooth} Let $(Y,\sw:Y \to \C)$ be a LG model of Fano type. Then the two LG Hodge numbers are the same: for non-negative integers $p,q \geq 0$,
        \begin{equation*}
            f^{p,q}(Y, \sw) = h^{p,q}(Y, \sw).
        \end{equation*}
    \end{conj}
    
\noindent Conjecture \ref{conj:KKP Hodge numbers smooth} holds if we impose the Hodge-Tate condition on $H^*(Y)$. We refer to Section \ref{sec:Extended Fano/LG} for a mirror symmetry motivation for this assumption.
    \begin{thm}\cite{ShamotoHodgeTate}\cite{Harderunpublished}\label{prop: Hodge-Tate implies KKP smooth setting}
        Let $(Y, \sw:Y \to \C)$ be a LG model of Fano type. Assume that Deligne's canonical mixed Hodge structure on $H^k(Y)$ is Hodge-Tate for all $k \geq 0$. Then the KKP conjecture holds. In other words, we have 
        \begin{equation*}
            f^{p,q}(Y, \sw) = h^{p,q}(Y, \sw)
        \end{equation*}
        for non-negative integers $p,q \geq 0$.
    \end{thm}
    
    \begin{proof}
        By \cite[Theorem 3.28]{ShamotoHodgeTate}, it is enough to show that the limiting mixed Hodge structure $H^k(Y, Y_{sm, \infty})$ is Hodge-Tate for all $k \geq 0$. Consider the long exact sequences of mixed Hodge structures
        \begin{equation*}
            \cdots \to H^{k-1}(Y_{sm, \infty}) \to H^k(Y, Y_{sm, \infty}) \to H^k(Y) \to H^k(Y_{Y_{sm}, \infty}) \to \cdots
        \end{equation*}
        Since the mixed Hodge structure on $H^k(Y)$ is Hodge-Tate, it is enough to show that $\Coker(H^{k-1}(Y) \to H^{k-1}(Y_{sm, \infty}))$ is Hodge-Tate. Consider the short exact sequence
        \begin{equation*}
            0 \to \Coker(H^{k-1}(Z) \to H^{k}(D_Z)) \to H^k_c(Y) \to \Ker(H^k(Z) \to H^k(D_Z)) \to 0
        \end{equation*}
        The Clemens-Schmid exact sequence implies that $\Coker(H^{k-1}(Z) \to H^{k}(D_Z))$ is the quotient of $\Ker(N_{k-1})$ by the image of the natural restriction map $\phi_Z:H^{k-1}(Z) \to H^{k-1}(Y_{sm, \infty})$. As the image of $\phi_Z$ coincides with the image of $\phi_Y:H^{k-1}(Y) \to H^{k-1}(Y_{sm, \infty})$, we argue that any non-zero class $\eta$ of type $(p,q)$ in $\Ker(N_{k-1})$ is in the image of $\phi_Y$ or Hodge-Tate (i.e. $p=q$). Now, take a non-zero class $\gamma \in H^{k-1}(Y_{sm, \infty})$ of type $(p,q)$. Then $N_{k-1}^i(\gamma) \neq 0$ while $N_{k-1}^{i+1}(\gamma)=0$ for some $i$. Also, the class $N_{k-1}^i(\gamma) \in \Ker(N_{k-1})$ is of type $(p-i, q-i)$ so that $p=q$ or it belongs to $\Image \phi_Y$. We finish the proof by showing that if $N_{k-1}^i(\gamma) \in \Image \phi_Y$, then $i=0$. Note that $\Image \phi_Y=\Image \phi_Z$ and $H^{k-1}(Z)$ admits the pure Hodge structure of weight $k-1$. If $i \neq 0$, then we get the contradiction to the fact that $N_{k-1}^i=0$.  
    \end{proof}

\subsection{The hybrid case}\label{sec:f-adpated general}
We extend the previous discussions to the hybrid setting.
\begin{defin}\label{def:compactified LG general}
	A compactified hybrid LG model of $(Y,h:Y \to \C^N)$ is a datum $((Z,D_Z), f:Z \to (\bbP^1)^N)$ where
		\begin{enumerate}
			\item $Z$ is a smooth projective variety and $f=(f_1, \dots, f_N) :Z \to (\bbP^1)^N$ is a projective morphism such that $f|_Y=h$;
			\item the complement of $Y$ in $Z$ is a simple normal crossing anti-canonical divisor $D_Z:=D_{f_1} \sqcup \cdots \sqcup D_{f_N}$ where
    	$D_{f_i}:=(f_i^{-1}(\infty))_{red}$ is the reduced pole divisor of $f_i$ for all $i=1, \dots, N$;
    	\item the morphism $f:Z \to (\bbP^1)^N$ is semi-stable at $(\infty, \dots, \infty)$. In other words, locally over infinity corner $(\infty, \dots, \infty)$, there are no reducible fibers except over the boundary $\mathrm{L}:=(\bbP^1)^N\setminus \C^N$.
		\end{enumerate}
		In particular, we call the compactified LG model $((Z, D_Z), f)$ \textbf{tame} if the pole divisor $f_i^{-1}(\infty)$ is reduced for all $i$. 
	\end{defin}
 
 \begin{defin}\label{def:Fanotype}
    A hybrid LG model $(Y, \omega, h:Y \to \C^N)$ is of \textbf{Fano type} if it admits a tame compactification.
\end{defin}

Note that each morphism $\Hat{f_i}=(f_1, \dots, \Hat{f_i}, \dots, f_N):Z \to (\bbP^1)^{N-1}$ is the compactification of the potential $\Hat{h_i}=(h_1, \dots, \hat{h_i}, \dots, h_N):Y \to \C^{N-1}$ for all $i=1, \dots, N$.

 \begin{rem}\label{rem:choice of volume form}
	\begin{enumerate}
	    \item The semi-stable condition follows from the local triviality of the hybrid LG potential $h:Y \to \C^N$ near infinity.
	    \item To align Definition \ref{def:compactified LG general} with Definition \ref{def:compactified LG classic}, we could incorporate the data of a volume form $\vol_Z$ and adjust the tameness condition to require that the order of $\vol_Z$ at each irreducible component of $D_Z$ is $-1$. However, as this information does not play a role in the subsequent discussion, we choose not to include it as a datum of the (tame) compactified LG model.
     
	\end{enumerate}
	\end{rem}
	\begin{example}(Example \ref{ex:conic+line} continued)\label{ex:conic+line cpt}
	   $(\bbP^2, C \cup L)$ Consider the Fano pair $(\bbP^2, C \cup L)$ where $C \cup L$ is the union of a smooth conic $C$ and a line $L$. Recall that the hybrid LG model is given by 
	   \begin{equation*}
		\begin{tikzcd}
		    Y:=\{x+y=az, z^2=bxy\} \subset \bbP^2_{x,y,z} \times \C_a \times \C_b \arrow[d] \\
		    \C_a \times \C_b
		\end{tikzcd}
		\end{equation*}
		First, compactify the ambient space to $\bbP^2 \times \bbP^1 \times \bbP^1$. Then we have the fibration $f:\bar{Y} \to \bbP^1 \times \bbP^1$ where $\bar{Y}$ is  the closure of $Y$. The boundary divisor is the wheel of 5 lines:
	   \begin{enumerate}
	       \item one vertical line $L_0$ over $(\infty, 0)$;
	       \item over $\{\infty\} \times \bbP^1$, there are two disjoint copies of lines $L_1$ and $L_2$ that intersect with $L_0$;
	       \item over $\bbP^1 \times \{\infty\} $, there are two copies of lines $L_3$ and $L_4$ that intersect with each other at one point over $(0, \infty)$;
	       \item $L_1$ (resp. $L_2$) intersects with $L_3$ (resp. $L_4$) at one point $p_{13}$ (resp. $p_{24}$) over $(\infty, \infty)$.
	   \end{enumerate}
	   However, $\bar{Y}$ is singular at the intersection points $p_{13}$ and $p_{24}$. We can desingularize $\bar{Y}$ to the compactifed LG model $f:Z \to \bbP^1 \times \bbP^1$. Then $D_Z$ is given by a wheel of 7 lines where the two points $p_{13}$ and $p_{24}$ are replaced by the exceptional divisors. We will see that it recovers the compactified LG model of the ordinary LG model $(Y, \sw:Y \to \C)$ (See Example \ref{ex: conic+line compact}).
	\end{example} 

	We fix the notations. For any non-empty subset $I=\{i_1, \dots, i_l\} \subset \{1, \dots, N\}$ and $1 \leq l \leq N$, let $h_I:=(h_{i_1}, \dots, h_{i_l}):Y \to \C^{|I|}$ be a (non-proper) hybrid LG model whose generic fiber is denoted by $Y_{I}:=h_I^{-1}(\rho_I)$ for some $\rho_I \in \C^{|I|}$. As before, we use the coordinates $(z_{i_1}, \dots, z_{i_l})$ for $\C^{|I|}$ where $(z_1, \dots, z_N)$ are the coordinates for $\C^N$. Also, we denote $\sqcup_{i \in I}Y_i$ the simple normal crossing of the generic fibers and denote $Y\{l\}=\sqcup_{|I|=l}Y_I$ the simple normal crossing union of $Y_{I}$'s with $|I|=l$ for $1\leq l\leq N$. Use the similar notations for the compactification $f:Z \to (\bbP^1)^N$.
	 
	For each index set $I \subset \{1, \dots N\}$, $1 
	\leq l\leq N$ and $a \geq 0$, we define subsheaves of the sheaf of the logarithmic forms $\Omega_Z^a(\log D_Z)$ as follows:
	\begin{equation*}
	\begin{aligned}
	    \Omega_Z^a(\log D_Z, \{f_i\}_{i \in I}):&=\{\eta \in \Omega_Z^a(\log D_Z)|df_i \wedge \eta \in \Omega_Z^{a+1}(\log D_Z)\text{ for all } i \in I\}\\
	    \Omega_Z^a(\log D_Z, f_I):&=\Span_{i\in I}(\{\eta \in \Omega_Z^a(\log D_Z)|df_i \wedge \eta \in \Omega_Z^{a+1}(\log D_Z)\})\\
	    \Omega_Z^a(\log D_Z, f_{(l)}):&=\Span_{|I|=l}\{\Omega_Z^a(\log D_Z, f_I)\} \\
	    & =\Span\{\eta \in \Omega_Z^a(\log D_Z)|df_i \wedge \eta \in \Omega_Z^{a+1}(\log D_Z) \text{ at least for $l$ } df_i's \}
	\end{aligned}
	\end{equation*}
\noindent Endowed with the de Rham differential, they form the subcomplexes of the logarithmic de Rham complexes, which are denoted by $(\Omega_Z^\bullet(\log D_Z, \{f_i\}_{i \in I}),d)$, $(\Omega_Z^\bullet(\log D_Z, f_I),d)$ and $(\Omega_Z^\bullet(\log D_Z, f_{(l)}),d)$, respectively. Due to the inclusion-exclusion principle, one can summarize the relations among those subcomplexes as follows.
	\begin{lem}\label{lem: Inclusion-Exclusion principle}
	For each $a \geq 0$, we have
	\begin{equation*}
	    \Span_{|I|=N+1-l}\{\Omega^a(\log D_Z, \{f_i\}_{i \in I})\}= \cap_{|I|=l}\Omega^a(\log D_Z, f_I) = \Omega_Z^a(\log D_Z, f_{(N+1-l)}).
	\end{equation*}
	\end{lem}

\begin{example}\label{ex:rank2case}(Rank $2$ case)
     For readability, we elaborate on the case of $N=2$. To a compactified LG model $((Z,D_Z), f:Z \to \bbP^1 \times \bbP^1)$, we can associate two different subcomplexes of the logarithmic de Rham complex $(\Omega^\bullet_Z(\log D_Z),d)$. First, consider the subcomplex $(\Omega^\bullet_Z(\log D_Z, f_1, f_2),d)$ which is preserved by the wedge product with both $df_1$ and $df_2$. This subcomplex is the pullback of the following diagram
	\begin{equation*}
	\begin{tikzcd}
	\Omega^\bullet_Z(\log D_Z, f_1, f_2) \arrow[r, dotted, "j_1"] \arrow[d, dotted, "j_2"] & \Omega^\bullet_Z(\log D_Z, f_1) \arrow[d, "i_1"] \\
	\Omega^\bullet_Z(\log D_Z, f_2) \arrow[r, "i_2"] & \Omega^\bullet_Z(\log D_Z) 
	\end{tikzcd}
	\end{equation*}
	where $i_1$ and $i_2$ are natural inclusions. By definition, it admits three differentials $d$, $\wedge df_1$ and $\wedge df_2$. Another subcomplex is the pushout of $j_1$ and $j_2$, denoted by $(\Omega^\bullet_Z(\log D_Z, f_{12}),d)$. Unlike the complex $(\Omega^\bullet_Z(\log D_Z, f_1, f_2),d)$, it only admits the de Rham differential.

	We perform a local computation of two de Rham subcomplexes $(\Omega^\bullet_Z(\log D_Z, f_1, f_2),d)$ and $(\Omega^\bullet_Z(\log D_Z, f_{12}), d)$. Recall that the complement of $Y$ in $Z$ is given by $D_Z=D_{f_1} \cup D_{f_2}$ where $D_{f_i}$ is the vertical boundary divisor of $f_i:Z \to \bbP^1$. Denote $D_{\infty}$ the intersection of $D_{f_1} \cap D_{f_2}$. For $p \in D_{\infty}$, we can find local analytic coordinates $z_1, \cdots, z_n$ centered at $p$ with $k_1+k_2\leq n$ such that in a neighborhood of $p$:
	\begin{itemize}
		\item the divisor $D_{f_1}$ is given by $\prod_{i=1}^{k_1}z_i=0$ and the potential $f_1$ is given by 
		\begin{equation*}
		f_1(z_1, \cdots, z_n)=\frac{1}{z_1^{m_1}\cdots z_{k_1}^{m_{k_1}}} 
		\end{equation*}
		for some $m_i \geq 1$. 
		\item the divisor $D_{f_2}$ is given by $\prod_{j=k_1+1}^{k_1+k_2}z_j=0$ and the potential $f_2$ is given by 
		\begin{equation*}
		f_2(z_1, \cdots, z_n)=\frac{1}{z_{k_1+1}^{m_{k_1+1}}\cdots z_{k_1+k_2}^{m_{k_1+k_2}}}
		\end{equation*}
		for some $m_j \geq 1$. 
	\end{itemize}
	By Lemma \ref{lem:local f-adapted smooth}, we know that the $f_i$-adapted de Rham complex $\Omega^a_Z(\log D_Z, f_i)$ is locally free of rank $\binom{n}{a}$ for all $a \geq 0$. Explicitly,
	\begin{equation*}
	\Omega_Z^a(\log D_Z, f_i)=\bigoplus_{p=0}^a\left[\frac{1}{f_i} \wedge^pW_i \oplus d\log f_i \wedge \left(\wedge^{p-1}W_i\right)\right]\bigotimes \wedge^{a-p}R_i
	\end{equation*}
	where $W_i$ is spanned by logarithmic 1-forms associated to the vertical part of $f_i:Z \to \bbP^1$ and $R_i$ is spanned by holomorphic 1-forms on $Y$ and logarithmic 1-forms associated to the horizontal part of $f_i:Z \to \bbP^1$. This local description allows one to describe $\Omega^a_Z(\log D_Z, f_1, f_2)$ for all $a \geq 0$. Explicitly, we have
	\begin{equation}\label{eq:local description}
	\begin{aligned}
	\Omega_Z^a(\log D_Z, f_1, f_2)=\bigoplus_{p+q=0}^a & \left[\frac{1}{f_1} \wedge^pW_1 \oplus d\log f_1 \wedge \left(\wedge^{p-1}W_1\right)\right] \bigotimes 
	\left[\frac{1}{f_2} \wedge^{q}W_2 \oplus d\log f_2 \wedge \left(\wedge^{q-1}W_2\right)\right] \\
	& \bigotimes \wedge^{a-p-q}R
	\end{aligned}
	\end{equation}
	where $R$ is spanned by holomorphic 1-forms on $Y$. Similarly, one can give  a local description of the complex $\Omega^\bullet_Z(\log D_Z, f_{12})$. Note that this is the subcomplex of the logarithmic de Rham complex $\Omega_Z^\bullet(\log D_Z)$ spanned by the forms preserved by either $df_1$ or $df_2$. Therefore, we have
	\begin{equation*}
	\begin{aligned}
	\Omega_Z^a(\log D_Z, f_{12})=\bigoplus_{p=0}^a & \left(\left[\frac{1}{f_1} \wedge^pW_1 \oplus d\log f_1 \wedge \left(\wedge^{p-1}W_1\right)\right]  +
	\left[\frac{1}{f_2} \wedge^{p}W_2 \oplus d\log f_2 \wedge \left(\wedge^{p-1}W_2\right)\right]\right) \\
	& \bigotimes \wedge^{a-p}R
	\end{aligned}
	\end{equation*}
\end{example}
    
\noindent In general, one can write down the local descriptions of $\Omega_Z^a(\log D_Z, \{f_i\}_{i \in I})$, $\Omega_Z^a(\log D_Z, f_I)$ and $\Omega_Z^a(\log D_Z, f_{(l)})$ similar as above, which leads to the following proposition.
	\begin{prop}\label{prop: degeneracy of Hodge filtration general case}
	The Hodge-to-de Rham spectral sequences for $(\Omega_Z^\bullet(\log D_Z, \{f_i\}_{i \in I}),d)$, $(\Omega_Z^\bullet(\log D_Z, f_I), d)$ and $(\Omega_Z^\bullet(\log D_Z, f_{(l)}),d)$ degenerate at the $E_1$-page. 
	\end{prop}
 \begin{proof}
 We first prove the rank 2 case. Recall that the strictness of Hodge filtrations on both  $R\Gamma(\Omega^\bullet_Z(\log D_Z, f_1), d)$ and $R\Gamma(\Omega^\bullet_Z(\log D_Z, f_2), d)$ was done in \cite{ESYIrregularHodge}(\cite{KKPbogomolov} for the tame case) by applying the method of Deligne-Illusie \cite{DeligneIllusiedef}. The main idea is to reduce the problem to a field of positive characteristic $p>0$, $\bbk$, and show the formality by constructing the global lifting of Frobenius morphism over $W_2(\bbk)$, the ring of Witt vectors of length 2 of $\bbk$. For the complex $(\Omega_Z^\bullet(\log D_Z, f_1, f_2),d)$, its local description (\ref{eq:local description}) allows one to apply the same argument used in \cite{ESYIrregularHodge}(\cite{KKPbogomolov} for the tame case). The only non-trivial part is to construct the gluing morphism between two choices of local liftings. However, since the horizontal part of boundary divisor with respect to $f_1$ is the vertical part of boundary divisor with respect to $f_2$ and vice versa, the choices of local liftings for the case of each $\Omega_{Z/\bbk}^\bullet(\log D_Z, f_i),d)$ are compatible. Therefore, the complex $R\Gamma(\Omega^\bullet_Z(\log D_Z, f_1, f_2), d)$ is strict. Also note that the pushout diagram 
		\begin{equation*}
		\begin{tikzcd}
		\Omega^\bullet_Z(\log D_Z, f_1, f_2) \arrow[r,  "j_1"] \arrow[d, "j_2"] & \Omega^\bullet_Z(\log D_Z, f_1) \arrow[d, dotted, "i_1"] \\
		\Omega^\bullet_Z(\log D_Z, f_2) \arrow[r, dotted, "i_2"] & \Omega^\bullet_Z(\log D_Z,  f_{12}) 
		\end{tikzcd}
		\end{equation*}
		induces a short exact sequence of the filtered complexes
		\begin{equation*}
		\begin{aligned}
		    0 \to (R\Gamma(\Omega^\bullet_Z(\log D_Z, f_1, f_2), d),F^\bullet) & \xrightarrow{(j_1, j_2)} \bigoplus_{i=1}^2(R\Gamma(\Omega^\bullet_Z(\log D_Z, f_i), d), F^\bullet) \\
		    & \xrightarrow{i_1-i_2} (R\Gamma(\Omega^\bullet_Z(\log D_Z, f_{12}), d), F^\bullet) \to 0
		\end{aligned}
		\end{equation*}
		Since the first two terms are strict and both $(j_1, j_2)$ and $H^*((j_1, j_2))$ are strict morphisms, we conclude that the cone complex $(R\Gamma(\Omega^\bullet_Z(\log D_Z, f_{12}), d), F^\bullet)$ is strict as well.
  
  The same argument can be iteratively applied to the higher rank case due to Lemma \ref{lem: Inclusion-Exclusion principle}. 
	\end{proof}
 
\begin{defin}
	Let $((Z,D_Z),f:Z \to (\bbP^1)^N)$ be a compactified hybrid LG model of $(Y, h:Y \to \C^N)$. For any index set $I \subset \{1, \dots, N\}$ and $1\leq l \leq N$, we define LG Hodge numbers as 
\begin{equation*}
	      \begin{aligned}
	            & f^{p,q}(Y, \{h_i\}_{i \in I}):=\dim_\C\bbH^q(Z, \Omega^p_Z(\log D_Z, \{f_i\}_I) \\
	           & f^{p,q}(Y,h_I):=\dim_\C \bbH^q(Z, \Omega_Z^p(\log D_Z, f_I)) \\
	           & f^{p,q}(Y,h_{(l)}):=\dim_\C \bbH^q(Z, \Omega_Z^p(\log D_Z, f_{(l)}))
	      \end{aligned}
	      \end{equation*}
for non-negative integers $p,q \geq 0$.
\end{defin}
\noindent Note that such LG Hodge numbers are independent of the choice of compactification \cite{YuIrregular, ESYIrregularHodge}

From now one, we assume that the compactified hybrid LG model $((Z,D_Z),f:Z \to (\bbP^1)^N)$ is tame. Then we can describe the limiting behavior of the complexes $(\Omega_Z^\bullet(\log D_Z, \{f_i\}_{i \in I}),d)$, $(\Omega_Z^\bullet(\log D_Z, f_I), d)$ and $(\Omega_Z^\bullet(\log D_Z, f_{(l)}),d)$.
	
	\begin{prop}\label{prop: limit of fadpated complex general case}
	\begin{enumerate}
	    \item The complex $(\Omega^\bullet_Z(\log D_Z, \{f_i\}_{i \in I}),d)$ is a well-defined limit of the relative de Rham complex $(\Omega^\bullet_Z(\log D_Z, \mathrm{rel}(\sqcup_{i \in I}f_i^{-1}(\rho_i)), d)$ as $\rho_i \to -\infty$ for $i \in I$. In particular, the Gauss-Manin parallel transport has a well defined limit as $\rho_I \to (-\infty, \dots ,-\infty)$ which identifies \\
	    $H^k(Y, \sqcup_{i \in I}h_i^{-1}(\rho);\C)$ with $\bbH^k(Z, \Omega^\bullet_Z(\log D_Z, \{f_i\}_{i \in I}))$ for $k \geq 0$.
	    \item The complex $(\Omega^\bullet_Z(\log D_Z, f_I),d)$ is a well-defined limit of the relative de Rham complex \\ $(\Omega^\bullet_Z(\log D_Z,  \mathrm{rel}(f_I^{-1}(\rho_I)), d)$ as $\rho_I \to (-\infty, \dots ,-\infty)$. In particular, the Gauss-Manin parallel transport has a well defined limit as  $\rho_I \to (-\infty, \dots ,-\infty)$ which identifies $H^k(Y, h^{-1}_I(\rho_I);\C)$ with $\bbH^k(Z, \Omega^\bullet_Z(\log D_Z, f_I))$ for $k \geq 0$.
	    \item The complex $(\Omega^\bullet_Z(\log D_Z, f_{(l)}),d)$ is a well-defined limit of the relative de Rham complex \\
	    $(\Omega^\bullet_Z(\log D_Z, \mathrm{rel}(\sqcup_{|I|=N+1-l}f_I^{-1}(\rho_I)), d)$ as $\rho_I \to (-\infty, \dots ,-\infty)$. In particular, the Gauss-Manin parallel transport has a well defined limit as  $\rho_I \to (-\infty, \dots ,-\infty)$ which identifies \\
	    $H^k(Y, \sqcup_{|I|=N+1-l}h_I^{-1}(\rho_I);\C)$ with $\bbH^k(Z, \Omega^\bullet_Z(\log D_Z, f_{(l)}))$ for $k \geq 0$. 
	\end{enumerate}
	\end{prop}
\noindent By Lemma \ref{lem: Inclusion-Exclusion principle}, the cohomology $\bbH^k(Z, \Omega^\bullet_Z(\log D_Z, f_{(l)}))$ can be considered as a limit of $H^k(Y, \cap_{|I|=l}h_I^{-1}(\rho_I);\C)$.
	\begin{proof}[Proof of Proposition \ref{prop: limit of fadpated complex general case}]
 Take a small polydisk centered at infinity $(\infty, \dots, \infty) \in (\bbP^1)^N$. Let $p:\mathcal{Z}:=Z \times \Delta^N \to \Delta^N$ be a proper family and $D_\mathcal{Z}:=D_Z \times \Delta^N$. For each case, we define the relative divisors
	\begin{enumerate}
	    \item $\Gamma:=\{(z, \rho_I) \in Z \times \Delta^N|\prod_{i \in I}(f_i(z)-\rho_i)=0\}$
	    \item $\Gamma:=(p\times f_I)^{-1}(\text{graph}( \Delta^N \hookrightarrow (\bbP^1)^N))$
	    \item $\Gamma:=\{(z, \rho) \in Z \times \Delta^N|\prod_{|I|=N+1-l}(f_I(z)-\rho_I)=0\}$
	\end{enumerate}
	They all intersect $D_\mathcal{Z}$ transversally and let $D_\Gamma=D_\mathcal{Z}\cap \Gamma$. We can construct the subsheaf of the relative meromorphic differential forms which vanishes along $D_\Gamma$ as follows:
	\begin{equation*}
	       \Omega^a_{\mathcal{Z}/\Delta^N}(\log D_\mathcal{Z}, \text{ rel} \mathbf{\Gamma}):=\ker(\Omega^a_{\mathcal{Z}/\Delta^N}(\log D_\mathcal{Z}) \to i_{\mathbf{\Gamma}*}\Omega^a_{\mathbf{\Gamma}/\Delta^N}(\log D_\mathbf{\Gamma})).
	\end{equation*}
	This is locally free and the complex $(\Omega^\bullet_{\mathcal{Z}/\Delta^N}(\log D_\mathcal{Z}, \text{ rel } \mathbf{\Gamma}),d)$ is preserved under the relative differential. Let's denote this complex by $\sE^\bullet_{\mathcal{Z}/\Delta^N}$. From the local computation, we have 
    \begin{enumerate}
	       \item \begin{equation*}
	       \begin{aligned}
	       (\sE^\bullet_{\mathcal{Z}/\Delta^N}|)_{Z \times \{\mathbf{\epsilon}\notin L\}}& =(\Omega^\bullet_Z(\log D_Z, \sqcup_{i \in I}f_i^{-1}(\rho_i), d)) \\
	       (\sE^\bullet_{\mathcal{Z}/\Delta^N}|)_{Z \times \{\mathbf{\epsilon}= 0\}}& =(\Omega^\bullet_Z(\log D_Z, \{f_i\}_{i \in I}),d)
	       \end{aligned}
	   \end{equation*}
	        \item \begin{equation*}
	       \begin{aligned}
	       (\sE^\bullet_{\mathcal{Z}/\Delta^N}|)_{Z \times \{\mathbf{\epsilon}\notin L\}}& =(\Omega^\bullet_Z(\log D_Z,  f_I^{-1}(\rho_I), d) \\
	       (\sE^\bullet_{\mathcal{Z}/\Delta^N}|)_{Z \times \{\mathbf{\epsilon}= 0\}}& =(\Omega^\bullet_Z(\log D_Z, f_I),d)
	       \end{aligned}
	   \end{equation*}
	        \item \begin{equation*}
	       \begin{aligned}
	       (\sE^\bullet_{\mathcal{Z}/\Delta^N}|)_{Z \times \{\mathbf{\epsilon}\notin L\}}& =(\Omega^\bullet_Z(\log D_Z, \sqcup_{|I|=N+1-l}f_I^{-1}(p), d) \\
	       (\sE^\bullet_{\mathcal{Z}/\Delta^N}|)_{Z \times \{\mathbf{\epsilon}= 0\}}& =(\Omega^\bullet_Z(\log D_Z, f_{(l)}),d)
	       \end{aligned}
	   \end{equation*}
	\end{enumerate}
	   where $L:=\Delta^N \setminus (\Delta^*)^N$. In order to finish the proof, we need to construct the Gauss-Manin connection on $R^ap_*\sE^\bullet_{\mathcal{Z}/\Delta^N}$
    \begin{equation*}
	       \nabla^{\text{GM}}: R^ap_*\sE^\bullet_{\mathcal{Z}/\Delta^N} \to R^ap_*\sE^\bullet_{\mathcal{Z}/\Delta^N} \otimes_{\cO_{\Delta^N}}\Omega^1_{\Delta^N}(\log L)
	   \end{equation*}
    Equivalently, we construct the subcomplex $\sE^\bullet_\mathcal{Z}$ of $(\Omega^\bullet_\mathcal{Z}(\log D_\mathcal{Z}),d))$ which fits into the short exact sequence of complexes
     \begin{equation}\label{SES:GM connection general}
	       0 \to \sE^\bullet_{\mathcal{Z}/{\Delta}^N}[-1] \otimes_{p^{-1}\cO_{\Delta^N}}p^{-1}\Omega_{{\Delta}^N}^1(\log L)  \to \sE^\bullet_\mathcal{Z} \to \sE^\bullet_{\mathcal{Z}/{\Delta}^N} \to 0
	   \end{equation}
	 since the connection $\nabla^{\text{GM}}$ is the connecting homomorphism in the long exact sequence of hyperderived
direct images associated with the short exact sequence (\ref{SES:GM connection general}). Define the subcomplex $\sE^\bullet_\mathcal{Z}$ to be
	   \begin{equation*}
 	        \ker\left(\frac{ \Omega^\bullet_\mathcal{Z}(\log (D_\mathcal{Z} \sqcup p^{-1}(L)))}{p^{-1}(\Omega^2_{{(\Delta^\ast)}^N}(\log L))\wedge \Omega^{\bullet-2}_\mathcal{Z}(\log (D_\mathcal{Z} \sqcup p^{-1}(L)))} \to i_{\mathbf{\Gamma}*}\frac{\Omega^\bullet_{\mathbf{\Gamma}}(\log D_\mathbf{\Gamma})}{p^{-1}(\Omega^2_{{(\Delta^\ast)}^N}(\log L)) \wedge\Omega^{\bullet-2}_{\mathbf{\Gamma}}(\log D_\mathbf{\Gamma}) }\right).
	   \end{equation*}
	   The above morphism is well-defined since the canonical morphism $\Omega^\bullet_\mathcal{Z}(\log (D_\mathcal{Z} \sqcup p^{-1}(L))) \to i_{\mathbf{\Gamma}*}\Omega^\bullet_{\mathbf{\Gamma}}$ is compatible with Koszul filtration. Then the short exact sequence is defined by the kernel of the following surjective morphism of short exact sequences (\ref{SES:GM connection general})
	\begin{equation*}
	       \begin{tikzcd}
	       0 \arrow[r] \arrow[d] & 0 \arrow[d] \\
	       \Omega^\bullet_{\mathcal{Z}/\Delta^N}(\log D_\mathcal{Z})[-1] \otimes p^{-1}\Omega^1_{\Delta^N}(\log L) \arrow[r] \arrow[d] & i_{\mathbf{\Gamma}*}\Omega^\bullet_{\mathbf{\Gamma}/\Delta^N}(\log D_\mathbf{\Gamma})[-1] \otimes p^{-1}\Omega^1_{\Delta^N}(\log L) \arrow[d]\\
	       \frac{\Omega^\bullet_\mathcal{Z}(\log D_\mathcal{Z} \sqcup p^{-1}(L))}{p^{-1}(\Omega^2_{\Delta^N}(\log L))\wedge \Omega^{\bullet-2}_\mathcal{Z}(\log D_\mathcal{Z} \sqcup p^{-1}(L))} \arrow[r] \arrow[d]
	        & i_{\mathbf{\Gamma}*}\frac{\Omega^\bullet_{\mathbf{\Gamma}}(\log D_\mathbf{\Gamma})}{p^{-1}(\Omega^2_{\Delta^N}(\log L)) \wedge\Omega^{\bullet-2}_{\mathbf{\Gamma}}(\log D_\mathbf{\Gamma}) } \arrow[d]\\
	    \Omega^\bullet_{\mathcal{Z}/\Delta^N}(\log D_\mathcal{Z}) \arrow[r] \arrow[d] &
	     i_{\mathbf{\Gamma}*}\Omega^\bullet_{\mathbf{\Gamma}/\Delta^N}(\log D_\mathbf{\Gamma}) \arrow[d] \\
	         0 \arrow[r] & 0
	       \end{tikzcd}
	   \end{equation*}
	\end{proof}
 \begin{example}(Rank $2$ case) We elaborate on this limiting behavior in the rank 2 case. Consider $(\Omega^\bullet_Z(\log D_Z, f_1, f_2), d)$ and $(\Omega^\bullet_Z(\log D_Z, f_{12}), d)$ as in Example \ref{ex:rank2case}. We begin with the complex \\ $(\Omega^\bullet_Z(\log D_Z, f_{12}), d)$. Let $\epsilon=(\epsilon_1, \epsilon_2)$ be a point near infinity and $Y_\epsilon$ be $f^{-1}(\epsilon)$. Then we have
	\begin{equation*}
	\begin{aligned}
	\Omega_Z^a(\log D_Z, \text{rel } Y_\epsilon):
	&=\ker(\Omega_Z^a(\log D_Z) \to i_{Y_\epsilon}\Omega^a_{Y_\epsilon}) \\
	&=\bigoplus_{p=0}^a\left[(\frac{1}{f_1}-\epsilon_1) \wedge^pW_1 + d\log f_1 \wedge \left(\wedge^{p-1}W_1\right)\right]\bigotimes \wedge^{a-p}R_1\\
	&+\left[(\frac{1}{f_2}-\epsilon_2) \wedge^pW_2 + d\log f_2 \wedge \left(\wedge^{p-1}W_2\right)\right]\bigotimes \wedge^{a-p}R_2 \\
	\end{aligned}
	\end{equation*}
	so when $\epsilon \to (0,0)$, this sheaf specializes to $\Omega^a_Z(\log D_Z, f)$. On the other hand, to describe $\Omega^a(\log D_Z, f_1, f_2)$ as a limit of a certain relative de Rham complex, we should consider the simple normal crossing union $Y_{1,\epsilon} \sqcup Y_{2, \epsilon}$ where $Y_{i,\epsilon}=f_i^{-1}(\epsilon_i)$. We have
	\begin{equation*}
	\begin{aligned}
	\Omega_Z^a(\log D_Z, \text{ rel} (Y_{1,\epsilon}\sqcup Y_{2, \epsilon})):
	&=\ker(\Omega_Z^a(\log D_Z) \to i_{Y_{1,\epsilon}}\Omega^a_{Y_{1,\epsilon}}\oplus i_{Y_{1,\epsilon}}\Omega^a_{Y_{1,\epsilon}}) \\
	&=\bigoplus_{p=0}^a\left[(\frac{1}{f_1}-\epsilon_1) \wedge^pW_1 + d\log f_1 \wedge \left(\wedge^{p-1}W_1\right)\right]\\
	&\bigotimes \left[(\frac{1}{f_2}-\epsilon_2) \wedge^pW_2 + d\log f_2 \wedge \left(\wedge^{p-1}W_2\right)\right]\bigotimes \wedge^{a-p}R \\
	\end{aligned}
	\end{equation*}
	so when $\epsilon \to (0,0)$, it specializes to $\Omega^a_Z(\log D_Z, f_1, f_2)$. Here, we use the assumption that $f:Z \to \bbP^1 \times \bbP^1$ is semi-stable at $(\infty, \infty) \in \bbP^1 \times \bbP^1$.
 \end{example}
\begin{example}
	We summarize the cohomology of several different complexes when $N=3$:
	\begin{equation*}
	    \begin{aligned}
	        &\bbH^k(Z,\Omega_Z^\bullet(\log D_Z, \{f_1, f_2, f_3\}))=\bbH^k(Z,\Omega_Z^\bullet(\log D_Z, f_{(3)})) \cong  H^k(Y, Y\{1\})=H^k(Y, Y_1+Y_2+Y_3)\\
	        &  \bbH^k(Z,\Omega_Z^\bullet(\log D_Z, \{f_i, f_j\}) \cong H^k(Y, Y_i+Y_j), \quad \bbH^k(\Omega_Z^\bullet(\log D_Z, f_{\{ij\}})) \cong H^k(Y, Y_{ij}) \\
	        & \bbH^k(Z,\Omega_Z^\bullet(\log D_Z, f_{(2)})) \cong H^k(Y, Y\{2\}) \cong H^k(Y, Y_{12}+Y_{13}+Y_{23}) \\
	        & \bbH^k(Z,\Omega_Z^\bullet(\log D_Z, f_{\{i\}})) \cong H^k(Y, Y_i) \\
	        & \bbH^k(Z,\Omega_Z^\bullet(\log D_Z, f_{\{123\}}))\cong \bbH^k(Z,\Omega_Z^\bullet(\log D_Z, f_{(1)})) \cong H^k(Y, Y\{3\}) = H^k(Y, Y_{123})
	    \end{aligned}
	\end{equation*}
	\end{example}
    
To state the extended version of Conjecture \ref{conj:KKP Hodge numbers smooth}, we introduce another set of Hodge LG numbers coming from the monodromies associated to $h:Y \to \C^N$. We provide the de Rham theoretic description by using simplicial methods (See Appendix \ref{Appendix MHC}).
    
    Fix $I \subset \{1, \dots, N\}$. For all $i \in I$, we write the intersection of the boundary divisors $D_{f_i}$'s as $D_{\infty, I}$ and set $g_i=1/f_i$.  For $p \in \Z_{\geq 0}$, $q_I=(q_1,\cdots, q_{|I|}) \in \Z^{|I|}_{\geq 0}$, we construct a multi-filtered complex $((A_{f_I}^{\bullet, \underline{\bullet}}),L_\bullet, F^\bullet, d_0, \{d_i\}_{i\in I})$ as follows:
        \begin{equation*}
        \begin{aligned}
         & A_{f_I}^{p,\underline{q}}:=\frac{\Omega_Z^{p+|q|+|I|}(\log D_Z)}{\sum_{i \in I}W_{f_i}(D_{f_i})}, \quad d_0=d_{dR}:A_{f_I}^{p,\underline{q}} \to A_{f_I}^{p+1, \underline{q}},  \quad d_i=g_i^{-1}dg_i:A_{f_I}^{p,\underline{q}} \to A_{f_I}^{p,\underline{q}+e_i} \\
        & W_r^{f_I}A_{f_I}^{p,\underline{q}}=W(D_Z)_{r+2|q|+|I|}\frac{\Omega_Z^{p+|q|+|I|}(\log D_Z)}{\sum_{i \in I}W_{f_i}(D_{f_i})}, \quad F^{p'}A_{f_I}^{p,\underline{q}}=\frac{F^{p'}\Omega_Z^{p+|q|+|I|}(\log D_Z)}{\sum_{i \in I}W_{f_i}(D_{f_i})} \\
         \end{aligned}
        \end{equation*}
        where $W_\bullet(-)$ is Deligne's canonical weight filtration associated to the divisor $(-)$. The induced filtrations of $W^{f_I}_\bullet$ and $F_{f_I}^\bullet$ on the total complex $(\Tot(A_{f_I}^{\bullet, \underline{\bullet}}), d=d_0 + \sum_{i \in I} d_i)$ are also denoted by the same. Note that the complex $(\Tot(A_{f_I}^{\bullet, \underline{\bullet}}), d)$ is quasi-isomorphic to $\Omega_Z^\bullet(\log D_Z) \otimes \cO_{D_\infty}$ under the morphism
        \begin{equation*}
            \begin{aligned}
                \nu_I:\Omega_Z^p(\log D_Z) \otimes \cO_{D_\infty} & \to A_{f_I}^{p,\underline{0}} \\
                \tau & \mapsto [(-1)^p\bigwedge_{i \in I}g^{-1}_idg_i \wedge \tau \mod \sum_{i \in I}W_{f_i}(D_{f_i})]
            \end{aligned}
        \end{equation*}
        where $\bigwedge_{i \in I}$ respects the increasing order of the indexes. It also provides the multi-filtered morphism  
        \begin{equation*}
            \nu_I:(\Omega_Z^\bullet(\log D_Z), W_\bullet, F^\bullet) \to (\Tot(A_{f_I}^{\bullet, \underline{\bullet}}), W^{f_I}_\bullet, F_{f_I}^\bullet)
        \end{equation*}
        which induces the long exact sequence of bi-filtered complexes
         \begin{equation}\label{eq:LES Y_I}
            \cdots \to H^k(Y, Y_{I,\infty}) \to H^k(Y) \to H^k(Y_{I,\infty}) \to H^{k+1}(Y,Y_{I, \infty}) \to \cdots 
        \end{equation}
         In particular, the kernel of $\nu$  is $(\Omega_Z^\bullet(\log D_Z, f_{I}),d)$ so that we have a canonical quasi-isomorphism from $(\Omega_Z^\bullet(\log D_Z, f_{I}),d)$ to the shifted cone $\cone(\nu)[-1]$. By introducing rational structures following \cite{FujisawaseveralvarII}\cite{ShamotoHodgeTate}, one can promote the bi-filtered complexes to the $\Q$-mixed Hodge complexes. Therefore, the long exact sequence (\ref{eq:LES Y_I}) becomes the sequence of $\Q$-mixed Hodge structures and the induced filtrations $W^{f_I}_\bullet$ (resp. $F_{f_I}^\bullet$) on both $ H^k(Y, Y_{I, \infty})$ and $ H^k(Y_{I, \infty})$ becomes weight (resp. Hodge) filtration. We will write the weight (resp. Hodge) filtration on both $ H^k(Y, Y_{I, \infty})$ and $ H^k(Y_{I, \infty})$ by $W^{h_I}_\bullet$ (resp. $ F_{h_I}^\bullet$) as it is independent of the choice of the tame compactification $(Z, f:Z \to (\bbP^1)^N)$.
    
        On the other hand, the mixed Hodge structure on $H^\bullet(Y, \sqcup_{i \in I}Y_i)$ can be  described by using a simplicial method as follows: Fix $K \subset I$. There exists a morphism of bi-filtered complexes
        \begin{equation*}
             \nu_{KI}:((\Tot(A_{f_K}^{\bullet, \underline{\bullet}}), d_K), W_\bullet^{f_K}, F^\bullet_{f_K}) \to ((\Tot{(A_{f_I}^{\bullet, \underline{\bullet}})}, d_I), W_\bullet^{f_I}, F^\bullet_{f_I})
        \end{equation*}
        induced by 
        \begin{equation*}
        \begin{aligned}
            A_{f_K}^{p, \underline{q}} & \to A_{f_I}^{p, (\underline{q}, 0)} \\
            \tau & \mapsto (-1)^p\bigwedge_{i \in I \setminus K} g^{-1}_idg_i \wedge \tau 
        \end{aligned}
        \end{equation*}
        where $\bigwedge_{i \in I \setminus K}$ respects the increasing order of the indexes. Let's construct a bi-filtered complex $(\{A_{[\bullet]}^\bullet\}, W_\bullet, F^\bullet)$ over the constant simplicial variety $\{Z_{[\bullet]}\}$ of length $|I|$ such that 
        \begin{enumerate}
            \item $Z_{[p]}=\bigoplus_{i=1}^{\binom{|I|}{p+1}}Z$;
            \item $A_{[p]}=\underset{\substack{K \subset I \\ |K|=p+1}}{\bigoplus} (\Tot(A_{f_K}^{\bullet, \underline{\bullet}}), W^{f_K}, F_{f_K})$;
            \item For $J \subset K \subset I$, the structure morphisms are given by $\nu_{JK}$.
        \end{enumerate}
         Note that the collection of morphisms $\{\nu_J\}_{J \subset I}$, defines the morphisms of bi-filtered complexes over the constant simplicial variety $\{Z_{[\bullet]}\}$
         \begin{equation*}
             \nu_{[\bullet]}:(\Omega_{Z_{[\bullet]}}^\bullet(\log D_Z), W_\bullet, F^\bullet) \to (\{A_{[\bullet]}^\bullet
         \}, W_\bullet, F^\bullet).
         \end{equation*}
        By construction, the pushforward of the kernel of $\{\nu_{[\bullet]}\}$ under the canonical anchor map $\pi_Z:\{Z_\bullet\} \to Z$ is isomorphic to the bi-filtered complex $(\Omega_Z^\bullet(\log D_Z, \{f_i\}_I), W_\bullet, F^\bullet)$. Again, equipped with the rational structures following \cite{ShamotoHodgeTate, FujisawaseveralvarII}, we end up getting the long exact sequence of $\Q$-mixed Hodge structures
        \begin{equation*}
            \cdots \to H^k(Y, \sqcup_{i \in I}Y_{i,\infty}) \to H^k(Y) \to H^k(\sqcup_{i \in I}Y_{i,\infty}) \to H^{k+1}(Y,\sqcup_{i \in I}Y_{i,\infty}) \to \cdots 
        \end{equation*}
      We write the weight (resp. Hodge) filtration on both $H^k(Y, \sqcup_{i \in I}Y_{i,\infty})$ and $H^k(\sqcup_{i \in I}Y_{i,\infty})$ by $W^{\{h_i\}_I}_\bullet$ (resp. $F_{\{h_i\}_I}^\bullet$). Finally, one can take a similar procedure to obtain the $\Q$-mixed Hodge structures on \\
$H^k(Y, Y\{N-l+1\})$, which is denoted by $(H^k(Y, Y\{N-l+1\}_\infty), W^{(N-1+l)}_\bullet, F_{(N-l+1)}^\bullet)$. Now we have Hodge-theoretic analogue of Theorem \ref{thm: A-side cubical diagram coh}.
\begin{thm}\label{thm:mainhodge}
Let $(Y, h:Y \to \C^N)$ be a hybrid LG model of Fano type. 
    \begin{enumerate}
        \item For $-n \leq a \leq n$, there exists a $\Q$-limiting mixed Hodge structure on $H^*(Y,h)=H^*(Y, Y\{1\})$ satisfying the functorial property: there exists a cubical diagram of $\Q$-limiting mixed Hodge structures
        \begin{equation*}
        \begin{aligned}
         H^{a+n-N}(Y\{N\}) &\xrightarrow{\psi_N} \bigoplus_{|I|=N-1}H^{a+n-|I|}(Y_I, h_I) \xrightarrow{\psi_{N-1}} \cdots \\
         &\xrightarrow{\psi_2} \bigoplus_{|I|=1} H^{a+n-|I|}(Y_I, h_I) \xrightarrow{\psi_1} H^{a+n}(Y, h)
        \end{aligned}
    \end{equation*} 
    where the limiting weight filtration is induced by the global monodromy $T$. 
    \item By applying Mayer-Vietoris sign rule \eqref{eq:MV sign rule}, this diagram becomes the $E_1$-page of the spectral sequence for the perverse filtration $P_\bullet^h$ on $(H^*(Y), W_\bullet)$ where $W_\bullet$ is Deligne's canonical weight filtration.
    \end{enumerate}
\end{thm}

\begin{proof}

By Proposition \ref{prop:two-perverse-same}, we have the following description of the perverse filtration $P^h_\bullet$:
\begin{equation*}
    P^h_{k+l}H^k(Y)=\Image(\gamma_l:H^k(Y, Y\{l\}) \to H^k(Y))
\end{equation*}
for $l=0, \dots, N$. Also, we have constructed the mixed Hodge structure on $H^k(Y,Y\{l\})$ such that the morphism $\gamma_l$ is a morphism of mixed Hodge structures, where the one on $H^k(Y)$ is Deligne's canonical mixed Hodge structure. This completes the proof. 
\end{proof}

	 \begin{rem}\label{Rem:HodgetoDerham general new proof}
    In this case, the $E_1$-degeneracy property of the Hodge-to-de Rham spectral sequences on \\
    $(\Omega_Z^\bullet(\log D_Z, \{f_i\}_{i \in I}),d)$, $(\Omega_Z^\bullet(\log D_Z, f_I), d)$ and $(\Omega_Z^\bullet(\log D_Z, f_{(l)}),d)$ also follows from the fact that $H^k(Y, \sqcup_{i \in I}Y_{i, \infty})$, $H^k(Y, Y_{I, \infty})$, and $(Y, Y\{N-l+1\}_\infty)$ form the mixed Hodge structures (especially, the degeneration of the Hodge filtration), respectively. 
    \end{rem}
	 
	\begin{defin}
        Let $(Y, h:Y \to \C^N)$ be a hybrid LG model of Fano type. For each $I \in \{1, \dots, N\}$ and $1 \leq l \leq N$, we define the LG Hodge numbers as
        \begin{equation*}
            \begin{aligned}
                & h^{p,q}(Y, \{h_i\}_I):=\dim_\C\Gr^{W^{\{h_i\}_I}}_{2p}H^{p+q}(Y, \sqcup_{i \in I}Y_{i,\infty};\C) \\
                & h^{p,q}(Y, h_I):=\dim_\C \Gr^{W^{h_I}}_{2p}H^{p+q}(Y, Y_{I, \infty};\C) \\
                & h^{p,q}(Y, h_{(l)}):=\dim_\C \Gr^{W^{(N-l+1)}}_{2p}H^{p+q}(Y, Y\{N-l+1\}_\infty;\C)
            \end{aligned}
        \end{equation*}
        for non-negative integers $p,q \geq 0$.
    \end{defin}
	
    \begin{conj}(\textbf{Extended KKP conjecture})\label{conj:Extended KKP} 
        \begin{enumerate}
            \item Let $(Y, h:Y \to \C^N)$ be a hybrid LG model of Fano type. For any $I \subset \{1, \dots, N\}$ and $1 \leq l \leq N$, there are identifications of the LG Hodge numbers
        \begin{equation*}
            \begin{aligned}
                & f^{p,q}(Y, \{h_i\}_I)=h^{p,q}(Y, \{h_i\}_I) \\
                & f^{p,q}(Y, h_I)=h^{p,q}(Y, h_I) \\
                & f^{p,q}(Y, h_{(l)})=h^{p,q}(Y, h_{(l)})
            \end{aligned}
        \end{equation*}
        for non-negative integers $p,q \geq 0$;
            \item Let $(Y, \sw:Y \to \C)$ be the associated ordinary LG model of the hybrid LG model $(Y, h:Y \to \C^N)$. Then there are identifications of LG Hodge numbers 
        \begin{equation*}
            f^{p,q}(Y, \sw)=f^{p,q}(Y, h_{(N)})
        \end{equation*}
        for non-negative integers $p,q \geq 0$.
        \end{enumerate}
    \end{conj}

   We explain the motivation behind this conjecture. The first part establishes a connection between the Hodge-theoretic data (on the left-hand side) and those originating from the categorical perspective (on the right-hand side), as discussed in Section \ref{sec:extendedfano lg 5.2}. The second part provides a Hodge-theoretic enhancement of the gluing property (Proposition \ref{prop:Gluing Property in general}), elucidating in what sense Conjecture \ref{conj:Extended KKP} extends the original KKP conjecture (Conjecture \ref{conj:KKP Hodge numbers smooth}). We should also note that, while the Hodge numbers are expected to be the same, their mixed Hodge structures are not necessarily isomorphic. Next, we will prove both parts of Conjecture \ref{conj:Extended KKP} under an assumption inspired by certain predictions in mirror symmetry.

    \begin{thm}\label{prop: Hodge-Tate implies KKP hybrid setting:General}
        Let $(Y, h:Y \to \C^N)$ be a hybrid LG model of Fano type. Assume that Deligne's canonical mixed Hodge structures on $H^k(Y)$ and $H^k(Y_I)$ are Hodge-Tate for all $I \subset \{1, \dots, N\}$, $k \geq 0$. Then the first part of Conjecture \ref{conj:Extended KKP} holds. In other words, for $1 \leq l \leq N$ and non-negative integers $p,q \geq 0$, we have
        \begin{equation*}
        \begin{aligned}
            & f^{p,q}(Y, \{h_i\}_I)=h^{p,q}(Y, \{h_i\}_I) \\
            &f^{p,q}(Y, h_I) = h^{p,q}(Y, h_I) \\
            &f^{p,q}(Y, h_{(l)}) = h^{p,q}(Y, h_{(l)})
        \end{aligned}
        \end{equation*}
        
    \end{thm}
    \begin{proof}
    We first prove the rank 2 case and apply the inductive argument for the general case. For the rank 2 case, it is enough to show that the Hodge-Tate conditions on $H^k(Y)$, $H^k(Y_1)$ and $H^k(Y_2)$ imply that the mixed Hodge structures on $H^k(Y)$,  $H^k(Y_{12, \infty})$ and $H^k(Y, Y_{1,\infty} \sqcup Y_{2, \infty})$ are Hodge-Tate. Consider the monodromy $T_1$ and $T_2$ introduced in Section \ref{sec:hybridLG}. They induce two increasing filtrations of the sub-mixed Hodge structures $W_\bullet^{h_1}$ and $W_\bullet^{h_2}$ on $(H^k(Y_{12, \infty}), W_\bullet^h, F^\bullet_h)$, respectively. As the weight filtration $W_\bullet^h$ is relative to both $W_\bullet^{h_1}$ and $W_\bullet^{h_2}$, it follows that $(H^k(Y_{12, \infty}), W_\bullet^h, F^\bullet_h)$ is Hodge-Tate if and only if for either $i=1$ or $i=2$, the induced limiting mixed Hodge structures on $(\Gr^{W^{h_i}}_aH^k(Y_{12, \infty}), W_\bullet^h, F^\bullet_h)$ are Hodge-Tate for all $a \geq 0$. For simplicity, we consider the case when $i=1$. The associated graded pieces $\Gr_{F_h}^p\Gr^{W^{h_1}}_aH^k(Y_{12,\infty})$ can be computed by taking the limiting mixed Hodge structure associated to the morphism $\sw_1:=h_2|_{Y_1}:Y_1 \to \C$, denoted by $(H^k(Y_{12, \infty_1}), W_\bullet^{\sw_1}, F^\bullet_{\sw_1})$. In other words, we have an isomorphism of the associated graded pieces
    \begin{equation*}
        \Gr_{F_{h_1}}^p\Gr^{W^{h_1}}_aH^k(Y_{12,\infty}) \cong \Gr_{F_{\sw_1}}^p\Gr^{W^{\sw_1}}_aH^k(Y_{12, \infty_1})
    \end{equation*}
     By Theorem \ref{prop: Hodge-Tate implies KKP smooth setting}, the assumption that the mixed Hodge structure on $H^k(Y_1)$ is Hodge-Tate implies that the limiting mixed Hodge structure on $H^k(Y_{12,\infty_1})$ is Hodge-Tate. Therefore,  $\Gr_{F_{\sw_1}}^p\Gr^{W^{\sw_1}}_aH^k(Y_{12,\infty})$ becomes trivial for $a \neq 2p$. Furthermore, it follows that  $\Gr_{F_h}^p\Gr^{W^h}_q\Gr^{W^{h_1}}_{2p}H^k(Y_{12,\infty})$ is trivial for all $q \neq 2p$ because 
    $(\Gr^{W^{h_1}}_{2p}H^k(Y_{12, \infty}), W_\bullet^h, F^\bullet_h)$ forms the limiting Hodge structure and $W_\bullet^h$ is of type $(-1,-1)$. Note that one can apply the parallel argument when $i=2$ to deduce the same conclusion. 

     To show that $H^k(Y, Y_{1, \infty} \sqcup Y_{2, \infty})$ is Hodge-Tate, we first recall that the mixed Hodge structure on $H^k(Y_{1, \infty}\sqcup Y_{2, \infty})$ is constructed in a way that it induces the long exact sequence of mixed Hodge structures
        \begin{equation*}
            \cdots \to H^k(Y_{1, \infty} \sqcup Y_{2, \infty}) \to H^k(Y_{1, \infty})\oplus H^k(Y_{2, \infty}) \to H^k(Y_{12, \infty}) \to \cdots
        \end{equation*}
        For $i=1,2$, the mixed Hodge structure $(H^k(Y_{i, \infty}), W_\bullet^{h_i}, F_{h_i}^\bullet)$ associated to $h_i:Y \to \C$ is relative to Deligne's canonical weight filtration $W_\bullet$. Since the canonical mixed Hodge structure $(H^k(Y_i), W_\bullet, F^\bullet)$ is Hodge-Tate by assumption, we apply the same argument to conclude that $(H^k(Y_{i,\infty}), W_\bullet^{h_i}, F^\bullet_{h_i})$ is Hodge-Tate.

        For the general case, the same argument above can be applied. To do this, one needs to show that the co-simplicial mixed Hodge complex whose cohomology of each simplicial term is Hodge-Tate induces the mixed Hodge structure of Hodge-Tate type on the level of cohomology. This is proven in Corollary \ref{cor:MHC hodge-Tate}. We leave the details of remaining proof to the reader. 
    \end{proof} 

\begin{example}
    A hybrid LG model $(Y, h:Y \to \C^2)$ constructed in the Example \ref{ex:conic+line}, \ref{ex:conic+line cpt} satisfies the assumption in Theorem \ref{prop: Hodge-Tate implies KKP hybrid setting:General}. The Hodge-Tateness of $H^k(Y)$ follows from that the compactification $Z$ is Hodge-Tate since this is a rational elliptic surface and two horizontal divisors we remove are $\bbP^1$.  
\end{example}

\begin{lem}\label{lem:hybridFanomeans}
    Let $(Y,h:Y \to \C^N)$ be a hybrid LG model of Fano type. Then the associated ordinary LG model $(Y,\sw:Y \to \C)$ is of Fano type. 
\end{lem}

\begin{proof}
    Let $((Z, D_Z), f:Z \to (\bbP^1)^N )$ be a tame compactification. We define a projective manifold $\bar{Z}$ via the pullback of the following diagram:
\begin{equation*}
    \begin{tikzcd}
     & Z \arrow[d, "f"] && \Bar{Z} \arrow[ll, "\Bar{\pi}"] \arrow[d, "{f_\pi}"] \\
    \C^N \arrow[r, hook] \arrow[d, "\Sigma"]& \bbP^1_{[a_1:b_1]} \times \cdots \times \bbP^1_{[a_N:b_N]} \arrow[d, dotted, "\bar{\Sigma}"] && \mathcal{P}((\bbP^1)^N) \arrow[ll, "\pi"] \arrow[dll, "\Bar{\Sigma}"] \\
    \C \arrow[r, hook] & \bbP^1 & 
    \end{tikzcd}
\end{equation*}
where $\bar{\Sigma}= \left[\sum a_i\prod_{j \neq i} b_j:b_1\cdots b_N\right]$ is a rational pencil and $\mathcal{P}(\bbP^1 \times \bbP^1)$ is a successive blow-up of the base locus of the pencil. Note that the base locus is given by the boundary $\mathsf{L}:=(\bbP^1)^N \setminus \C^N$ over which $Z$ is reduced by assumption, and the blow-up depends on the choice of order of the irreducible component in the base locus. 
\end{proof}
    
\begin{example}\label{ex: conic+line compact}(Example \ref{ex:conic+line cpt} continued) 
Consider the hybrid LG model $(Y, h:Y \to \C)$ that is constructed as a mirror to the pair $(\bbP^2, C \cup L)$ where $C \cup L$ is the union of a smooth conic $C$ and a line $L$. In Example \ref{ex:conic+line cpt}, we have constructed the compactification $f:Z \to \bbP^1 \times \bbP^1$ whose boundary divisor is the wheel of 7 lines. Pulling back $Z$ under the blow up of $\bbP^1 \times \bbP^1$ at $p=(\infty, \infty)$, the boundary divisor $D_Z$ becomes the wheel of 9 lines. Note that we further blow up the space $\Bl_p(\bbP^1 \times \bbP^1)$ along the base locus, denoted by $q$, due to the multiplicity of the base locus. By pulling back again under this blow-up, we get a well-known compactification of the ordinary LG model $(Y, \sw:Y \to \C)$. Recall that the ordinary LG model $(Y, \sw:Y \to \C)$ is obtained by removing two horizontal divisors from the mirror of the pair $(\bbP^2, E)$ where $E$ is a smooth elliptic curve. In this case, these two divisors is the pre-image of the exceptional divisor of $\Bl_q\Bl_p(\bbP^1 \times \bbP^1) \to \Bl_p(\bbP^1 \times \bbP^1)$ under $f_\pi$. 
\end{example}

If we further assume that Deligne's mixed Hodge structures on $H^k(Y, Y_{sm})$ is Hodge-Tate, one can prove the second part of the conjecture. 
\begin{thm}\label{thm:extended KKP second}
    Let $(Y,h:Y \to \C^N)$ be a hybrid LG model of Fano type. Assume that Deligne's mixed Hodge structure on $H^k(Y, Y_{sm})$ is Hodge-Tate and the limiting mixed Hodge structure on $H^k(Y, h)=H^k(Y, Y\{N\}_\infty)$ is Hodge-Tate. Then the second part of Conjecture \ref{conj:Extended KKP} holds. In other words, for non-negative integers $p,q \geq 0$, we have
    \begin{equation*}
        f^{p,q}(Y,\sw)=f^{p,q}(Y, h_{(N)})
    \end{equation*}
    
\end{thm}
\begin{proof}
    First, observe that there exists an isomorphism $H^*(Y, Y_{sm}) \cong H^*(Y, h)$ that is compatible with the monodromy $T$. Given our assumption that the hybrid LG model $(Y, h)$ is of Fano type, this implies that its ordinary LG model $(Y, \sw)$ is of Fano type as well due to Lemma \ref{lem:hybridFanomeans}, and consequently, the monodromy $T$ induces the weight monodromy filtration $W^T_\bullet$. On $H^*(Y,h)$, such a weight filtration constitutes the limiting mixed Hodge structure, so it is enough to show that 
    \[
    \Gr_F^p\Gr^{W^T}_{q}H^*(Y,Y_{sm})=0
    \]
    except for $q=2p$. This follows from the same technique that we used in the proof of Theorem \ref{prop: Hodge-Tate implies KKP hybrid setting:General}. Namely, Deligne's weight filtration $W_\bullet$ is compatible with $W^T_\bullet$ so that the associated graded pieces $\Gr^WH^*(Y,Y_{sm})$ admits a limiting mixed Hodge structure with respect to $W_\bullet^T$. Since the monodromy operator is of type $(-1,-1)$ and Deligne's mixed Hodge structure is assumed to be Hodge-Tate, we achieve the claim. 
\end{proof}
    
    \subsubsection{The perverse filtration}
    Let $(Y, h:Y \to \C^N)$ be a hybrid LG model of Fano type and $((Z, D_Z),f:Z \to (\bbP^1)^N)$ be its tame compactification. For $a \geq 0$, we define an increasing filtration on the sheaf of the logarithmic de Rham complex $\Omega^a_Z(\log D_Z)$ as follows:
    \begin{equation*}
        G^f_l\Omega_Z^a(\log D_Z):=\begin{cases}
        0 & l < 0 \\
        \Omega_Z^a(\log D_Z) & l \geq N \\
        \Omega_Z^a(\log D_Z, f_{(N-l)}) & 0 \leq l \leq N
        \end{cases}
    \end{equation*}
    This filtration $G^f_\bullet$ induces the filtration in cohomology
    \begin{equation*}
        G^f_lH^k(Y;\C)=\Image(\bbH^k(Y, G^f_{l-k}\Omega_Z^\bullet(\log D_Z)) \to \bbH^k(Z, \Omega_Z^\bullet (\log D_Z)))
    \end{equation*}
    We have showed that the filtration $G_\bullet^f$ does not depend on the choice of the compactification. Furthermore, it is equivalent to the perverse filtration $P^h_\bullet$ associated to the affinitzation map $h:Y \to \C^N$ due to Proposition \ref{prop:two-perverse-same}.
    
    \begin{thm}\label{thm: two fil same general}
        For $k \geq 0$, the two filtrations $G^f_\bullet$ and $P^h_\bullet$ on $H^k(Y)$ are the same.
    \end{thm}
    \begin{thm}\label{thm:multiplicativity}
    For $k_1, k_2 \geq 0$, the perverse filtration $P_\bullet$ associated to a hybrid LG model $(Y, h:Y \to \C^N)$ satisfies the following relation under the cup product
    \begin{equation*}
        \cup: P_{k_1+a_1}H^{k_1}(Y;\C) \otimes P_{k_2+a_2}H^{k_2}(Y;\C) \to P_{k_1+k_2+\min(a_1, a_2)}H^{k_1+k_2}(Y;\C)
    \end{equation*}
    where $0 \leq a_1, a_2 \leq N$.
    \end{thm}

    \begin{proof}
    Let $\Delta_Z:Z \to Z\times Z$ be a diagonal morphism and $D_{Z\times Z}$ be the boundary divisor $D_Z \times Z \cup Z \times D_Z$. We consider the morphism $f \times f:Z \times Z \to (\bbP^1)^N \times (\bbP^1)^N$ and write $f_{i,j}$ for the $(i,j)$-th component of $f \times f$.
    Recall that $H^\bullet(Y;\C)=\bbH^\bullet(Z, \Omega_Z^\bullet(\log D_Z))$ and the cup product is induced by the composition of the logarithmic de Rham complexes
    \begin{equation*}
    \Omega_Z^\bullet(\log D_Z) \otimes \Omega_Z^\bullet(\log D_Z) \xrightarrow{kun} \Omega_{Z\times Z}^\bullet(\log D_{Z\times Z}) \xrightarrow{\Delta_Z^*} \Omega_Z^\bullet(\log D_Z)
    \end{equation*}
    where $kun$ is the Künneth morphism. First, we show that the morphism $kun$ is indeed the morphism of filtered complexes 
    \begin{equation*}
        kun:(\Omega_Z^\bullet(\log D_Z), G_\bullet^f) \otimes (\Omega_Z^\bullet(\log D_Z), G_\bullet^f) \to (\Omega_{Z\times Z}^\bullet(\log D_{Z\times Z}), G_\bullet^{f\times f})
    \end{equation*}
    If $\alpha \in G^f_a\Omega_Z^{r}(\log D_Z)$ and $\beta \in G^f_b\Omega_Z^{s}(\log D_Z)$ are simple, then $\alpha \wedge \beta$ is preserved by at least $N-a+N-b$ number of $df_{i,j}\wedge$'s, so in $G^{f\times f}_{a+b}\Omega_Z^{r+s}(\log D_{Z\times Z})$. Next, take a simple form $\gamma \in G^{f \times f}_c\Omega_{Z\times Z}^t(\log D_{Z \times Z})$. By construction, $\gamma$ is preserved by at least $2N-c$ number of $df_{i,j}\wedge$'s. We say $\gamma$ is preserved by $(N-c_1, N-c_2)$ $df_{i,j} \wedge$'s if it is in the image of $G^f_{c_1}\Omega^\bullet_Z(\log D_Z) \otimes G^f_{c_2}\Omega^\bullet_Z(\log D_Z)$ under the Künneth morphism. Now, it is straightforward to check that the image $\Delta_Z^*(\gamma)$ is preserved by at least $\max(N-c_1, N-c_2)$ and at most $\min(2N-c_1-c_2,N)$ $df_i$'s. 
    \end{proof}
    
    \begin{rem}
    Recall that the definition of the perverse filtration $P_\bullet^h$ does not involve the shift on the complex $Rh_*{\C_Y}$. Therefore, the multiplicativity conjecture of the perverse filtration $P^h_\bullet$ is the following form (See \cite[Section 1.4.1]{de2012topology}):
    \begin{equation*}
        \cup: P^h_{l_1}H^{k_1}(Y) \otimes P^h_{l_2}H^{k_2}(Y) \to P^h_{l_1+l_2-N}H^{k_1+k_2}(Y)
    \end{equation*}
    while the weaker compatibility is already known as follows:
    \begin{equation*}
        \cup: P^h_{l_1}H^{k_1}(Y) \otimes P^h_{l_2}H^{k_2}(Y) \to P^h_{l_1+l_2}H^{k_1+k_2}(Y).
    \end{equation*}
    
    \end{rem}
    \begin{cor}
    If $a_1=N$, then the cup product is compatible with the perverse filtration. In other words, we have
    \begin{equation*}
        \cup: P_{k_1+N}H^{k_1}(Y) \otimes P_{k_2+a_2}H^{k_2}(Y) \to P_{k_1+k_2+a_2}H^{k_1+k_2}(Y). 
    \end{equation*}
    \end{cor}

\section{Extended Fano/LG correspondence}\label{sec:Extended Fano/LG}
In this section, we discuss mirror symmetry between a Fano pair $(X,D)$ and its mirror hybrid LG model $(Y, \omega, h:Y \to \C^N)$, called the \textit{extended Fano/LG correspondence}. Specifically, our focus is on $B$-to-$A$ mirror symmetry, treating the Fano pair $(X,D)$ as a $B$-side object and the hybrid LG model $(Y, \omega, h:Y \to \C^N)$ as an $A$-side object. A primary objective of this section is to propose a relative mirror symmetry conjecture ($B$-to-$A$), which is expected to recover the mirror P=W conjecture. We begin with an overview of the classical Fano/LG correspondence and subsequently extend our discussion to the hybrid case.

\subsection{Classical Fano/LG correspondence}\label{sec:Fano/LG}

We review a classical version of the Fano/LG correspondence. We refer more details to \cite{AurouxFanoSYZ}\cite{KKPHodgemirror}\cite{KKPbogomolov}.

\begin{defin}\label{def:ordinary LG model}
	A \textbf{symplectic Landau-Ginzburg (LG) model} is a triple 
	$(Y, \omega, \sw:Y \to \C)$ where
	\begin{enumerate}
		\item $(Y, \omega)$ is $n$-dimensional complex K\"ahler Calabi-Yau manifold with a K\"ahler form $\omega \in \Omega^2(Y)$;
		\item $\sw:Y \to \C$ is a (surjective) holomorphic map such that $\sw:Y \to \C$ is a locally trivial symplectic fibration near infinity with smooth fibers. In other words, there exists a compact subset $K \subset \C$ where $\sw:Y \to \C$ is locally trivial over $\C \setminus K$.
	\end{enumerate}
	We call $\sw:Y \to \C$ a \textbf{Landau-Ginzburg (LG) potential}.
\end{defin}

Let $(X,D)$ be a Fano pair where $X$ is a compact $n$-dimensional complex K\"ahler Fano manifold (or smooth projective variety) and $D$ is an anti-canonical divisor. We say that a LG model $(Y, \omega, \sw:Y \to \C)$ is mirror to the pair $(X,D)$ if it satisfies the following mirror relations:
\begin{equation*}
    \begin{aligned}
    & X  \longleftrightarrow  (\mathsf{w}:Y \to \C) \\
& D  \longleftrightarrow  Y_{sm} \\
& U \longleftrightarrow Y
    \end{aligned}
\end{equation*}
We write it as $(X,D)|(Y, \omega, \sw:Y \to \C)$ and refer to it as a Fano mirror pair. These abstract correspondences can be made explicit once the relevant mathematical objects for each mirror pair are specified.

\begin{rem}
\begin{enumerate}
    \item The properness of the LG potential $\sw:Y \to \C$ is tied to the choice of the anti-canonical divisor $D$. According to mirror symmetry predictions, if $D$ is smooth, the LG potential $\sw:Y \to \C$ is expected to be proper. If $D$ is not smooth, the LG potential $\sw:Y \to \C$ is no longer expected to be proper. This phenomenon of the non-properness of the LG potential is a primary motivation for the hybrid LG model.
    \item In \cite{KKPbogomolov}, the Fano/LG correspondence considers a choice of the defining anti-canonical section $s_X \in |K_X^{-1}|$ and a holomorphic volume form on $Y$, $\vol_Y$. Since these choices are not crucial for our current discussion, we omit them in our definition.
\end{enumerate}
\end{rem}

\subsubsection{B-side objects}\label{sec: Classical B-side}
Let $(X,D)$ be a Fano pair and $\iota:D \hookrightarrow X$ be the natural inclusion. We consider the derived pushforward functor between the (dg enhancement of) bounded derived category of coherent sheaves of $D$ and $X$
\begin{equation*}
    \iota_*:D^b\Coh(D) \to D^b\Coh(X)
\end{equation*}
Note that $\iota_*$ is compatible with the Serre functor in the sense that $\iota_* \circ S_X|_D=S_X \circ \iota_*$ where $S_X=(-)\otimes K_X[n]$ and $S_X|_D=(-)\otimes K_X|_D[n]$. We will also say that the functor $i_*$ is compatible with the pair of auto-equivalences $(S_X|_D, S_X)$. Moreover, the categorical localization \cite{Seidelnatrual08} of $\iota_*$ becomes equivalent to the derived category of the complement $U=X\setminus D$, $D^b\Coh(U)$. In case that $D$ is smooth, we apply Hochschild homology to the functor $\iota_*$ to obtain a cohomological object we should consider. In other words, since $HH_a(D^b\Coh(D))=\bigoplus_{p-q=a}H^{p,q}(D)$ and $HH_a(D^b\Coh(X))=\bigoplus_{p-q=a}H^{p,q}(X)$ for $-n \leq a \leq n$ due to the HKR isomorphism \cite{HKR}, we have 
\begin{equation}\label{eq:Gysin B-side Hocshchild}
    \iota_!:\bigoplus_{p-q=a}H^{p,q}(D) \to \bigoplus_{p-q=a}H^{p,q}(X)
\end{equation}
Furthermore, the pair $(S_X|_D, S_X)$ induces the automorphism of the map (\ref{eq:Gysin B-side Hocshchild}) whose logarithm is given by the cup product of $(c_1(K_X|_D), c_1(K_X))$. We simply denote this automorphism by $S$. Since $X$ is Fano, it induces the monodromy weight filtration centered at $n-|a|$, denoted by $W^S_\bullet$. Then the associated graded pieces are given by 
\begin{equation*}
    \begin{aligned}
        &\Gr^{W^S}_{2(n-p)}\iota_!:H^{p,q}(D) \to H^{p,q}(X) \quad \text{for } a \geq 0 \\
        &\Gr^{W^S}_{2(n-q)}\iota_!:H^{p,q}(D) \to H^{p,q}(X) \quad \text{for } a \leq 0
    \end{aligned}
\end{equation*}

\subsubsection{A-side objects}\label{sec: Classical A-side}
Let $(Y, \omega, \sw:Y \to \C)$ be a symplectic LG model. Due to the intricate nature of the general treatment of symplectic geometry for LG models, we will provide a brief overview, specifically focusing on the case where $(Y,\omega)$ is a Liouville manifold. Additional details can be found in \cite[Appendix A]{ASKhovanov}\cite{ASLefschetz}.

\begin{defin}\label{def:horizontal admissible}

An admissible Lagrangian associated to a symplectic fibration $\sw:Y \to \C$ is a (possibly non-compact) conical Lagrangian $L$ in $Y$ such that $\sw(L)$ is contained in a union of a compact subset and (possibly multiple) radial rays away from the negative real axis.  We say $L$ is \textit{horizontally admissible} if $L$ is proper over $\sw(L)$.
\end{defin}

To an admissible Lagrangian $L$, one can associate a subset $D_L \subset (-\pi, \pi)$ of directions of $L$ near $\infty$. In order to define Floer theory of admissible Lagrangians, we allow non-compact Hamiltonian perturbations as well as choose "counterclockwise" direction near $\infty$. For admissible Lagrangians $K$ and $L$, we say $K>L$ if $\theta_k>\theta_L$ for any $\theta_k \in D_K$ and $\theta_L \in D_L$. Consider $\mathfrak{A}$ be a directed $A_\infty$-category whose
\begin{itemize}
    \item Objects are horizontally admissible Lagrangian branes of $(Y, \omega, \sw)$;
    \item Morphism spaces are \begin{equation*}
     \Hom_\mathfrak{A}(K,L)=\begin{cases}
     CF^\bullet(K,L) & K>L \\
     \C<e^+_L> & K=L \\
     0 & otherwise \\
     \end{cases}
    \end{equation*}
    where $e^+_L$ is a strict unit.
\end{itemize}
This is directed in the sense that the set of objects of $\mathfrak{A}$, is a poset with $\Hom_\mathfrak{A}(K,L)=0$ unless $K>L$. We also define the $A_\infty$ structure in the usual way \cite[Appendix A.8]{ASKhovanov}. Moreover, in this category, $L$ is not quasi-isomorphic to its perturbation $\phi_\epsilon L$ where $\phi_\epsilon$ does not cross the negative real axis because $\Hom_\mathfrak{A}(\phi_\epsilon L, L) \neq 0$ while $\hom_\mathfrak{A}(L, \phi_\epsilon L)=0$. One achieves a quasi-isomorphism between $L$ and $\phi_\epsilon(L)$ by inverting all quasi-units in the cohomology of morphism spaces. Let $Z \subset H^0(\mathfrak{A})$ be a collection of all quasi-units. 

\begin{defin}\cite[Definition A.11]{ASKhovanov} \cite{ASLefschetz}
The Fukaya-Seidel category $\FS(Y,\sw)$ is defined to be a $A_\infty$-localization of $\mathfrak{A}$ at $Z$.
\end{defin}

\noindent We consider the cup functor from the Fukaya category $\Fuk(Y_{sm})$ to the Fukaya-Seidel category $\FS(Y,\sw)$ 
\begin{equation*}
    \cup:\Fuk(Y_{sm}) \to \FS(Y,\sw)
\end{equation*}
which is given by the trajectory of a parallel transport along a U-shaped curve enclosing all critical values. Consider the Serre functor $S_Y:=\phi_{2\pi}^{-1}[n]$ on $\Fuk(Y,\sw)$ where $\phi_{2\pi}$ is defined by wrapping once counter-clockwise. Let $\mu\in Aut(\Fuk(Y_{sm}))$ be the global monodromy induced by a parallel transport counter-clockwise along a large enough loop and we write $S_Y|_{Y_{sm}}=\mu^{-1}[n-1]$. Then the cup functor is compatible with the Serre functor in the sense that $\cup \circ S_Y|_{Y_{sm}} = S_Y \circ \cup$. We will also say that the functor $\cup$ is compatible with the pair of auto-equivalences $(S_Y|_{Y_{sm}}, S_Y)$. Moreover, the categorical localization \cite{Seidelnatrual08} of $\cup$ becomes equivalent to the wrapped Fukaya category of $Y$, $\mathcal{W}(Y)$. For later use, we denote $D^\pi(-)$ the split-closed derived category.

In case that $\sw:Y \to \C$ is a Lefschetz fibration, we can apply Hochschild homology to obtain a relevant cohomological object. As a vector space over $\C$, we have $HH_a(\Fuk(Y_{sm})) \cong H^{n-1+a}(Y_{sm})$ and $HH_a(\FS(Y, \sw)) \cong H^{n-1+a}(Y_{sm})$ for $-n \leq a \leq n$. The induced morphism $HH_a(\cup)$ becomes the connecting homomorphism of the long exact sequence of the pair $(Y, Y_{sm})$
\begin{equation*}
    \rho:H^{n-1+a}(Y_{sm}) \to H^{n+a}(Y, Y_{sm})
\end{equation*}
Moreover, the the pair $(\mu, \phi)$ induces the pair of operators given by taking the monodromy around infinity. We denote it by $T$ and the induced monodromy weight filtration by $W^T_\bullet$. Motivated by the case of the Lefschetz fibration, we will consider the connecting homomorphism $\rho$ with the filtration $W^T_\bullet$ as a cohomological object on the A-side. 

\begin{rem}\label{rem:exactness}
When dealing with situations where $(Y, \omega)$ is not exact, as considered in the main body, the theory of Fukaya-Seidel category has not been fully developed yet except in the case of Lefschetz fibration, as seen in the series of works of Seidel (e.g. \cite{SeidelLefV}). Moreover, it is worth noting that one should work over some Novikov field $\bbk$ not just $\C$. We leave this as a black box and only provide conjectural arguments regarding this topic. 
\end{rem}

\subsubsection{Mirror symmetry}\label{sec:Classical Fano/LG-mirror}
Based on the previous discussion, we introduce two versions of the mirror symmetry conjecture for a Fano mirror pair $(X,D)|(Y, \omega, \sw:Y \to \C)$. Here, we first present an ad-hoc version of the homological mirror symmetry (HMS) conjecture which is only valid for the case that $(Y, \omega)$ is a Liouville manifold.
\begin{conj}\label{conj:HMS classic}
There exists a commutative diagram of derived categories 
\begin{equation}\label{eq:rel HMS classic}
    \begin{tikzcd}
        D^b\Coh(D) \arrow[r,"\iota_*"] \arrow[d, "\cong"]& D^b\Coh(X) \arrow[d, "\cong"] \\
        D^\pi\Fuk(Y_{sm}) \arrow[r, "\cup"]& D^\pi\FS(Y,\sw)
    \end{tikzcd}
\end{equation}
which is compatible with auto-equivalences $(S_X|_D, S_X)$ and $(S_Y|_{Y_{sm}}, S_Y)$.
\end{conj}

\noindent As discussed before, the induced operations of the Serre functors $S_X$ and $S_Y$ on the Hochschild homology would correspond to the automorphisms $S$ and $T$, respectively. Since the automorphism $S$ unipotent, it is natural to expect that the monodromy operation $T$ is also unipotent. Geometrically, it corresponds to the existence of a tame compactification which we will discuss in Section \ref{sec:deform}. We refer the reader to \cite[Remark 2.5]{KKPbogomolov} for more details.

\begin{conj}\label{conj:CMS classic}
    For $-n \leq a \leq n$, there is an isomorphism of two sequences of cohomology groups
    \begin{equation}\label{eq:CMS classic}
    \begin{tikzcd}
        \bigoplus_{p-q=a} H^{p,q}(D) \arrow[r,"\iota_!"] \arrow[d,"\cong"]&  \bigoplus_{p-q=a} H^{p,q}(X) \arrow[d,"\cong"] \\
        H^{n-1+a}(Y_{sm})\arrow[r, "\rho"]&  H^{n+a}(Y,Y_{sm})
    \end{tikzcd}
\end{equation}
    which is compatible with the automorphisms $S$ and $T$. In particular, this diagram is induced from the diagram \eqref{eq:rel HMS classic} by applying Hochschild homology $HH_a(-)$. 
\end{conj}

\begin{example}\label{ex:Mirror symmetry of CP1}
 Let $(X,D)$ be $(\CP^1, \{0\} \cup \{\infty\})$. The mirror LG model $(Y, \sw)$ is given by Laurent polynomial $\sw(z)=z+\frac{1}{z}$. This is a Lefschetz fibration with two critical values at $\pm1 \in \C$. A generic fiber is two points which gives rise to objects (thimbles) $\Delta_0, \Delta_1$ in Fukaya-Seidel category $FS(Y, \sw)$. Note that by perturbing $\Delta_0$, one can see that it intersects with $\Delta_1$ at one point. In other words $HF^\bullet(\Delta_0, \Delta_1)=\C[z]$ and this implies that both thimbles $\Delta_0$ and $\Delta_1$ correspond to $\cO_X$ and $\cO_X(1)$, respectively.
\end{example}

\begin{example}\label{ex:Mirror Symmetry surface}
Let $X$ be a del Pezzo surface of degree $0 \leq d \leq 9$ and $D$ be a smooth anti-canonical divisor. The compactified mirror LG model is an elliptic fibration $f: Z \to \bbP^1$ over $\bbP^1$ with $3+d$ singularities near the origin and the wheel of $9-d$ lines at infinity. By removing the fiber at infinity $f^{-1}(\infty)$, we have a genuine LG model $(Y,\omega, \sw: Y \to \C)$. In \cite{AKOHMSDelpezzo}, the authors essentially prove the HMS. The mirror symmetry of del Pezzo surfaces with a smooth anti-canonical fiber is recently generalized to log Calabi-Yau surfaces with toric boundaries in \cite{HomologicalMirrorlogcYcluster}. For such a pair $(X,D)$ with $D$ having $N$ irreducible components, a mirror LG potential $\sw:Y \to \C$ is given by a certain Lefschetz fibration whose generic fiber is a $N$-punctured elliptic curve \cite{MirrorlogcYclusterI}. 
\end{example}

\begin{example}
    In \cite{GSveryaffine}, the authors have proved Conjecture \ref{conj:HMS classic} for a pair $(X,D)$ where $X$ is a quasi-projective Fano and $D$ is the toric boundary divisor. The mirror LG model is given by Hori-Vafa construction. The main ingredient is the microlocal sheaf theory \cite{Nadmicrolocal} and its connection with partially wrapped Fukaya categories \cite{GPSMicrolocal}.
\end{example}

When $D$ is smooth and $Y_{sm}$ is proper, Conjecture \ref{conj:CMS classic} can be regarded as a cohomological lift of the mirror P=W conjecture. The commutative diagram (\ref{eq:CMS classic}) can be regarded as an equivalence between a part of the $E_1$-page of the spectral sequence for the weight filtration on $H^\bullet(U)$ (the top row) and that of the perverse Leray filtration associated to $\sw:Y \to \C$ on $H^\bullet(Y)$ (the bottom row). In this case, $T$ is unipotent and the associated monodromy weight filtration $W_\bullet^T$ is centered at $n-|a|$ as well as constitutes the limiting mixed Hodge structure (Proposition \ref{prop: mixed hodge structure f-adapted smooth case}). By taking the associated graded pieces $\Gr^{W^{S/T}}_{2(n-p)}$ for $a \geq 0$, we have
\begin{equation}\label{eq:CMS classic}
    \begin{tikzcd}
         H^{p,q}(D) \arrow[r,"\iota_!"] \arrow[d,"\cong"]&   H^{p,q}(X) \arrow[d,"\cong"] \\
        \Gr^{W^T}_{2(n-p)}H^{n-1+a}(Y_{sm})\arrow[r, "\rho"]&  \Gr^{W^T}_{2(n-p)}H^{n+a}(Y,Y_{sm})
    \end{tikzcd}
\end{equation}
We do the same for $\Gr^{W^{S/T}}_{2(n-q)}$ for $a \leq 0$, and then we get the following isomorphisms:
	\begin{equation*}
	\begin{aligned}
	&\Gr_F^p\Gr^W_{p+q+i}H^{p+q}(U)\cong \Gr^{W}_{2(n-p)}\Gr^P_{n+p-q+i}H^{n+p-q}(Y) \quad \text{for } p-q \geq 0\\
	&\Gr_F^q\Gr^W_{p+q+i}H^{p+q}(U)\cong \Gr^{W}_{2(n-q)}\Gr^P_{n+p-q+i}H^{n+p-q}(Y) \quad \text{for } p-q \leq 0
	\end{aligned}
	\end{equation*}
for $i=0,1$. Here $W_\bullet$ is Deligne's canonical weight filtration. In particular, if the mixed Hodge structure on $H^k(Y)$ is assumed to be of Hodge-Tate type for all $k \geq 0$, that is the same assumption we impose for Theorem \ref{prop: Hodge-Tate implies KKP smooth setting}, then Conjecture \ref{conj:CMS classic} implies the mirror P=W conjecture (Conjecture \ref{conj:mirrorP=W}).

\subsection{Extended Fano/LG Correspondence}\label{sec:extendedfano lg 5.2}
We first propose the definition of symplectic hybrid LG model. Then we propose two formulations of mirror symmetry for a Fano mirror pair \\
$(X,D)|(Y,\omega,h:Y\to \C^N)$: categorical and cohomological.

\begin{defin}
	A \textbf{symplectic hybrid Landau-Ginzburg (LG) model of rank $N$} is a triple 
	$(Y, \omega, h=(h_1, h_2, \dots, h_N):Y \to \C^N)$ where 
	\begin{enumerate}
		\item $(Y, \omega)$ is $n$-dimensional complex K\"ahler Calabi-Yau manifold with a K\"ahler form $\omega \in \Omega^2(Y)$;
		\item $h:Y \to \C^N$ is a proper surjective holomorphic map such that 
		\begin{enumerate}
		    \item (Local trivialization) there exists a constant $R>0$ such that for any non-empty subset $I \subset \{1, \dots, N\}$, the induced map $h_I:Y \to \C^{|I|}$ is a symplectic fibration over the region $B_I:=\{|z_i| > R | i \in I\}$ compatible with $\omega$.  Furthermore,  $j \notin I$, we have $v(h_j)=0$ over $B_I$ for any horizontal vector field $v \in T^{h_I}Y$ (with respect to the two form $\omega$).
		    \item (Compatibility) for $I \subset J$, such local trivializations are compatible under the natural inclusions $B_J \times \C^{N-|J|} \subset B_I \times \C^{N-|I|} \subset \C^N$. 
		\end{enumerate}

	\end{enumerate}
\end{defin}
\noindent A major difference from Definition \ref{def: weak hybrid LG model general} is that we require each induced locally trivial fibrations near infinity $h_I:Y \to \C^{|I|}$ is indeed a symplectic fibration. In other words, over the region $B_I$, a fiber $h_I^{-1}(b_I)$ is a symplectic manifold whose symplectic form is the pullback of $\omega$ under the natural inclusion.

For a symplectic hybrid LG model, one can also promote to the gluing property (Proposition \ref{prop:Gluing Property in general}) in the symplectic category. We will demonstrate how this is used to understand a relative homological mirror symmetry in the rank 2 case (see Section \ref{sec:extended fano lg rank2}). 

We study mirror relations for a Fano mirror pair $(X,D)|(Y, \omega, h:Y \to \C^N)$ (Definition \ref{def:mirrorhybrid}) from two perspectives: categorical and cohomological. A key features is the functoriality among the different mirror relations, which explains how one can see the mirror P=W conjecture from the extended Fano/LG correspondence. 

We introduce some notations. Let $\square$ be a cubical category whose objects are finite subsets of $\N$ and morphisms $Hom(I,J)$ consists of a single element if $I \subset J$ and else is empty. Given a category $\mathsf{C}$, we define a cubical object to be a contravariant functor $F:\square \to \mathsf{C}$, which we also call a cubical diagram of categories. For a cubical object $F$ and $I \subset \N$, we write
\begin{equation*}
    \begin{aligned}
        & F_I:=F(I)  \\
        & d_{IJ}:=F(I \to J): X_J \to X_I, \quad I \subset J
    \end{aligned}
\end{equation*}
We also define a morphism of cubical objects in an obvious way. The ambient category $\mathsf{C}$ we consider is a $(\infty, 1)$ category of $\bbk$-linear differential graded/$A_\infty$ categories where
\begin{itemize}
    \item $\bbk$ is $\C$ or Novikov field;
    \item each object $C$ is smooth and proper \cite{KSAinfty};
    \item for such objects $C_1$ and $C_2$, a morphism $F:C_1 \to C_2$ is spherical \cite{SphericalDG},
\end{itemize}
We denote such a category by $\mathsf{C}_\bbk$. We also consider a category of finite dimensional vector spaces over $\bbk$, denoted by $\mathsf{Vect}_\bbk$. 
\subsubsection{$B$-side objects}
Associated to a Fano pair $(X,D)$, we define the cubical object $\mathfrak{Coh}(X,D)$ in $\mathsf{C}_\C$ by
\begin{equation*}
    \begin{aligned}
        & \mathfrak{Coh}(X,D)_I=D^b\Coh(X_I), \quad \mathfrak{Coh}(X,D)_\emptyset =D^b\Coh(X)\\
        & \mathfrak{Coh}(X,D)_{IJ}=\iota^J_{I*}:D^b\Coh(D_J) \to D^b\Coh(D_I).
    \end{aligned}
\end{equation*}
where $\iota^J_I:D_J \hookrightarrow D_I$ is the natural inclusion and $D^b\Coh(-)$ is the dg enhancement of bounded derived category of coherent sheaves on $(-)$. On $D^b\Coh(X)$, for each $I$, there is an auto-equivalence $\Phi^X_{\Sigma_{i \in I}D_i}$ given by tensoring with the line bundle $\cO_X(\Sigma_{i \in I}D_i)$. For any smooth subvariety $Z$ of $X$ and $\diamond=\Sigma_{i \in I}D_i$, we write $\Phi^Z_\diamond$ for the functor $(-) \otimes \cO_X(\diamond)|_Z$ on $D^b\Coh(Z)$. Then for every $\diamond=\Sigma_{i \in I}D_i$, the cubical object $\mathfrak{Coh}(X,D)$ comes with the auto-equivalence $\Phi_\diamond$ where $\Phi_\diamond(I)=\Phi^{D_I}_\diamond:D^b\Coh(D_I) \to D^b\Coh(D_I)$. Also, the adjunction formula implies that $(\Phi^{D_I}_{\Sigma_{i \notin I}D_i})^{-1}[n-|I|]$ becomes the Serre functor $S_{D_I}$ on $D^b\Coh(D_I)$ for all $I$. 

By applying Hochschild homology $HH_a(-)$ $(-n \leq a \leq n)$ to the cubical object $\mathfrak{Coh}(X,D)$, the HKR isomorpism yields the cubical object $\mathfrak{HH}_a(X,D)$ in $\mathsf{Vect}_\C$,
\begin{equation*}
    \begin{aligned}
        &\mathfrak{HH}_a(X,D)_I=\bigoplus_{p-q=a}H^{p,q}(D_I) \\
        &\mathfrak{HH}_a(X,D)_{IJ}=\iota^J_{I!}:\bigoplus_{p-q=a}H^{p,q}(D_J) \to \bigoplus_{p-q=a}H^{p,q}(D_I)
    \end{aligned}
\end{equation*}
where $\iota_!$'s are the Gysin morphisms. Also, it respects the automorphism, denoted by $S_I$, induced by the auto-equivalence $\Phi_{\Sigma{i \in I}}D_i$ for each $I$. Note that the logarithm of each action is nilpotent because it is given by taking the cup product with the first Chern class $c_1(\Sigma_{i \in I}D_i)$. We write the induced monodromy weight filtration as $W^{S_I}_\bullet$. When $I=\{1, \dots, N\}$, we simply write $S$ for $S_I$. We will take the associated graded pieces $\Gr^{W^S}_{2(n-p)}$ for $a \geq 0$ and $\Gr^{W^S}_{2(n-p)}$ for $a\leq 0$ of the filtration $W^S_\bullet$ to get the Hodge-graded pieces (see Corollary \ref{cor: cms implies pw general}).

\subsubsection{$A$-side objects}
On the A-side, we would expect that there exists a cubical diagram of Fukaya categories mirror to $\mathfrak{Coh}(X,D)$. Let $(Y, \omega, h:Y \to \C^N)$ be a symplectic hybrid LG model. Recall that we have constructed a cubical object $\mathfrak{HH}_a(Y,h)$ in Section \ref{sec:hybridLG},
\begin{equation*}
    \begin{aligned}
        &\mathfrak{HH}_a(Y,h)_I=H^{n+a-|I|}(Y_I, Y_{I,sm}) \\
        &\mathfrak{HH}_a(Y,h)_{IJ}=\rho^J_{I}:H^{n+a-|J|}(Y_J,Y_{J,sm}) \to H^{n+a-|I|}(Y_I,Y_{I,sm})
    \end{aligned}
\end{equation*}
for $-n \leq a \leq n$. Also, there exist monodromy operators $T_I$ given by the composition of $T_i$'s for $i \in I$, where $T_i$ is the monodromy action induced by a large loop along the $i$-th coordinate hyperplane near infinity. On the $A$-side, the main issue is to construct its categorical lift.

\begin{conj}\label{conj:A-side diagram HMS general}
    Let $(Y, \omega, h:Y \to \C^N)$ be a hybrid LG model. There exists a well-defined cubical object $\mathfrak{Fuk}(Y,h)$ in $\mathsf{C}_\bbk$
    \begin{equation*}
    \begin{aligned}
        &\mathfrak{Fuk}(Y,h)_I:=\FS(Y_I, \sw_I) \\
        &\mathfrak{Fuk}(Y,h)_{IJ}:=\cup^J_I:\FS(Y_J, \sw_J) \to \FS(Y_I, \sw_I)
    \end{aligned}
    \end{equation*}
    satisfying the following properties:
    \begin{enumerate}
        \item for each $I$, the monodromy operation $T_I$ gives rise to the auto-equivalence $\Psi_{T_I}$ on $\mathfrak{Fuk}(Y,h)$ such that $(\Psi^{Y_{I^\circ}}_{T_{I}})^{-1}[n-|I|]$ is the Serre functor of $\FS(Y_I, \sw_I)$ where $I^\circ$ is the complement of $I$.
        \item for each $I$, by applying Hochschild homology $HH_a(-)$ for $-n \leq a \leq n$, the cubical diagram $\mathfrak{Fuk}(Y,h)$ with the auto-equivalence $\Psi_{T_I}$ gives rise to the cubical diagram $\mathfrak{HH}_a(Y,h) \otimes \bbk$ with the monodromy operator $T_I$. 
    \end{enumerate}
\end{conj}
\subsubsection{Mirror symmetry}

\begin{conj}\label{conj:HMS general}
    Let $(X,D)|(Y, \omega, h:Y \to \C^N)$ be a Fano mirror pair as above. There is an isomorphism of cubical objects in $\mathsf{C}_\bbk$
    \begin{equation*}
        \mathfrak{Coh}(X,D) \cong D^\pi\mathfrak{Fuk}(Y,h)
    \end{equation*}
    which is compatible with auto-equivalences $\Phi_{\Sigma_{i \in I}D_i}$ and $\Psi_{T_I}$ for each $I$. 
\end{conj}

\begin{example}
    To the author's knowledge, the HMS statement that is closest to the above conjecture is the work of A.Hanlon and J.Hicks \cite{HanlonHicksFunctorialityandhomological} in the case of a smooth projective toric variety $X_\Sigma$. The A-side category is given by the partially wrapped Fukaya category $\mathcal{W}((\C^*)^n,\mathfrak{f}_\Sigma)$ \cite{GPSSectorialdescent}\cite{sylvanpartially} stopped at FLTZ skeleton. In case that $X$ is smooth Fano, this category turns out to be equivalent to the Fukaya-Seidel category of a LG model. They prove the HMS conjecture as a quasi-equivalence
    \begin{equation*}
       D^b\Coh(X) \cong \mathcal{W}((\C^*)^n, \mathfrak{f}_\Sigma) 
    \end{equation*}
     which is compatible with auto-equivalences $\Phi_{\Sigma_{i \in I}D_i}$ and $\Psi_{T_I}$ for all $I$. Functoriality on the A-side is partially established via Lagrangian correspondences on the objects. 
\end{example}

\begin{conj}\label{conj: CMS general}. Let $(X,D)|(Y,\omega, h:Y \to \C^N)$ be a Fano mirror pair. For $-n \leq a \leq n$, there exists an isomorphism of cubical objects in $\mathsf{Vect}_\C$
\begin{equation*}
    \mathfrak{HH}_a(X,D) \cong \mathfrak{HH}_a(Y,h)
\end{equation*}
which is compatible with the automorphisms $S_I$ and $T_I$ for each $I$.
\end{conj}

Conjecture \ref{conj: CMS general} implies that the monodromy actions $T_i$ associated to the hybrid LG model are unipotent for $i=1, \dots, N$. Geometrically, this means that $h_i$ admits a tame compactification $f_i:Z_i \to \bbP^1$ whose pole divisor is reduced. The following corollary implies that Conjecture \ref{conj: CMS general} can be seen as a cohomological lift of the mirror P=W conjecture (Conjecture \ref{conj:mirrorP=W}). 

\begin{cor}\label{cor: cms implies pw general}
    If the mixed Hodge structure on $H^k(Y)$ is of Hodge-Tate type for all $k \geq 0$, then Conjecture \ref{conj: CMS general} implies the mirror P=W conjecture (Conjecture \ref{conj:mirrorP=W}).
\end{cor}

\begin{proof}
    Applying the Mayer-Vietoris sign rule \eqref{eq:MV sign rule}, we obtain the isomorphism between the $E_1$-page of the spectral sequence for the weight filtration on $H^*(X)$ and that for the perverse Leray filtration on $H^*(Y)$. Since both spectral sequences degenerate at the $E_2$-page, we have 
   
\begin{equation*}
	\begin{aligned}
	&\Gr_F^p\Gr^W_{p+q+i}H^{p+q}(U)\cong \Gr^{W}_{2(n-p)}\Gr^P_{n+p-q+i}H^{n+p-q}(Y) \quad \text{for } p-q \geq 0\\
	&\Gr_F^q\Gr^W_{p+q+i}H^{p+q}(U)\cong \Gr^{W}_{2(n-q)}\Gr^P_{n+p-q+i}H^{n+p-q}(Y) \quad \text{for } p-q \leq 0
	\end{aligned}
	\end{equation*}
for $i=0,1,\dots, N$. Now if we assume that the mixed Hodge structure on $H^k(Y)$ is of Hodge-Tate type for all $k \geq 0$, we have the mirror P=W conjecture (Conjecture \ref{conj:mirrorP=W}).
\end{proof}

\begin{rem}
    One may want to work with a general flag $\mathcal{F}$ of subvarieties in $Y$ whose induced flag filtration is the same as the perverse filtration, instead of introducing a cubical object. This approach is more natural from a topological viewpoint. However, the mirror interpretation of each flag in $\mathcal{F}$ is not clear unless the rank 2 case. Once we choose $Y_{sm}$ as the first flag, even though we have the local fibration $h_{Y_{sm}}:Y_{sm} \to \C^{N-1}$, a choice of the second flag is not canonical because we could not compose the summation map with $h_{Y_{sm}}$. Also, the author does not know a proof of its lift to the category of mixed Hodge structures for computing of the weight filtration $W_\bullet$ on $H^*(Y)$.
\end{rem}

\subsection{Rank 2 case}\label{sec:extended fano lg rank2}
We first summarize the previous discussion specializing to the rank 2 case. Let $(X,D)|(Y, \omega, h:Y \to \C^2)$ be a Fano mirror pair of rank 2. Conjecture \ref{conj:HMS general} claims that there is an equivalence between two diagrams of categories:
	\begin{equation*}
	\begin{tikzcd}
    D^b\Coh(D_{12})\arrow{r}{\iota^{12}_{1*}}\arrow{d}{\iota^{12}_{2*}}
    &  D^b\Coh(D_1)
    \arrow{d}{\iota^1_{{X}*}} \\
  D^b\Coh(D_2) \arrow{r}{\iota^2_{{X}*}} 
   & D^b\Coh(X) 
  \end{tikzcd} \qquad
	\begin{tikzcd}
	        D^\pi\Fuk(Y_{12}) \arrow[r, "\cup^{12}_1"] \arrow[d, "\cup^{12}_2"] & D^\pi\FS(Y_1, \sw_1) \arrow[d, "\cup^{1}_{Y}"] \\
	        D^\pi\FS(Y_2, \sw_2) \arrow[r, "\cup^{1}_{Y}"] & D^\pi\FS(Y, \sw)
	    \end{tikzcd}
    \end{equation*}
  which is compatible with the auto-equivalences $\Phi_{D_i}$ and $\Psi_{T_i}$ for each $i=1,2$. By taking Hochschild homology $HH_a(-)$, for $-n<a<n$, this equivalence is expected to produce the equivalence between two diagrams of cohomology groups 
  \begin{equation*}
  \begin{tikzcd}
    \bigoplus_{p-q=a}H^{p,q}(D_{12}) \arrow[r, "\iota^{12}_{1!}"] \arrow[d, "\iota^{12}_{2!}"] & \bigoplus_{p-q=a}H^{p,q}(D_1) \arrow[d, "\iota^1_{X!}"] \\
    \bigoplus_{p-q=a}H^{p,q}(D_2) \arrow[r, "\iota^2_{X!}"] & \bigoplus_{p-q=a}H^{p,q}(X) 
    \end{tikzcd}
    \begin{tikzcd}
    H^{a+n-2}(Y_{12})\arrow[r, "\rho^{12}_1"] \arrow[d, "\rho^{12}_2"] & H^{a+n-1}(Y_1, Y_{12}) \arrow[d, "\rho^1_{Y}"] \\
    H^{a+n-1}(Y_2, Y_{12})\arrow[r, "\rho^2_{Y}"] & H^{a+n}(Y, Y_{sm})
    \end{tikzcd}
    \end{equation*}
    and the auto-equivalences $\Phi_{D_i}$ and $\Psi_{T_i}$ induce the unipotent endomorphism $S_i$ and $T_i$, respectively. By applying the Mayer-Vietoris sign rule \eqref{eq:MV sign rule}, we have the filtered isomorphism of the diagrams
\begin{equation*}
\begin{tikzcd}
\bigoplus_{p-q=a}H^q(D_{12}, \Omega_{D_{12}}^p) \arrow{r}{(\iota^{12}_{1!}, -\iota^{12}_{2!})} \arrow{d}{\cong} & [2ex]
 \bigoplus_{i=1}^2 \bigoplus_{p-q=a}H^q(D_i, \Omega_{D_i}^p) \arrow{r}{\iota^1_{X!}+\iota^2_{X!}} \arrow{d}{\cong} & [4ex] \bigoplus_{p-q=a}H^q(X, \Omega_{X}^p) \arrow{d}{\cong}\\
H^{a+n-2}(Y_{12}) \arrow{r}{(\rho^{12}_1, -\rho^{12}_2)} &  \bigoplus_{i=1}^2 H^{a+n-1}(Y_i, Y_{12}) \arrow{r}{\rho_1 + \rho_2} & H^{a+n}(Y,Y_{sm})
\end{tikzcd}
\end{equation*}
Taking the associated graded pieces $\Gr_{2(n-p)}^{W^{S/T}}$ for $a \geq 0$ and $\Gr_{2(n-q)}^{W^{S/T}}$ for $a \leq 0$, it induces
\begin{equation*}
	\begin{aligned}
	&\Gr_F^p\Gr^W_{p+q+i}H^{p+q}(U)\cong \Gr^{W}_{2(n-p)}\Gr^P_{n+p-q+i}H^{n+p-q}(Y) \quad \text{for } p-q \geq 0\\
	&\Gr_F^q\Gr^W_{p+q+i}H^{p+q}(U)\cong \Gr^{W}_{2(n-q)}\Gr^P_{n+p-q+i}H^{n+p-q}(Y) \quad \text{for } p-q \leq 0
	\end{aligned}
	\end{equation*}
for $i=0,1,2$. If one further assumes that the mixed Hodge structure on $H^k(Y)$ is of Hodge-Tate type for all $k \geq 0$, then we obtain the mirror P=W conjecture. 

While the argument on the $B$-side is concrete, the argument on the $A$-side is speculative. We try to deliver the idea why such an argument would work. Temporarily, we drop the properness condition on a hybrid LG model $(Y,\omega, h:Y \to \C^2)$ and suppose that $(Y,\omega)$ is a Louville manifold and the associated LG models $(Y,\sw), (Y_1, \sw_1)$ and $(Y_2, \sw_2)$ are all Lefschetz fibrations. One way to construct the cup functor from $\FS(Y_i, \sw_i) \to \FS(Y, \sw)$ for $i=1,2$ is to introduce the intermediate Fukaya-Seidel category $\FS(Y_{sm}, h_{Y_{sm}})$. Note that the gluing property implies that $\FS(Y_{sm}, h_{Y_{sm}})$ has the semi-orthogonal decomposition $\langle\FS(Y_1, \sw_1), \FS(Y_2, \sw_2)\rangle$. If we take $\FS(Y,\sw)$ to be partially wrapped Fukaya category \cite{GPSSectorialdescent} associated to a LG model, then the cup functor $\cup^i_{Y}:\FS(Y_i, \sw_i) \to \FS(Y,\sw)$ is given by the compositions 
	\begin{equation}\label{eq:gluing+composition}
		\FS(Y_{i},\sw_{i} ) \xrightarrow{\iota_i} \FS(Y_{sm}, h_{Y_{sm}}) \xrightarrow{\rho} \mathcal{W}(Y_{sm}) \xrightarrow{\cup^{sm}_Y} \FS(Y, \sw)
	\end{equation}
	where $\iota_i$ comes from the semi-orthogonal decomposition of $\FS(Y_{sm}, h_{Y_{sm}})$ and $\rho$ is the localization. By applying Hochschild homology $HH_a(-)$ on the sequence (\ref{eq:gluing+composition}), we obtain the morphism $\rho_Y^i:H^{a+n-1}(Y_1,Y_{12}) \to H^{a+n}(Y,Y_{sm}$ as the composition of two morphisms
	\begin{equation*}
	    H^{a+n-1}(Y_1, Y_{12}) \to H^{a+n-1}(Y_{sm}, Y_{12}) \xrightarrow{\rho^{sm}_Y} H^{a+n}(Y, Y_{sm})
	\end{equation*}
	where the first morphism is the canonical inclusion under the gluing isomorphism  and the second morphism $\rho^{sm}$ is the connecting homomorphism of the long exact sequence of cohomology groups of the triple $(Y_{12}, Y_{sm}, Y)$. Especially for the second one, we refer to \cite[Example 1.4]{GPScovariantly} for the computation of Hochschild homology of the partially wrapped Fukaya category. 
\appendix
\section{Hodge theory of simplicial varieties}\label{sec:Appendix}

\subsection{Review of Deligne's Hodge theory}\label{Appendix MHC}
We recollect some aspects of Deligne's Hodge theory \cite{TheoriedeHodge2}\cite{TheoriedeHodge} that are used in the main part of the paper. We mainly follow the framework in \cite{ElzeinHodge}. We only deal with the coefficients $\Q$ and $\C$ for simplicity. 

\subsubsection{A mixed Hodge complex}
We first introduce the general framework of the $\Q$-mixed Hodge structures, which is called \textbf{$\Q$-mixed Hodge complex} (MHC for short).

\begin{defin}(MHC) A $\Q$-mixed Hodge complex (MHC) $K$ consists of the following data:
\begin{enumerate}
    \item a complex $K_\Q$ of $\Q$-modules whose cohomology $H^k(K_\Q)$ is a $\Q$-module of finite type;
    \item a filtered complex $(K_\Q, W)$ of $\Q$-vector spaces with an increasing filtration $W$;
    \item a bi-filtered complex $(K_\C, W, F)$ of $\C$-vector spaces with an increasing (resp. decreasing) filtration $W$ (resp. $F$) and an isomorphism
    \begin{equation*}
        \alpha:(K_\Q, W) \otimes \C \xrightarrow{\cong} (K_\C, W)
    \end{equation*}
    satisfying the following axioms: for all $n, i$, the cohomology $H^i(-)$ of the triple 
    \begin{equation*}
        (\Gr^W_n(K_\Q), \Gr^W_n(K_\C, F), \Gr^W_n(\alpha))
    \end{equation*}
    form a $\Q$-Hodge structure of weight $n+i$.
\end{enumerate}

We write a $\Q$-MHC as a quadruple $K=(K_\Q, (K_\Q, W), (K_\C, W, F), \alpha)$.

\end{defin}\label{def: cmhc}
\begin{defin}(Cohomological MHC)
A $\Q$-cohomological mixed Hodge complex (cohomological MHC) $K$ on a topological space $X$ consists of the following data:
\begin{enumerate}
    \item a complex $K_\Q$ of sheaves of $\Q$-modules such that $H^k(X, K_\Q)$ are $\Q$-module of finite type;
    \item a filtered complex $(K_\Q, W)$ of sheaves of $\Q$-vector spaces with an increasing filtration $W$;
    \item a bi-filtered complex of sheaves $(K_\C, W, F)$ of $\C$-vector spaces on $X$ with an increasing (resp. decreasing) filtration $W$ (resp. $F$) and an isomorphism 
    \begin{equation*}
        \alpha:(K_\Q, W) \otimes \C \xrightarrow{\cong} (K_\C, W)
    \end{equation*}
    \item the induced quadruple $(R\Gamma(X, K_\Q), (R\Gamma(X, K_\Q),W), (R\Gamma(X, K_\C),W,F), R\Gamma(\alpha))$ forms a $\Q$-mixed Hodge complex.
\end{enumerate}
\end{defin}

By definition, for a given $\Q$-cohomological MHC $K=(K_\Q, (K_\Q, W), (K_\C, W, F), \alpha)$, we have the induced $\Q$-MHC 
\begin{equation*}
    R\Gamma K=(R\Gamma K_\Q, R\Gamma(K_\Q, W), R\Gamma(K_\C, W, F), R\Gamma(\alpha))
\end{equation*}

\begin{thm}(Deligne)
The cohomology of a $\Q$-mixed Hodge complex carries a $\Q$-mixed Hodge structure. 
\end{thm}

\begin{example}
One of the main geometric examples we are interested in is a $\Q$-mixed Hodge structure on the cohomology of smooth quasi-projective variety $U$. One needs to choose a good compactification $(X,D)$ of $U$ and construct a $\Q$-cohomological MHC over $X$ that yields a $\Q$-mixed Hodge structure on $H^\bullet(U)$. This is achieved by considering the logarithmic de Rham complexes $\Omega_Z^\bullet(\log D)$. See Section \ref{sec:weight} for more details.
\end{example}

We introduce one more notion which is frequently used in the main context of the paper.

\begin{defin}
A $\Q$-mixed Hodge complex $K=(K_\Q, (K_\Q, W), (K_\C, W, F), \alpha)$ is Hodge-Tate if the induced $\Q$-mixed Hodge structure on $H^i(\Gr_n^W(K))$ is Hodge-Tate for all $n$.
\end{defin}

\subsubsection{A simplicial variety}

In order to construct a mixed Hodge structure for a singular/non-compact algebraic variety $S$, Deligne considers a resolution of $S$ by smooth varieties that recovers the cohomological data of $S$. This is achieved by employing a simplicial method, as originally presented by Deligne in \cite{TheoriedeHodge}.

Let $\Delta$ be a simplicial category whose objects are the subsets $\Delta_n:=\{0, 1, 2, \dots, n\}$ of integers for $n \in \N$ and morphisms are non-decreasing maps. For $0 \leq i \leq n+1$, we define the $i$-th face map $\delta_i:\Delta_n \to \Delta_{n+1}$ to be the one that does not hit the element $i \in \Delta_{n+1}$. 

\begin{defin}\label{def:simplicial space}
Let $\Top$ be the category of topological spaces.
\begin{enumerate}
    \item A simplicial space is defined to be a contravariant functor $X:\Delta \to \Top$ and denoted by $X=\{X_{[\bullet]}\}$.
    \item A morphism of two simplicial spaces is defined to be a natural transformaion of the functors.
    \item A sheaf $F_{[\bullet]}$ on a simplicial space $X_{[\bullet]}$ is the collection of sheaves $\{F_{[n]}\}$ on each $X_{[n]}$ such that 
    \begin{enumerate}
        \item for each $f:\Delta_n \to \Delta_m$ with $X(f):X_{[m]} \to X_{[n]}$, we have the induced pullback morphism $F(f)_*:X(f)^*:F_{[n]} \to F_{[m]}$;
        \item for each $f:\Delta_n \to \Delta_m$ and $g:\Delta_r \to \Delta_n$, we have $F(f \circ g)_*=F(f)_*\circ F(g)_*$.
    \end{enumerate}
    \item A morphism of sheaves $\phi:F \to G$ on a simplicial space $X$, consists of morphisms $\phi_{[n]}:F_{[n]} \to G_{[n]}$ which are compatible with natural morphisms $F(f)_\bullet$ and $G(f)_\bullet$ for all $f:\Delta_n \to \Delta_m$ .
\end{enumerate}
\end{defin}

As variants of Definition \ref{def:simplicial space}, one can define the notion of complexes of sheaves $\{K_{[q]}^p\}_{p,q}$ on the simplicial space $X$  where $p$ is the degree of the complex and $q$ is the simplicial degree. Morphisms between two such complexes of sheaves can be similarly defined.

\begin{defin}
A simplicial space $X_{[\bullet]}$ is called smooth (resp. compact) if each topological space $X_n$ is smooth (resp. compact) for all $n \geq 0$.
\end{defin}

\begin{example}(Constant simplicial space)\label{ex:const simplicial}
Given a topological space $S$, there is the associated simplicial space $\{S_{[\bullet]}\}$ where all $S_{[n]}$'s are $S$ and $S(f)=Id$ for all $f:\Delta_n \to \Delta_m$. 
\end{example}
\begin{defin}
    Let $\mathcal{V}ar_\mathbb{K}$ be the category of algebraic varieties over $\mathbb{K}=\Q (\text{ or } \C)$. A simplicial space $X_{[\bullet]}$ is called a simplicial variety if it belongs to a simplicial object in $\mathcal{V}ar_\mathbb{K}$
\end{defin}

\begin{example}\label{ex:snc}
Let $D$ be a simple normal crossing variety with irreducible components $D_1, \dots, D_N$. We write $D_{[i]}$ for the disjoint union of $(i+1)$-th intersections of irreducible components. There is a canonical simplicial variety associated to $D$ 
\begin{equation*}
    D_{[\bullet]}:=(D_{[0]} \xleftarrow{(\delta_0, \delta_1)} D_{[1]} \xleftarrow{(\delta_0, \delta_1, \delta_2)} D_{[2]} \leftarrow \cdots \leftarrow D_{[N-1]} \leftarrow 0 \cdots)
\end{equation*}
Note that each $D_{[i]}$ is smooth and there is a canonical morphism of simplicial varieties $a:D_{[\bullet]} \to D$ induced by the inclusions. This is called the simplicial (or Mayer-Vietoris) resolution of the normal crossing variety $D$.
\end{example}

The constant simplicial space, introduced in Example \ref{ex:const simplicial}, allows us to do sheaf theory on a simplicial space. Any sheaf $F$ on the constant simplicial space $S_\bullet$ induces the co-simplicial sheaf $\{F^n\}$ on $S$ with the face maps $\delta_i:F^n \to F^{n+1}$. In particular, if $F$ is abelian, then the co-simplicial sheaf $\{F^n\}$ becomes the complex whose differential is $d=\sum_i (-1)^i\delta_i:F^n \to F^{n+1}$. Similarly, the complex of abelian sheaves $\{K_{[q]}^p\}$ on $S_\bullet$ induces the bicomplex $K^{p,q}$ where $q$ is the co-simplicial degree. One can consider the associated simple complex $sK$
\begin{equation*}
    (sK)^n:=\bigoplus_{p+q=n}K^{p,q}, \quad d(x^{p,q})=d_K(x^{p,q})+\sum_i (-1)^i\delta_ix^{p,q}
\end{equation*}

\noindent Consider a simplicial space $X_{[\bullet]}$ over $S$, where the augmentation map $a:X_{[\bullet]} \to S$ is regarded as a morphism of simplicial spaces. For example, when $S$ is a point, the derived pushforward $Ra_*(-)$ defines the notion of global sections of simplicial sheaves/complexes of simplicial sheaves on $X$. Also, given a sheaf $\cF$ on $S$, the augmentation map $a:X_{[\bullet]} \to S$ gives the  natural morphism $\cF \to Ra_*a^*(\cF)$ on $S$. 
\begin{defin}
    The morphism $a:X_{[\bullet]} \to S$ is of \textit{cohomological descent} if the natural morphism $\phi(a):\cF \to Ra_*a^*\cF$ is an isomorphism for all abelian sheaves $F$ on $S$.
\end{defin}

\begin{thm}\label{thm:simplicial resol}(Deligne)
    Let $U$ be a quasi-projective algebraic variety. There exists an augmented smooth simplicial variety $a:U_{[\bullet]} \to U$ such that 
    \begin{enumerate}
        \item the morphism $a:U_{[\bullet]} \to U$ is cohomological descent;
        \item there exists a simplicial compact smooth variety $X_{[\bullet]}$ containing a simplicial normal crossing divisor $D_{[\bullet]}$ whose complement is $U_{[\bullet]}$. 
    \end{enumerate}
\end{thm}
Theorem \ref{thm:simplicial resol} implies that we are able to construct a $\Q$-mixed Hodge structure on the cohomology of $U$ by taking a simplicial resolution of $U$. 

\subsubsection{Hodge theory on a simplicial space}

One can generalize the notion of a $\Q$-cohomological mixed Hodge complex over a simiplicial space $X_{[\bullet]}$ by considering the same data $(K_\Q, (K_\Q, W), (K_\C, W, F), \alpha)$ over $X_{[\bullet]}$ and requiring that over each $X_{[n]}$, the restriction becomes the $\Q$-cohomological MHC (See Definition \ref{def: cmhc}). We refer to it as a co-simplicial $\Q$-cohomological mixed Hodge complex over $X_{[\bullet]}$. By taking the global section functor on each $X_{[n]}$, we obtain the co-simplicial version of $\Q$-MHC, denoted by 
\begin{equation*}
    R\Gamma_{[\bullet]} K=(R\Gamma_{[\bullet]} K_\Q, R\Gamma_{[\bullet]}(K_\Q, W), R\Gamma_{[\bullet]}(K_\C, W, F), R\Gamma_{[\bullet]}(\alpha))
\end{equation*}

\noindent Note that the underlying complex $K_\Q$ (or $R\Gamma_{[\bullet]}K_\Q$) is bi-graded. It means that one needs to take the associated simple complex to obtain $\Q$-MHC. A key is to define the relevant increasing filtration $sW$ and decreasing filtration $sF$ on $sK$ to get a $\Q$-MHC. 

\begin{defin}
Let $K=(K_\Q, (K_\Q, W), (K_\C, W, F), \alpha)$ be a co-simplicial $\Q$-mixed Hodge complex and $sK$ be the associated simple complex. We define three filtrations on $sK$:
\begin{equation*}
\begin{aligned}
    & L^r(sK)=s(K^p_{[q]})_{q \geq r} \quad F^r(sK)^i:=\bigoplus_{p+q=i}F^rK^p_{[q]} \\
    & \delta(W, L)_n(sK)^i:=\bigoplus_{p+q=i}W_{n+q}K^p_{[q]}
\end{aligned}
\end{equation*}
We call $L^\bullet$, $F^\bullet$, and $\delta(W, L)_\bullet$ a co-simplicial, Hodge, and weight filtration, respectively.
\end{defin}

\begin{thm}\label{thm:filtration L}\cite[Theorem 8.1.15]{TheoriedeHodge}
    Let $K=(K_\Q, (K_\Q, W), (K_\C, W, F), \alpha)$ be a co-simplicial $\Q$-mixed Hodge complex. Then, 
    \begin{enumerate}
        \item the complex $(sK, \delta(W, L), F)$ is a $\Q$-mixed Hodge complex;
        \item the filtration $L$ on $sK$ is a filtration of $\Q$-mixed Hodge complexes;
        \item the filtration $L$ on $H^*(sK)$ is a filtration of $\Q$-mixed Hodge structures.
    \end{enumerate}
\end{thm}

\begin{cor}\label{cor:MHC hodge-Tate}
   Let $K=(K_\Q, (K_\Q, W), (K_\C, W, F), \alpha)$ be a co-simplicial $\Q$-mixed Hodge complex whose restriction $K_{[n]}$ is MHC of Hodge-Tate type for all $n$. Then, the associated $\Q$-mixed Hodge complex $(sK_\Q, (sK_\Q, \delta(W, L)), (sK_\Q, \delta(W, L), F), \alpha)$ is Hodge-Tate.
\end{cor}
\begin{proof}
    Note that the associated graded pieces $\Gr^n_L(sK)$ is the same as the restriction $K_{[n]}$ as $\Q$-mixed Hodge complexes. By assumption, the cohomology $H^*(\Gr^n_L(sK))$ forms a $\Q$-mixed Hodge structure of Hodge-Tate type. Moreover, $H^*(\Gr^n_L(sK))$ is the $E_1$-page of the spectral sequence of the filtration $L$ on $H^*(sK)$. The conclusion follows from (3) of Theorem \ref{thm:filtration L}
\end{proof}

We provide the description of the $E_1$-page of the spectral sequence associated to the weight filtration $\delta(W, L)$. Consider the double complex
\begin{equation}\label{eq:d_1 simplicial}
    \begin{tikzcd}
        H^{q-(n+1)}(\Gr^W_{n+1}K^\bullet_{[n+p+1]}) \arrow[r, "\partial"] & H^{q-n}(\Gr^W_{n}K^\bullet_{[n+p+1]}) \arrow[r, "\partial"] & H^{q-(n-1)}(\Gr^W_{n-1}K^\bullet_{[n+p+1]}) \\
        H^{q-(n+1)}(\Gr^W_{n+1}K^\bullet_{[n+p]}) \arrow[r, "\partial"] \arrow[u, "d"] & H^{q-n}(\Gr^W_{n}K^\bullet_{[n+p]}) \arrow[r, "\partial"] \arrow[u, "d"] &  H^{q-(n-1)}(\Gr^W_{n-1}K^\bullet_{[n+p]}) \arrow[u, "d"]
    \end{tikzcd}
\end{equation}
where $\partial$ is a connecting morphism and $d$ is the alternating sum of the face maps. Then the associated double complex is the $E_1$-page of the spectral sequence for the weight filtration. In other words, we have 
\begin{equation*}
    {}_{\delta(W,L)}E_1^{p,q}(sK_\Q)=\bigoplus_nH^{q-n}(\Gr^W_nK^*_{[n+p]})
\end{equation*}
and $d_1:{}_{\delta(W,L)}E_1^{p,q} \to {}_{\delta(W,L)}E_1^{p+1,q}$ is given by the induced differential. 

\begin{example}(Example \ref{ex:snc} continued)
    Let $D$ be a simple normal crossing variety introduced in Example \ref{ex:snc}. A co-simplicial $\Q$-cohomological MHC on $D_{[\bullet]}$ is given by 
    \begin{equation*}
        (\Q_{D_{[\bullet]}}, (\Q_{D_{[\bullet]}}, \tau), (\Omega^\bullet_{D_{[\bullet]}}, W, F), \alpha:(\Q_{D_{[\bullet]}}, W) \otimes \C \cong (\Omega^\bullet_{D_{[\bullet]}}, W)
    \end{equation*}
    where $\tau$ is the standard filtration. Then the $E_1$-spectral sequence for the weight filtration $\delta(W, L)$ is 
    \begin{equation*}
    \begin{aligned}
        & E_1^{p,q}=\bbH^{p+q}(D, \Gr^W_{-p}(s\Omega^\bullet_{D_{[\bullet]}}) \cong H^q(D_{[p]};\Q) \\
        & d_1=\sum_{i \leq p+1}(-1)^i\delta_i^*:E_1^{p,q} \to E^{p+1,q}
    \end{aligned}
    \end{equation*}
\end{example}

\begin{example}\label{ex:mhc U}
    Let $U$ be a quasi-projective variety. By Theorem \ref{thm:simplicial resol}, we consider a simplicial resolution $a:U:={U_{[\bullet]}} \to U$ with a compactification $(X_{[\bullet]}, D_{[\bullet]})$ of $U_{[\bullet]}$. Let $j_{[\bullet]}:U_{[\bullet]} \hookrightarrow X_{[\bullet]}$ be the open simplicial embedding. Then
    \begin{equation*}
        (Rj_{[\bullet],*}\Q, (Rj_{[\bullet],*}\Q, \tau), (\Omega^\bullet_{[X_{[\bullet]}}(\log D_{[\bullet]}), W, F))
    \end{equation*}
    is a $\Q$-cohomological MHC on $X_{[\bullet]}$. Furthermore, the
    double complex (\ref{eq:d_1 simplicial}) is given by 
    \begin{equation*}
    \begin{tikzcd}
        H^{q-(2n+2)}(D^{n+1}_{[n+p+1]}) \arrow[r, "G"] & H^{q-(2n)}(D^{n}_{[n+p+1]}) \arrow[r, "G"] & H^{q-(2n-2)}(D^{n-1}_{[n+p+1]}) \\
        H^{q-(2n+2)}(D^{n+1}_{[n+p]}) \arrow[r, "G"] \arrow[u, "d"] & H^{q-(2n)}(D^{n}_{[n+p]}) \arrow[r, "G"]\arrow[u, "d"] &  H^{q-(2n-2)}(D^{n-1}_{[n+p]}) \arrow[u, "d"]
    \end{tikzcd}
\end{equation*}
where $D_{[\bullet]}^n$ denotes the disjoint union of intersections of $n$ components of the normal crossing variety $D_{[\bullet]}$ and $G$ is the Gysin map. Therefore, the $E_1$-spectral sequence for the weight filtration $\delta(W,L)$ is given by 
    \begin{equation*}
    \begin{aligned}
        & E_1^{p,q} \cong \bigoplus_n H^{q-2n}(D_{[n+p]}^n;\Q)  \\
        & d_1=G+(-1)^p d:E_1^{p,q} \to E^{p+1,q}
    \end{aligned}
    \end{equation*}
\end{example}

\subsection{The combinatorial perverse filtration}\label{sec:comb perv}
Let $Y$ be a smooth quasi-projective variety with a proper affinization map $h=(h_1, \cdots, h_N):Y \to \C^N$. Note that a general flag which induces the perverse filtration on the cohomology group of $Y$ comes from the ordered collection of hyperplanes $\mathfrak{H}=\{H_1, \cdots, H_N\}$ in $\C^N$. Due to the genericity condition, we can further assume that any flags induced from the collection $\mathfrak{H}=\{H_1, \cdots, H_N\}$ give rises to the same perverse filtration. It allows one to consider a new filtration on the cohomology of $Y$, which in turn becomes equivalent to the perverse filtration associated to $h:Y \to \C^N$. To define such a filtration, we first recall the notations that was used in the main part of the paper. For any index set $I=\{i_1, i_2, \cdots, i_m\} \subset \{1, 2, \cdots, N\}$, we define 
\begin{equation*}
\begin{aligned}
    &H_I:=H_{i_1} \cap \cdots \cap H_{i_m} \qquad Y_I:=h^{-1}(H_I) \\
    &H(I):=\Sigma_{j \not\in I} H_I \cap H_j 
    \qquad Y(I):=h^{-1}(H(I)) 
\end{aligned}
\end{equation*}
    For each $l \in \{1, 2, \cdots, N\}$, we consider two simple normal crossing subvarieties, $H\{l\}:=\sqcup_{|I|=l}H_I$ and $Y\{l\}:=h^{-1}(H\{1\})$. This defines a flag of subvarieties $\mathcal{F}_\mathfrak{H}$ in $Y$ given as follows:
    \begin{equation*}
       \mathcal{F}_\mathfrak{H}: Y\{N\} \subset Y\{N-1\} \subset \cdots \subset Y\{1\} \subset Y\{0\}:=Y
    \end{equation*}
    \begin{defin} For $k \geq 0$,
    the flag filtration on the cohomology $H^k(Y;\Q)$ associated to the flag $\mathcal{F}_\mathfrak{H}$ is given by 
    \begin{equation*}
        L^\mathfrak{H}_{k+i}H^k(Y;\Q):=\ker(H^k(Y;\Q) \xrightarrow{J_i} H^k(Y\{i\};\Q))
    \end{equation*}
    where $J_i$ is the canonical restriction for all $i \geq 0$.
    \end{defin}
    
    We claim that the filtration $L_\bullet^\mathfrak{H}$ is equal to the perverse filtration associated to $h:Y \to \C^N$. 
    
    \begin{prop}\label{prop:combinatorial perverse}
       Let $h:Y \to \C^N$ be a proper affinization morphism and $\mathfrak{H}$ be a collection of hyperplanes as above. Then the flag filtration $L_\bullet^\mathfrak{H}$ on $H^k(Y;\Q)$ is the perverse Leray filtration associated to $h:Y \to \C^N$ for all $k \geq 0$.
    \end{prop}
    \begin{proof}
          Consider the following diagram
    \begin{equation*}
        \begin{tikzcd}
    H^k(Y;\Q) \arrow[r, "J_i"] \arrow[rd, swap, "\Phi_i"] & H^k(Y\{i\};\Q) \arrow[d, "\Psi_i"] \\
    & \bigoplus_{|I|=i}H^k(Y_I;\Q) \\
        \end{tikzcd}
    \end{equation*}
    where $\Phi_i$ is the direct sum of the canonical restrictions and $\Psi_i$ is the canonical restriction to each irreducible component of $Y\{i\}$. By definition, the kernel of $\Phi_i$ is identified with $P_{k+i}H^k(Y;\Q)$. Therefore, it is clear that $L^\mathfrak{H}_{k+i}H^k(Y;\Q) \subset P_{k+i}H^k(Y;\Q)$. The other direction follows from Lemma \ref{lem:combinatorial perverse}.
    \end{proof}

\begin{lem}\label{lem:combinatorial perverse}
        Let $U$ be a smooth quasi-projective variety and $V=\sqcup V_i$ be a normal crossing divisor of $U$ where all $V_i$s are smooth. Assume that the kernel of the canonical restriction $H^k(U;\Q) \to H^k(V_i;\Q)$ is the same for all $i$. Consider the following diagram
        \begin{equation*} 
        \begin{tikzcd}
    H^k(U;\Q) \arrow[r, "J"] \arrow[rd, swap, "\Phi"] & H^k(V;\Q) \arrow[d, "\Psi"] \\
    & \bigoplus_i H^k(V_i;\Q)
        \end{tikzcd}
    \end{equation*}
     where all the morphisms are the canonical restrictions. Then we have $\ker(\Phi) \subset \ker(J)$.
    \end{lem}
    \begin{proof}
        Note that each cohomology group admits the $\Q$-mixed Hodge structures and all the morphisms respect those structures. By taking the associated graded pieces of the weight filtration, $\Gr^W_s(-)$, we have the following commutative diagram
        \begin{equation}\label{eq:diagram cubical}
        \begin{tikzcd}
    \Gr^W_sH^k(U;\Q) \arrow[r, "\Gr^W_s(J)"] \arrow[rd, swap, "\Gr^W_s(\Phi)"] & \Gr^W_sH^k(V;\Q) \arrow[d, "\Gr^W_s(\Psi)"] \\
    & \bigoplus_i \Gr^W_sH^k(V_i;\Q) 
        \end{tikzcd}
    \end{equation}
 We will show that $\ker(\Gr^W_s(\Phi)) \subset \Ker(\Gr^W_s(J))$ for all $s \geq k$, which proves Lemma \ref{lem:combinatorial perverse}. The diagram (\ref{eq:diagram cubical}) can be deduced from the diagram of the $E_1$-page of the  associated spectral sequences. To write down the spectral sequence, we borrow some notations from Example \ref{ex:mhc U}. For a simplicial variety $Y:={Y_{[\bullet]}}$, we take a good compactification of $Y$, that is a pair of simplicial varieties $(X_{Y,[\bullet]}, D_{Y,[\bullet]})$ such that each pair $(X_{Y,[n]}, D_{Y,[n]})$ is a good compactification of $Y_{[n]}$. Since $D_{Y,[n]}$ is the simple normal crossing divisor in $X_{Y,[n]}$, we denote the disjoint union of $m$-th intersections of irreducible components of $D_{Y, [n]}$ by $D_{Y,[n]}^m$. Then the $E_1$-page of the spectral sequence for the weight filtration on $H^*(Y)$ is given by 
    \begin{equation*}
        (E_1^{p,q}(Y)=\bigoplus_nH^{q-2n}(D^n_{Y,[n+p]}), d_1) 
    \end{equation*}
    where the differential $d_1$ is given by the sum of the Gysin maps $G$ and the alternating sum of the face maps:
    \begin{equation}\label{eq:diagram spseq nc}
        \begin{tikzcd}
            E_1(Y)^{p-1,q} & E_1(Y)^{p,q} & E_1(Y)^{p+1,q} \\
         H^{q+2p-2}(D^{-p+1}_{Y,[0]}) \arrow[r, "G"] \arrow[draw=none]{d}{\bigoplus}\arrow[rd, "\sum_i(-1)^i\delta_i"]& H^{q+2p}(D^{-p}_{Y,[0]}) \arrow[r, "G"] \arrow[draw=none]{d}{\bigoplus} \arrow[rd, "\sum_i(-1)^i\delta_i"]& H^{q+2p+2}(D^{-p-1}_{Y,[0]}) \arrow[draw=none]{d}{\bigoplus} \\
            H^{q+2p-4}(D^{-p+2}_{Y,[1]}) \arrow[r, "G"] &
             H^{q+2p-2}(D^{-p+1}_{Y,[1]}) \arrow[r, "G"] & 
              H^{q+2p}(D^{-p}_{Y,[1]})  
        \end{tikzcd}
    \end{equation}
    Note that the spectral sequence converges at the $E_2$-page and we have 
    \begin{equation*}
        E_2^{p,q}(Y)=\Gr^W_qH^{p+q}(Y)
    \end{equation*}
    Coming back to our situation, let $X_U$ be a smooth compactification of $U$ with the simple normal crossing boundary $D_U$. It induces the good compactification of $V$ viewed as a smooth simplicial variety, denoted by $(X_{V, [\bullet]}, D_{V, [\bullet]})$. Consider the diagram of the $E_1$-pages of the spectral sequences associated to the weight filtrations 
    \begin{equation*}
        \begin{tikzcd}
    E_1^{\bullet,\bullet}(U) \arrow[r, "J"] \arrow[rd, "\Phi"] & E_1^{\bullet,\bullet}(V) \arrow[d, "\Psi"] \\
    & \bigoplus_{i=1}^N E_1^{\bullet,\bullet}(V_i) \\
        \end{tikzcd}
    \end{equation*}
    which induces the diagram (\ref{eq:diagram cubical}). Now, take $\alpha \in \ker(\Gr_q^W\Phi:\Gr_q^WH^{p+q}(U;\Q) \to \bigoplus_i\Gr_q^WH^{p+q}(V_i;\Q))$. By assumption, in the $E_1^{p,q}(U)$-term, this is represented by the element in the top part of the diagram \label{eq:diagram spseq nc}$(\alpha, \cdots, \alpha) \in H^{q+2p}(D_{U,[0]}^{-p})=\bigoplus_{i=1}H^{q+2p}(D_U^{-p})$ and the other terms are 0. Then $E_1^{p,q}(J)(\alpha)$ sits in $H^{q+2p}(D_{V,[0]}^{-p})$. On the other hand, in the differential $d_1:E_1^{p-1,q}(V) \to E_1^{p,q}(V)$, the Gysin morphism $G:H^{q+2p-2}(D_{V,[0]}^{-p+1}) \to H^{q+2p}(D_{V, [0]}^{-p})$ is the direct sum of the same Gysin morphism $\bigoplus_{i=1}^N G:\bigoplus_{i=1}^NH^{q+2p-2}(D_{V_i,[0]}^{-p+1}) \to \bigoplus_i^N H^{q+2p}(D_{V_i, [0]}^{-p})$ in the differential $d_1:\bigoplus_{i=1}^NE_1^{p-1,q}(V_i) \to \bigoplus_{i=1}^NE_1^{p,q}(V_i)$. Since $\Gr_q^W(\Phi)(\alpha)=0$, the element $E_1^{p,q}(J)(\alpha)$ is in the image of $G:H^{q+2p-2}(D_{V,[0]}^{-p+1}) \to H^{q+2p}(D_{V, [0]}^{-p})$.

    \end{proof}

\noindent The $E_1$ page of the spectral sequence for the filtration $L^{\mathfrak{H}}_{\bullet}$ is given as follows:
    \begin{equation}\label{eq:spectral combperv}
        H^{k-N}(Y\{N\}) \xrightarrow{d_1} H^{k-N+1}(Y\{N-1\}, Y\{N\}) \xrightarrow{d_1} \cdots \xrightarrow{d_1} H^{k-1}(Y\{1\}, Y\{2\}) \xrightarrow{d_1} H^k(Y, Y\{1\})
    \end{equation}
    where each map $d_1$ is the connecting homomorphism for the corresponding tuple. This sequence is the same as the sequence \ref{eq:E_1-spectral flag} introduced in Section \ref{sec:hybridLG}. By applying the Mayer-Vietoris arguments, the sequence (\ref{eq:spectral combperv}) becomes equivalent to the following Čech type complex
    \begin{equation*}
         H^{k-N}(Y\{N\}) \xrightarrow{d_1} \bigoplus_{|I|=N-1}H^{k-N+1}(Y_I, Y(I)) \xrightarrow{d_1} \cdots \xrightarrow{d_1} \bigoplus_{|I|=1} H^{k-1}(Y_I, Y(I)) \xrightarrow{d_1} H^k(Y, Y\{1\})
    \end{equation*}
    where $d_1$ is the alternating sum of the relevant connecting homomorphisms with the same sign rule in (\ref{eq:MV sign rule}).

    \begin{rem}
        The combinatorial description of the perverse filtration depends on the assumption that the collection of hyperplanes $\mathfrak{H}$ is sufficiently generic to yield Theorem \ref{thm:cataldo}. However, this assumption does not hold in general for the coordinate hyperplanes of a hybrid LG potential. In such cases, the combinatorial description follows from the gluing property (Proposition \ref{prop:Gluing Property in general}). The author acknowledges the anonymous referee for highlighting this aspect.
    \end{rem}

\printbibliography

\end{document}